\documentclass[twoside,11pt]{article}
\usepackage{amsmath, mathtools,mathrsfs}
\usepackage{smile}
\usepackage[linesnumbered,ruled,vlined]{algorithm2e}
\usepackage{tikz}
\usepackage{bbm}
\usetikzlibrary{shapes}
\usetikzlibrary{plotmarks}
\usepackage{xcolor}
\usepackage{capt-of}
\usepackage{fullpage}
\usepackage[colorlinks=true,
            linkcolor=blue,
            urlcolor=blue,
            citecolor=blue]{hyperref}

\usepackage[colorinlistoftodos]{todonotes}

\SetKwInput{KwInput}{Input}               
\SetKwInput{KwOutput}{Output}            
\SetKwRepeat{KwRepeat}{repeat}{until}

\theoremstyle{definition}
\newtheorem*{definition*}{Definition}

\newtheorem*{remark*}{Remark}

\begin{document}
\begin{center}

{\LARGE Minimaxity in Robust Nonparametric Regression: from Heavy-Tailed to Sub-Gaussian Assumptions}
\vspace*{1.5em}

{
    \renewcommand{\thefootnote}{\fnsymbol{footnote}}
    \large
    
    Yikun Li \hspace*{1em} Matey Neykov
}
\vspace*{1.5em}

Department of Statistics and Data Science, Northwestern University
\vspace*{1em}

\texttt{yikunli2028@u.northwestern.edu ~  mneykov@northwestern.edu}

\end{center}

\begin{abstract}
We study nonparametric regression over star-shaped function classes under strong adversarial contamination. The clean (original) observations $(\tilde{X} _{i},\tilde{Y} _{i})_{i=1} ^{n}$ are i.i.d. and satisfy $\tilde{Y} _{i}=f ^{*}(\tilde{X} _{i})+\xi _{i}$, where $f ^{*}\in \mathcal{F}$, $\tilde{X}$ is distributed as $\mathbb{P} _{X}$, and the noise is conditionally centered. An adversary may inspect the clean sample and replace at most $\lfloor \epsilon n\rfloor$ observation pairs. We assume that $\mathcal{F}\subset L _{2}(\mathbb{P} _{X})$ is star-shaped with finite diameter $d$.

Under a Bernstein-type moment condition on function increments and a bounded conditional second moment for the noise, we establish the matching lower and upper bounds of order $\left((\eta ^{*}) ^{2} \vee \epsilon \right)\wedge d ^{2}$, uniformly over admissible star-shaped function classes, therefore characterizing the minimax rate of estimation error, where
\begin{equation*}
    \eta ^{*}:=\sup\left\{\eta >0: n\eta^2\leq\log \mathcal M _{\mathcal{F}}^{\text{loc}}(\eta,c)\right\}.
\end{equation*}
Here $\log \mathcal{M} _{\mathcal{F}} ^{\text{loc}}(\eta,c)$ denotes the local metric entropy of $\mathcal{F}$ in $L _2(\mathbb{P} _{X})$, and $c$ is a sufficiently large constant. Under a uniform relative sub-Gaussian condition on function increments and noise, the lower bound is adjusted to $\left((\eta ^{*})^{2}\vee \epsilon ^{2}\right)\wedge d ^{2}$ while the upper bound improves to $\left((\eta ^{*})^{2}\vee \epsilon ^{2}\log ^{2}\left(\frac{e}{\epsilon }\right)\right)\wedge d ^{2}$, yielding minimax optimality up to logarithmic factors.
\end{abstract}

\tableofcontents

\newpage

\section{Introduction}\label{section: introduction}

\subsection{Problem Formulation}\label{subsection: problem formulation}

Let $\mathcal{X}\subseteq\mathbb{R}^{p}$ be a measurable set, and let $\tilde{X}_{1},\ldots,\tilde{X}_{n}$ be independent and identically distributed according to a probability distribution $\mathbb{P}_{X}$ on $\mathcal{X}$. The uncontaminated responses satisfy
\begin{equation}\label{eq: model of the original data}
    \tilde{Y}_{i}
    = f^{*}(\tilde{X}_{i})+\tilde{\xi }_{i},
    \qquad i=1,\ldots,n,
\end{equation}
where $f^{*}\in\mathcal{F}$ and $\mathcal{F}$ is a given star-shaped class of measurable functions from $\mathcal{X}$ to $\mathbb{R}$ (see Definition~\ref{def_star_shaped_set}). We assume that $\mathcal{F}\subseteq L_{2}(\mathbb{P}_{X})$ and that $\mathcal{F}$ has finite $L_{2}(\mathbb{P}_{X})$-diameter (note that this does not necessarily imply that $\mathcal{F}$ is bounded in sup-norm). The noise-covariate pairs $(\tilde{\xi _{i}},\tilde{X}_{i}), i=1,2,\dots,n$ are independent, identically distributed, and have a common joint distribution $\mathbb{P}_{\xi ,X}$. The noise $\tilde{\xi }$ is allowed to depend on the covariate $\tilde{X}$, but it always satisfies that $\mathbb{E}\left[\tilde{\xi }\left\lvert\right. \tilde{X}\right]=0$. We denote the family of such distributions as $\mathcal{P}_{\xi ,X}$ and therefore $\mathbb{P}_{\xi ,X}\in \mathcal{P}_{\xi ,X}$. Such a setting includes the scenario that $\tilde{\xi }$ is independent of $\tilde{X}$ and $\mathbb{E}\tilde{\xi }=0$. In the following argument, we assume $\mathbb{P}_{X}$ is an arbitrary fixed distribution and write $\mathbb{P}_{\xi ,X}$ and $\mathcal{P}_{\xi ,X}$ as $\mathbb{P}_{\xi }$ and $\mathcal{P}_{\xi }$ for simplicity. The moment and tail conditions defining $\mathcal{P}_{\xi}$ will be specified for each result.

Write $\tilde{D}_{n}:=\left((\tilde{X}_{i},\tilde{Y}_{i})\right)_{i=1}^{n}$ for the uncontaminated sample. Before the sample is observed, an adaptive adversary may replace at most $\lfloor\epsilon n\rfloor$ observations by arbitrary elements of $\mathcal{X}\times\mathbb{R}\times \mathbb{R}$, where $\epsilon\in[0,1/2)$ is a known upper bound on the contamination proportion. The adversary may use the entire uncontaminated sample, the distributions $\mathbb{P}_{X}$ and $\mathbb{P}_{\xi}$, and the regression function $f^{*}$ when selecting the observations to replace and their replacement values.

More precisely, let $\mathfrak{C}_{n,\epsilon}$ denote the collection of admissible measurable, possibly randomized, contamination maps. For each $\mathcal{C}\in\mathfrak{C}_{n,\epsilon}$, the observed sample $D_{n}:=\left((X_{i},Y_{i})\right)_{i=1}^{n}=\mathcal{C}(\tilde{D}_{n})$ satisfies
\begin{equation*}
    \sum\limits _{i=1}^{n}\mathbf{1}_{\left\{(X _{i},Y _{i})\neq (\tilde{X}_{i},\tilde{Y}_{i})\right\}}\le \lfloor \epsilon n \rfloor 
\end{equation*}
almost surely. Thus, all covariates, responses, and the noise may be corrupted. This contamination mechanism is known as the \textit{strong contamination model} \cite{diakonikolas2023algorithmic,oliveira2025}.

For $f,g\in L_{2}(\mathbb{P}_{X})$, define
\begin{equation}\label{eq: l_2 p_x distance}
    \lVert f-g\rVert_{L_{2}(\mathbb{P}_{X})}:=\left(\int_{\mathcal{X}}\left[f(x)-g(x)\right]^{2}d\mathbb{P}_{X}(x)\right)^{1/2}.
\end{equation}
The corresponding norm is
$\lVert f\rVert_{L_{2}(\mathbb{P}_{X})}
:=\lVert f-0\rVert_{L_{2}(\mathbb{P}_{X})}$.
The expression in \eqref{eq: l_2 p_x distance} is a pseudometric on functions and a metric modulo $\mathbb{P}_{X}$-almost-sure equality. Throughout the paper, estimation accuracy is measured by the squared $L_{2}(\mathbb{P}_{X})$ loss.

\begin{remark*}[]
    We do not assume the full knowledge of $\mathbb{P}_{X}$ to exactly evaluate $L _{2}(\mathbb{P}_{X})$-norm of $f \in \mathcal{F}$, as such exact evaluation is unnecessary. For the algorithms and results in this paper to hold, we only require an oracle to access the value of an equivalent norm to the $L _{2}(\mathbb{P}_{X})$-norm. Here by equivalent norm we mean a norm $\|\cdot\|$, which agrees with $\|\cdot\|_{L_2(\mathbb{P}_X)}$ up to universal constants that are independent of $n,\mathcal{F},\epsilon $. A common example is the ordinary linear regression where $\tilde{Y}_{i}=\tilde{X}_{i}^{\top } \beta +\tilde{\xi }_{i}$, $\tilde{X}_{i}$ are independently, identically distributed with $\mathbb{E}\tilde{X}_{i}=0$ and $\text{Cov}\left(\tilde{X} _{i}\right)=\bm{\Sigma }$. In this scenario, having an oracle to $L _{2}(\mathbb{P}_{X})$-norm is equivalent to knowing the covariance $\bSigma$. Under the common assumption that we have nonzero lower and upper bounds $\kappa ,\kappa ^{-1}$ such that $\kappa \le \sigma _{\min}(\bm{\Sigma })\le \sigma _{\max}(\bm{\Sigma })\le \kappa ^{-1}$ the $L_2(\mathbb{P}_X)$ norm is equivalent to the standard Euclidean norm. The full knowledge of $\bm{\Sigma }$ is hence unnecessary. Nonetheless, for notational simplicity, we assume for ease of exposition that $\mathbb{P}_{X}$ is known in the presentations of our results and proofs. They can be extended to the aforementioned scenario with minimal efforts. 
\end{remark*}

Given the contaminated sample $D_{n}$, we seek to estimate $f^{*}$ using a measurable estimator $\hat{f}$ taking values in $\mathcal{F}$. For a specified noise family $\mathcal{P}_{\xi}$, define the minimax risk by
\begin{equation}\label{eq: minimax risk}
    \begin{aligned}
        \mathcal{R}_{n,\epsilon}(\mathcal{F};\mathcal{P}_{\xi}):=\inf_{\hat{f}}\,\sup_{f^{*}\in\mathcal{F}}\,\sup_{\mathbb{P}_{\xi}\in\mathcal{P}_{\xi}}\,\sup_{\mathcal{C}\in\mathfrak{C}_{n,\epsilon}}\mathbb{E}\left\lVert\hat{f}\left(\mathcal{C}(\tilde{D}_{n})\right)-f^{*}\right\rVert_{L_{2}(\mathbb{P}_{X})}^{2}.
    \end{aligned}
\end{equation}
The expectation is taken over the uncontaminated sample and any randomization used by the estimator or the adversary. The infimum ranges over estimators that may depend on the known quantities $\mathcal{F}$, $\mathbb{P}_{X}$, $\epsilon$, and $\mathcal{P}_{\xi}$, but not on the unknown regression function, noise distribution, or contamination mechanism. Dependence on the fixed design distribution $\mathbb{P}_{X}$ is suppressed in the notation.

Structural constraints on a regression function provide a way to incorporate prior information into estimation. Sparsity, smoothness, monotonicity, and bounds on a reproducing kernel Hilbert space norm are examples of such constraints. The corresponding function classes are often star-shaped, including both convex classes and nonconvex classes arising from sparsity restrictions. At the same time, heavy-tailed noise and corrupted observations motivate the development of estimators whose guarantees remain valid beyond an uncontaminated model.

Sparse linear regression provides a motivating example of a nonconvex star-shaped class. For a sparsity level $1\leq s<p$ and a radius $R>0$, consider
\begin{equation*}
    \mathcal{F}_{s,R}:=\left\{f_{\beta}:x\mapsto x ^{\top } \beta:\beta\in\mathbb{R}^{p},\ \left\lVert\beta\right\rVert_{0}\leq s,\ \left\lVert\beta\right\rVert_{2}\leq R\right\}.
\end{equation*}
This class is star-shaped about zero because multiplication by a scalar in $\left[0,1\right]$ preserves both constraints, although the class is not convex. Under a standard sub-Gaussian design ($\operatorname{Cov}\left(x\right)=\mathbf{I}_{p}$), $\left\lVert f_{\beta}-f_{\gamma}\right\rVert_{L_{2}\left(\mathbb{P}_{X}\right)}=\left\lVert\beta-\gamma\right\rVert_{2}$, so the coefficient constraint also ensures a finite $L_{2}\left(\mathbb{P}_{X}\right)$-diameter. \cite{liu2020} study sparse linear regression under Gaussian designs when an adversary may arbitrarily corrupt a fraction of the covariate--response pairs. This example motivates allowing nonconvex structural constraints together with corruption of entire observations.

Smooth nonparametric regression provides another example. H\"older norm balls on a compact domain are convex and contain the zero function, and are therefore star-shaped. Such classes describe response surfaces whose variation is controlled by a smoothness constraint. \cite{du2018} study robust nonparametric regression over H\"older and Sobolev classes under Huber's contamination model with a fixed-grid design. Their results illustrate the interaction between smoothness and contamination in determining estimation accuracy.

Shape-constrained regression provides a complementary setting in which prior information is expressed through qualitative restrictions. For example, the class of coordinatewise nondecreasing functions from $\left[0,1\right]^{p}$ to $\left[-B,B\right]$, where $B>0$, is convex and contains zero, hence is star-shaped. \cite{han2019} establish risk bounds for multivariate isotonic regression under both fixed and random designs in the uncontaminated setting. Studying such classes under contamination allows one to investigate how monotonicity constraints interact with the presence of unreliable observations.

Kernel regression supplies a further family of examples. Let $\mathcal{H}$ be a reproducing kernel Hilbert space with kernel $k$. The norm ball $\mathcal{F}_{\mathcal{H},R}:=\left\{f\in\mathcal{H}:\left\lVert f\right\rVert_{\mathcal{H}}\leq R\right\}$ is convex and star-shaped about zero. If $\sup_{x\in\mathcal{X}}k\left(x,x\right)\leq\kappa^{2}$, the reproducing property gives $\sup_{x\in\mathcal{X}}\left|f\left(x\right)\right|\leq\kappa R$ for every $f\in\mathcal{F}_{\mathcal{H},R}$. \cite{liuLian2026} study conditional mean estimation in reproducing kernel Hilbert spaces with heavy-tailed errors and arbitrary contamination of responses, demonstrating the relevance of robustness to both distributional tails and corrupted data in nonlinear regression.

These examples motivate a unified study of regression over star-shaped function classes under explicit distributional assumptions. We consider a design distribution and a strong contamination mechanism under which an adversary may inspect the entire uncontaminated sample before replacing a prescribed fraction of observations. Our objective is to understand how the geometry of the function class, the noise assumptions, and the contamination proportion determine achievable squared $L_{2}\left(\mathbb{P}_{X}\right)$ risk bounds, and to establish their sharpness for each admissible star-shaped classes.

\subsection{Related Literature}\label{subsection: related literature}

Our work connects the literature on local-entropy characterizations of minimax risk with robust regression under heavy-tailed observations and adversarial contamination.

\paragraph{Local entropy and constrained estimation.}

The information-theoretic bounds in \cite{yang1999} relate nonparametric minimax rates to metric entropy under suitable regularity conditions, with applications including regression. These results provide a classical foundation for using the complexity of a function class to characterize estimation difficulty. For convex sub-Gaussian function classes, \cite{mendelson2017} develops regression bounds governed by local covering parameters under appropriate assumptions on the noise. This work demonstrates the importance of local geometry in determining attainable prediction accuracy and motivates our use of local packing entropy. The minimax risk in the Gaussian sequence model under bounded convex constraints is characterized in \cite{neykov2023} through a critical radius defined by local metric entropy. Its multistage packing construction provides a geometric precursor to subsequent estimation procedures for constrained nonparametric regression. Under Gaussian noise and suitable boundedness or increment conditions, \cite{10.1214/25-EJS2419} characterize the minimax squared $L _{2}(\mathbb{P}_{X})$-risk over bounded star-shaped function classes through local entropy. This result provides the uncontaminated benchmark for the present problem, which additionally considers finite-variance noise and adaptive replacement of sample pairs.

\paragraph{Regression under heavy tails and contamination.}

The decision-theoretic framework in \cite{chen2016} develops robust testing procedures and minimax lower bounds based on a total-variation modulus of continuity under Huber's contamination model. This framework explains why the contamination contribution depends on the statistical model and the loss, a consideration that also guides the formulation of lower bounds in our setting. The median-of-means tournaments in \cite{lugosi2020tournaments} attain optimal accuracy--confidence tradeoffs for regression from independent, uncontaminated observations under appropriate moment and norm-equivalence assumptions. Their use of robust pairwise comparisons provides a methodological foundation for estimation with heavy-tailed observations. The median-of-means procedures in \cite{lecue2020} provide regression and regularization guarantees under weak moment assumptions in the presence of arbitrary outliers. Their bounds involve localized quadratic and multiplier complexities together with an outlier-dependent term, providing a useful comparison for our entropy-based analysis. Specializing the moment-dependent regression bound in \cite{oliveira2025} yields squared-error contamination terms of order $\epsilon$ under bounded second multiplier moments and $\epsilon ^{2}\log ^{2}\left(e/\epsilon\right)$ under sub-exponential multiplier tails. Its quadratic-risk analysis treats convex classes satisfying a small-ball condition and uses localized empirical-process complexities, making the local-entropy formulation for star-shaped classes with known $\mathbb{P}_{X}$ a central point of comparison with our work.

\paragraph{Robust estimation under star-shaped constraints.}

For robust mean estimation under star-shaped constraints, \cite{10.1214/25-AOS2576} combine local entropy with a squared-error contamination term of order $\epsilon ^{2}$ for Gaussian noise and $\epsilon ^{2}\log\left(1/\epsilon\right)$ for general unknown sub-Gaussian noise, with the noise scale held fixed. These results provide a closely related geometric characterization in a vector location model and motivate the corresponding investigation for random-design regression. Under a bounded covariance assumption and an additional sample-size condition involving the supremum of local entropy, \cite{peng2026} characterize robust mean estimation over star-shaped constraints with a squared-error contamination term linear in $\epsilon$. This heavy-tailed extension provides a useful comparison for how moment assumptions affect the contamination penalty, although its vector location model differs from the regression model considered here.

\paragraph{Structured regression models.}

Under fixed-design regression with Huber contamination, the local-median procedures in \cite{du2018} achieve minimax rates over H\"older and Sobolev classes. These smoothness classes provide concrete nonparametric examples motivating a general theory for star-shaped regression classes. Regression-depth estimators in \cite{gao2020} attain minimax rates under Huber contamination for several regression models, including sparse linear and nonparametric regression. These results provide model-specific robust benchmarks for our study of general star-shaped classes under strong contamination. Computationally efficient procedures for sparse linear regression under adversarial corruption of both covariates and responses are developed in \cite{liu2020}, with guarantees for Gaussian designs. Since bounded sparse parameter classes are star-shaped and generally nonconvex, this setting illustrates the relevance of allowing star-shaped constraints. Covariate filtering in \cite{pensia2025} yields computationally efficient robust linear-regression procedures under heavy-tailed covariates, heavy-tailed errors, and adversarial contamination. This work complements our geometric perspective by emphasizing tractable estimation within a specific regression family. The kernel-based procedures in \cite{liuLian2026} combine adaptive Huber loss with reproducing kernel Hilbert space methods to obtain estimation guarantees under heavy-tailed noise and response contamination. Kernel norm balls are convex and hence star-shaped, providing another nonparametric motivation for the class geometry considered here. The present work studies regression over bounded star-shaped classes under strong contamination and a known design distribution $\mathbb{P}_{X}$. Our focus is on local-entropy risk bounds under Bernstein-type increment conditions and finite-variance noise, together with sharper contamination dependence under relative sub-Gaussian increment assumptions and sub-Gaussian noise. In the heavy-tailed regime, the upper bound applies to every admissible class, while sharpness of the contamination term is understood for each star-shaped classes satisfying the assumptions.

\subsection{Our Contributions}\label{subsection: our contributions}

Building on the local-entropy characterization of uncontaminated regression in \cite{10.1214/25-EJS2419}, we investigate how strong contamination and distributional assumptions affect estimation over bounded star-shaped function classes. Applying the generalized Fano inequality to an admissible Gaussian noise submodel yields a lower bound for the uncontaminated problem, which remains valid under strong contamination because the adversary may leave the sample unchanged. This lower bound serves as a benchmark for the following contributions, with implicit constants depending on the moment and tail parameters specified in our assumptions.

\begin{itemize}
    \item \textbf{Robust estimation under finite-variance noise.} Under a Bernstein-type moment condition on the increments of $\mathcal{F}$ and a bounded conditional second moment assumption on $\xi$, we construct the tree-based procedures in Algorithms \ref{algorithm: construction of the pruned tree} and \ref{algorithm: theoretical algorithm} and establish the upper bound
    \begin{equation*}
        \mathcal{R}_{n,\epsilon}(\mathcal{F};\mathcal{P}_{\xi})\lesssim \left(\left(\eta ^{*}\right)^{2}\vee \epsilon\right)\wedge d ^{2}.
    \end{equation*}
    This bound holds for every star-shaped function class satisfying the stated assumptions. We further prove that the linear dependence on $\epsilon$ cannot be improved for each star-shaped class $\mathcal{F}$ uniformly over the admissible noise, therefore establishing the optimality of the contamination term.

    \item \textbf{Sharper bounds under sub-Gaussian assumptions.} Under relative sub-Gaussian increment and boundedness conditions on $\mathcal{F}$ and a sub-Gaussian assumption on $\xi$, we modify the procedures above to obtain the sharper upper bound
    \begin{equation*}
        \mathcal{R}_{n,\epsilon}(\mathcal{F};\mathcal{P}_{\xi})\lesssim \left(\left(\eta ^{*}\right)^{2}\vee \epsilon ^{2}\log ^{2}\left(e/\epsilon\right)\right)\wedge d ^{2}.
    \end{equation*}
    We complement this result with the lower bound
    \begin{equation*}
        \mathcal{R}_{n,\epsilon}(\mathcal{F};\mathcal{P}_{\xi})\gtrsim \left(\left(\eta ^{*}\right)^{2}\vee \epsilon ^{2}\right)\wedge d ^{2},
    \end{equation*}
    thereby characterizing the minimax risk up to a logarithmic factor in the contamination term.
\end{itemize}

\subsection{Organization}\label{subsection: organization}

This paper is structured as follows. Section \ref{section: general lower bound} establishes a general minimax lower bound in terms of the local metric entropy of the function class. Section \ref{section: minimax rate under heavy-tailed assumptions} studies function classes satisfying a Bernstein-type condition with noise having bounded conditional second moments. We construct a robust tree-based estimator and establish matching minimax upper and lower bounds under strong adversarial corruption. Section \ref{section: minimax rate under sub-gaussian assumptions} considers sub-Gaussian function increments and noise, obtaining a sharper dependence on the corruption proportion with upper and lower bounds that agree up to a logarithmic factor. Section \ref{section: examples} illustrates our results through sparse linear regression, low-rank trace regression, and H\"older smooth regression, including explicit calculations of their local entropy and resulting minimax rates. Section \ref{section: discussion} summarizes the paper and introduces some possible future research avenues beyond this work. The proofs of the lemmas, theorems, and corollaries in Sections \ref{section: general lower bound}, \ref{section: minimax rate under heavy-tailed assumptions}, and \ref{section: minimax rate under sub-gaussian assumptions} are deferred to Appendices \ref{section: proof of lower bound of minimax rate}, \ref{section: proof of minimax rate under heavy-tailed assumptions}, and \ref{section: proof of minimax rate under sub-gaussian assumptions}, respectively.

\subsection{Notation and Definitions}\label{subsection: notations and definitions}

We introduce the notation and definitions used throughout this paper. Let $n$ denote the sample size and $p$ the dimension of the covariates $\tilde{X}_{i}$. The parameter $\epsilon$ denotes the maximum fraction of observations that the adversary $\mathcal{C}$ may replace. We define $\tilde{X}:=(\tilde{X}_{1},\dots,\tilde{X}_{n})$, $\tilde{Y}:=(\tilde{Y}_{1},\dots,\tilde{Y}_{n})$, $\tilde{\xi }=(\tilde{\xi }_{1},\dots,\tilde{\xi }_{n})$ and $(X,Y)=\mathcal{C}(\tilde{X},\tilde{Y})$. For a positive integer $n$, write $[n]:=\left\{1,2,\dots,n\right\}$. For real numbers $a$ and $b$, define $a\wedge b:=\min\left\{a,b\right\}$ and $a\vee b:=\max\left\{a,b\right\}$. For nonnegative quantities $a$ and $b$, the notation $a\lesssim_{c,\dots}b$ means that $a\leq Kb$ for some constant $K>0$ depending only on the quantities indicated in the subscript. We write $a\gtrsim_{c,\dots}b$ if $b\lesssim_{c,\dots}a$.

For $f,g\in\mathcal{F}$, let $\left\lVert f-g\right\rVert_{L_{2}\left(\mathbb{P}_{X}\right)}$ denote their $L_{2}\left(\mathbb{P}_{X}\right)$-distance. Let $\mathbb{P}_{n}$ denote the empirical distribution of the observed covariates $X_{1},\dots,X_{n}$. The corresponding empirical $L_{2}$-distance is defined by
\begin{equation*}
    \left\lVert f-g\right\rVert_{L_{2}\left(\mathbb{P}_{n}\right)}:=\left(\frac{1}{n}\sum_{i=1}^{n}\left[f\left(X_{i}\right)-g\left(X_{i}\right)\right]^{2}\right)^{1/2}.
\end{equation*}
For $f\in L_{2}\left(\mathbb{P}_{X}\right)$ and $\delta\geq0$, let $B_{2}\left(f,\delta\right)$ denote the closed $L_{2}\left(\mathbb{P}_{X}\right)$-ball centered at $f$ with radius $\delta$:
\begin{equation*}
    B_{2}\left(f,\delta\right):=\left\{g\in L_{2}\left(\mathbb{P}_{X}\right):\left\lVert f-g\right\rVert_{L_{2}\left(\mathbb{P}_{X}\right)}\leq\delta\right\}.
\end{equation*}
We denote the $L_{2}\left(\mathbb{P}_{X}\right)$-diameter of $\mathcal{F}$ by
\begin{equation*}
    d:=\operatorname{diam}\left(\mathcal{F}\right):=\sup_{f,g\in\mathcal{F}}\left\lVert f-g\right\rVert_{L_{2}\left(\mathbb{P}_{X}\right)}.
\end{equation*}
For two probability distributions $\mathbb{P}$ and $\mathbb{Q}$, we use $\operatorname{KL}(\mathbb{P}\left\lVert\right. \mathbb{Q}):=\mathbb{E}_{\mathbb{P}}\log \left(\frac{d \mathbb{P}}{d \mathbb{Q}}\right)$ to represent the Kullback-Leibler (KL) divergence and $\operatorname{TV}(\mathbb{P},\mathbb{Q}):=\sup\limits _{A \in \mathscr{F}}\left|\mathbb{P}(A)-\mathbb{Q}(A)\right|$ the total variation distance between $\mathbb{P} $ and $\mathbb{Q}$, where $\mathscr{F}$ is the corresponding sigma-algebra.    

The following definitions of star-shapedness, packing numbers, and local entropy are fundamental throughout the analysis.

\begin{definition}[Star-shaped set]\label{def_star_shaped_set}
A function class $\mathcal{F}$ is star-shaped with respect to $f_{c}\in\mathcal{F}$ if the line segment joining $f_{c}$ to every $f\in\mathcal{F}$ is contained in $\mathcal{F}$; that is,
\begin{equation*}
    \left(1-t\right)f_{c}+tf\in\mathcal{F},\qquad \forall f\in\mathcal{F},\ \forall t\in\left[0,1\right].
\end{equation*}
Any such $f_{c}$ is called a star center of $\mathcal{F}$. The class $\mathcal{F}$ is called star-shaped if it has at least one star center.
\end{definition}

\begin{definition}[$\delta$-packing and packing number]\label{def_packing}
Let $\delta>0$. A finite subset $\left\{f_{1},\dots,f_{M}\right\}\subseteq\mathcal{F}$ is a $\delta$-packing of $\mathcal{F}$ with respect to a norm or seminorm $\left\lVert\cdot\right\rVert$ if
\begin{equation*}
    \left\lVert f_{i}-f_{j}\right\rVert\geq\delta,\qquad \forall i,j\in[M],\ i\neq j.
\end{equation*}
The packing number, denoted by $\mathcal{M}\left(\left\lVert\cdot\right\rVert,\delta,\mathcal{F}\right)$, is the supremum of the cardinalities of all finite $\delta$-packings of $\mathcal{F}$, with value $+\infty$ if these cardinalities are unbounded. When the norm or seminorm is clear from the context, we abbreviate this notation to $\mathcal{M}\left(\delta,\mathcal{F}\right)$.
\end{definition}

\begin{definition}[Local packing number and local entropy]\label{def_local_entropy}
For $\eta>0$ and a fixed universal constant $c>0$, define the local packing number of $\mathcal{F}$ by
\begin{equation}\label{eq: definition of local entropy}
    \mathcal{M}_{\mathcal{F}}^{\mathrm{loc}}\left(\eta,c\right):=\sup_{f\in\mathcal{F}}\mathcal{M}\left(\eta/c,B_{2}\left(f,\eta\right)\cap\mathcal{F}\right),
\end{equation}
where the packing numbers are computed with respect to the $L_{2}\left(\mathbb{P}_{X}\right)$-distance. The quantity $\log\mathcal{M}_{\mathcal{F}}^{\mathrm{loc}}\left(\eta,c\right)$ is called the local metric entropy of $\mathcal{F}$. Throughout the paper, $c$ is reserved for the constant in this definition and is taken to be sufficiently large.
\end{definition}

An important property of the local packing number is that for any star-shaped function class $\mathcal{F}$, $\mathcal{M}_{\mathcal{F}}^{\text{loc}}(\eta ,c)$ is nonincreasing with respect to $\eta $, as verified by the following lemma.

\begin{lemma}[Nonincreasingness of $\mathcal{M}(\eta ,c)$, Lemma 1.4 from \cite{10.1214/25-AOS2576}]\label{lemma: nonincreasingness of the local entropy}
    For any star-shaped class $\mathcal{F}$ and constant $c>0$, $\mathcal{M}_{\mathcal{F}}^{\text{loc}}(\eta ,c)$ is nonincreasing with respect to $\eta $.
\end{lemma}

Lemma \ref{lemma: nonincreasingness of the local entropy} guarantees the well-definedness of the critical radius, defined as follows.

\begin{definition}[Critical radius]\label{def_critical_radius}
Given the star-shaped $\mathcal{F}$ and $n$, the critical radius of $\mathcal{F}$ is defined as
    \begin{equation}\label{eq: definition of critical radius}
    \eta^{*}\left(\mathcal{F}\right):=\sup\left\{\eta>0:n\eta^{2}\leq\log\mathcal{M}_{\mathcal{F}}^{\mathrm{loc}}\left(\eta,c\right)\right\},
    \end{equation}
with the convention that $\eta^{*}\left(\mathcal{F}\right)=0$ if the set is empty. We suppress the dependence on $\mathcal{F}$ when it is clear from the context.
\end{definition}

\section{General Lower Bound}\label{section: general lower bound}



In this section, we recall Fano's method to derive the general minimax lower bound of \eqref{eq: minimax risk}. The resulting lower bound for the uncontaminated problem serves as a benchmark for our analysis and remains valid under strong contamination because the adversary may leave the sample unchanged. Fano's method has been widely applied to statistical problems, including constrained mean estimation \cite{10.1214/25-AOS2576,neykov2023}, support recovery in sparse linear regression \cite{Wainwright2007InformationTheoreticLO,4595378}, and support recovery in sparse principal component analysis \cite{10.1214/08-AOS664}.

We state below the version of the generalized Fano inequality used in this paper. We omit the proof and refer the reader to Sections 15.3 and 15.4 of \cite{Wainwright_2019} for a detailed treatment and the necessary information-theoretic background.

\begin{theorem}[Generalized Fano inequality]\label{theorem: fano inequality}
    Let $\mathcal{F}\subset L _{2}(\mathbb{P}_{X})$, and let $\left\{f _{1},\dots,f _{m}\right\}\subset \mathcal{F}$ be a finite $\delta$-packing with $m\ge 2$ and $\delta>0$, so that
    \begin{equation*}
        \left\lVert f _{i}-f _{j}\right\rVert_{L _{2}(\mathbb{P}_{X})}\ge \delta,\qquad \text{for all distinct }i,j\in [m].
    \end{equation*}
    Fix an admissible noise distribution $\mathbb{P}_{\xi}$, and let $J$ be uniformly distributed over $[m]$. Independently of $J$, draw $\tilde{X} _{1},\dots,\tilde{X} _{n}$ independently from $\mathbb{P}_{X}$ and $\xi _{1},\dots,\xi _{n}$ independently from $\mathbb{P}_{\xi}$, with the design variables and noise terms mutually independent. Conditional on $J=j$, generate the responses according to
    \begin{equation*}
        \tilde{Y} _{i}=f _{j}(\tilde{X} _{i})+\xi _{i},\qquad i\in [n].
    \end{equation*}
    Then the minimax risk under squared $L _{2}(\mathbb{P}_{X})$ loss satisfies
    \begin{equation*}
        \mathcal{R}_{n,\epsilon}(\mathcal{F};\mathcal{P}_{\xi})\ge \frac{\delta ^{2}}{4}\left(1-\frac{I(\tilde{X},\tilde{Y};J)+\log 2}{\log m}\right),
    \end{equation*}
    where $I(\tilde{X},\tilde{Y};J)$ denotes the mutual information between the full uncontaminated sample $(\tilde{X},\tilde{Y})$ and the index $J$ under the joint distribution specified above.
\end{theorem}

By Fano's inequality above, to derive a nontrivial lower bound on $\mathcal{R}_{n,\epsilon}(\mathcal{F};\mathcal{P}_{\xi})$, we need to control the mutual information $I(\tilde{X},\tilde{Y};J)$. We therefore restrict our attention to the uncontaminated Gaussian submodel with $\xi _{i}\sim\mathcal{N}\left(0,\sigma ^{2}\right)$, where $\sigma>0$ is chosen so that this noise distribution satisfies the stated assumptions. Since the minimax risk \eqref{eq: minimax risk} takes a supremum over admissible noise distributions and corruption mechanisms, any lower bound for this submodel also applies to the full model. We also note that such inclusion of $\sigma $ is compatible with our later additional assumptions introduced in Section \ref{section: minimax rate under heavy-tailed assumptions} and Section \ref{section: minimax rate under sub-gaussian assumptions}.

The following lemma from \cite{10.1214/25-EJS2419} provides a useful control given this additional assumption. We note that Gaussianity is used only in this lower-bound argument; it is not imposed on the noise distribution in our upper-bound results.

\begin{lemma}[Lemma 2.2 of \cite{10.1214/25-EJS2419}]\label{lemma: control on mutual information}
    Given the same settings as Theorem \ref{theorem: fano inequality} and additionally assuming that $\tilde{\xi }\sim \mathcal{N}(0,\sigma ^{2})$ for some known $\sigma $, we have
    \begin{equation*}
        I(\tilde{X},\tilde{Y};J)\le \frac{n}{2 \sigma ^{2}}\max\limits _{j \in [m]}\left\lVert f _{j}-h\right\rVert_{L _{2}(\mathbb{P}_{X})}^{2},\qquad \forall h \in \mathcal{F}.
    \end{equation*}
\end{lemma}

Combining Theorem \ref{theorem: fano inequality} and Lemma \ref{lemma: control on mutual information}, we come to the following lower bound on the minimax risk \eqref{eq: minimax risk}.

\begin{lemma}[Lemma 2.3 of \cite{10.1214/25-EJS2419}]\label{lemma: general minimax lower bound}
    Given the conditions stated in Theorem \ref{theorem: fano inequality} and Lemma \ref{lemma: control on mutual information}, for any $\eta $ satisfying $\log \mathcal{M}_{\mathcal{F}}^{\text{loc}}(\eta ,c)>4 \left(\frac{n \eta ^{2}}{2 \sigma ^{2}}\vee \log 2\right)$, we have 
    \begin{equation}\label{eq: general lower bound}
        \mathcal{R}_{n,\epsilon}(\mathcal{F};\mathcal{P}_{\xi})\ge \frac{\eta ^{2}}{8c ^{2}},
    \end{equation}
    where we recall that $c$ is the constant from Definition \ref{def_local_entropy}.
\end{lemma}

Lemma \ref{lemma: control on mutual information} follows from an inequality between the mutual information and the KL divergence and a direct calculation of KL divergence between Gaussian random variables. Lemma \ref{lemma: general minimax lower bound} is the result of a direct application of Theorem \ref{theorem: fano inequality} and Lemma \ref{lemma: control on mutual information}. For completeness, we include the proofs in Appendix \ref{subsection: proof of control on mutual information} and Appendix \ref{subsection: proof of general minimax lower bound}.

\section{Minimax rate under Bernstein-Type Function Class and Heavy-Tailed Noise}\label{section: minimax rate under heavy-tailed assumptions}

In this section, we study estimation under a Bernstein-type moment condition on the increments of $\mathcal{F}$ and a finite-variance assumption on the noise. Specifically, we assume that there exists a known constant $L>0$ such that, for every integer $q\ge 2$ and every $f,g\in\mathcal{F}$,
\begin{equation}\label{eq: bernstein condition}
    \left\lVert f-g\right\rVert_{L _{2q}(\mathbb{P}_{X})}^{2q}\le \frac{q!}{2}L ^{q-1}\left\lVert f-g\right\rVert_{L _{2}(\mathbb{P}_{X})}^{2}.
\end{equation}
Here, for $q\ge 1$, the $L _{q}(\mathbb{P}_{X})$-norm is defined by
\begin{equation*}
    \left\lVert f-g\right\rVert_{L _{q}(\mathbb{P}_{X})}:=\left(\int _{\mathcal{X}}\left|f(x)-g(x)\right|^{q}d\mathbb{P}_{X}(x)\right)^{1/q}.
\end{equation*}
For the noise, we retain the assumptions that the conditional expectation of $\tilde{\xi }$ given $\tilde{X}$ is constantly zero, and additionally assume that there exists a known constant $\sigma>0$ such that
\begin{equation}\label{eq: finite conditional second moment of xi}
    \mathbb{E}\left[\tilde{\xi }^{2}\left\lvert\right. \tilde{X}\right]\le \sigma ^{2},\qquad \text{a.s.}
\end{equation}
No higher-order moment or exponential-tail assumption is imposed on $\tilde{\xi }$. Denote the collection of distributions of such $\tilde{\xi }$ as $\mathcal{\mathcal{P}}_{\xi }^{\text{b}}$. We might wonder whether $\mathbb{E}\tilde{\xi }^{2}\le \sigma ^{2}$ for some $\sigma $ is sufficient. However, as we remark later, this is insufficient and \eqref{eq: finite conditional second moment of xi} is necessary for our argument.

We also note here that \eqref{eq: bernstein condition} actually implies an upper bound on the diameter of $\mathcal{F}$. By setting $q=2$ in \eqref{eq: bernstein condition}, we have 
\begin{equation*}
    \left\lVert f-g\right\rVert_{L _{4}(\mathbb{P}_{X})}^{4}\le L \left\lVert f-g\right\rVert_{L _{2}(\mathbb{P}_{X})}^{2}.
\end{equation*}
However, from Jensen's inequality, we also have 
\begin{equation*}
    \left\lVert f-g\right\rVert_{L _{2}(\mathbb{P}_{X})}^{4}\le \left\lVert f-g\right\rVert_{L _{4}(\mathbb{P}_{X})}^{4}.
\end{equation*}
Combining the two displays above, we know that $\left\lVert f-g\right\rVert_{L _{2}(\mathbb{P}_{X})}^{2}\le L, \forall f,g \in \mathcal{F}$, which shows that $d ^{2}\le L$.

Condition \eqref{eq: bernstein condition} is a joint assumption on the function class $\mathcal{F}$ and the design distribution $\mathbb{P}_{X}$. It controls the moments of the squared increments $\left(f(X)-g(X)\right)^{2}$ and yields Bernstein-type concentration for their empirical averages, thereby relating empirical squared distances to their population counterparts; see \cite{10.1214/25-EJS2419}. Under the standing bounded-diameter assumption, this condition is implied by either uniform sup-norm boundedness of $\mathcal{F}$ or a relative sub-Gaussian increment condition introduced later in Section \ref{section: minimax rate under sub-gaussian assumptions}. However, \eqref{eq: bernstein condition} does not require either of these stronger assumptions, as illustrated by the examples in \cite{10.1214/25-EJS2419}.

Although \eqref{eq: bernstein condition} is weaker than the relative sub-Gaussian increment condition, it still controls moments of all orders and implies exponential integrability of the squared increments. Its greater generality concerns the relationship between higher moments and the $L _{2}(\mathbb{P}_{X})$-distance: it does not require the sub-Gaussian scale of an increment to be proportional to its $L _{2}(\mathbb{P}_{X})$-norm. Accordingly, the term ``heavy-tailed'' in this subsection refers to the permitted noise distributions.

The (conditional) second-moment assumption on $\xi$ allows distributions with polynomially decaying tails and, in particular, distributions with an infinite fourth moment. Control of the conditional second moment of $\xi $ nevertheless ensures that, for every $f,g\in\mathcal{F}$,
\begin{equation*}
    \mathbb{E}\left[\tilde{\xi }^{2}\left(f(\tilde{X})-g(\tilde{X})\right)^{2}\right]=\mathbb{E}\left[\left(f(\tilde{X})-g(\tilde{X})\right)^{2}\mathbb{E}\left[\tilde{\xi } ^{2}\left\lvert\right. \tilde{X}\right]\right]\le \sigma ^{2}\left\lVert f-g\right\rVert_{L _{2}(\mathbb{P}_{X})}^{2}.
\end{equation*}
This provides second-moment control of the noise--function products that arise in regression comparisons without requiring higher-order noise moments.

\subsection{Lower Bound}\label{subsection: lower bound under the heavy-tailed setting}

We now argue in the following lemma that $\epsilon \sigma ^{2}\wedge d ^{2}$ is a valid lower bound of $\mathcal{R}_{n,\epsilon}(\mathcal{F};\mathcal{P}_{\xi})$ under the star-shapedness of $\mathcal{F}$, strong adversary $\mathcal{C}$, and Conditions \eqref{eq: bernstein condition}, \eqref{eq: finite conditional second moment of xi}.

\begin{lemma}[Lower bound regarding $\epsilon $]\label{lemma: lower bound regarding epsilon under the bernstein's condition}
    For the nonparametric regression problem formulated in Section \ref{subsection: problem formulation}, with the presence of the strong adversary $\mathcal{C}$ and Conditions \eqref{eq: bernstein condition}, \eqref{eq: finite conditional second moment of xi}, we have 
    \begin{equation}\label{eq: lower bound regarding epsilon under the bernstein's condition}
        \mathcal{R}_{n,\epsilon}(\mathcal{F};\mathcal{P}_{\xi})\gtrsim \epsilon \sigma ^{2}\wedge d ^{2}.
    \end{equation}
\end{lemma}

The proof, deferred to Appendix \ref{subsection: proof of lower bound regarding epsilon under the bernstein's condition}, is based on a two-point construction. Star-shapedness allows us to select two functions in $\mathcal{F}$ whose squared $L_2(\mathbb{P}_X)$-distance is of order $\epsilon\sigma^2\wedge d^2$. We then construct two regression models with these regression functions and a common conditional noise variance bounded by $\sigma^2$. The mixture probabilities depend on the covariate, ensuring that the conditional moment constraints hold while the two observations, under a suitable coupling, differ with marginal probability at most $\epsilon/2$. With probability bounded below by an absolute constant, the number of distinguishing observations does not exceed $\lfloor\epsilon n\rfloor$. The strong adversary can then replace their responses by common values, making the coupled corrupted samples identical. A two-point risk comparison therefore forces any estimator to incur squared error of order $\epsilon\sigma^2\wedge d^2$ under at least one of the two models. This argument uses only star-shapedness and the conditional moment constraints on the noise.

Combining Lemma \ref{lemma: lower bound regarding epsilon under the bernstein's condition} with the general lower bound \eqref{eq: general lower bound} yields
\begin{equation}\label{eq: lower bound under the bernstein condition}
    \mathcal{R}_{n,\epsilon}(\mathcal{F};\mathcal{P}_{\xi})\gtrsim _{\sigma }\left(\left(\eta ^{*}\right)^{2}\vee \epsilon\right)\wedge d ^{2}.
\end{equation}
Recall in \eqref{eq: minimax risk} that we define $\mathcal{R}_{n,\epsilon }(\mathcal{F},\mathcal{P}_{\xi })$ as the worst-case squared estimation error over all noise distributions that belong to $\mathcal{P}_{\xi }$. The following example demonstrates that, this general lower bound cannot be improved even when we restrict our scope to some good noise distributions. Or equivalently, we are free to expand our cover of the lower bound from some well-behaved noise distributions to the general distribution class $\mathcal{P}_{\xi }$ that satisfies the basic conditions $\mathbb{E}\left[\tilde{\xi }\left\lvert\right. \tilde{X}\right]=0$.

\begin{example}[Optimality of $\epsilon $]\label{example: optimality of epsilon under finite second moment}
    Let $\tilde{X} _{i}$ be distributed as a uniform distribution on $[0,1]$. The function class $\mathcal{F}=\left\{\theta \mathbf{1}_{\left\{x \in [0,t]\right\}}: \theta \in [0,1],t \in [0,1]\right\}$. Clearly $\mathcal{F}$ is a star-shaped function class. Also, since $\left|f-g\right|\le 1,\forall f,g \in \mathcal{F}$, we have 
    \begin{equation*}
        \mathbb{E}\left|f-g\right|^{2p}\le \mathbb{E}\left|f-g\right|^{2}\le \frac{p!}{2}\mathbb{E}\left|f-g\right|^{2}.
    \end{equation*}
    Therefore, the moment increment condition \eqref{eq: bernstein condition} holds for $\mathcal{F}$ with $L=1$. Let $\tilde{Y}=f ^{*}(\tilde{X})+\tilde{\xi }$ for some $f ^{*}\in \mathcal{F}$ and $\xi \sim \mathcal{N}(0, 1)$.

    For a fixed and known contamination level $\epsilon $, let $A=[0,\epsilon /2]$, $f _{0}(x)=0, f _{1}(x)=\mathbf{1}_{A}(x)$. Then the squared $L _{2}(\mathbb{P}_{X})$ distance between $f _{0}$ and $f _{1}$ is 
    \begin{equation*}
        \left\lVert f _{1}-f _{0}\right\rVert_{L _{2}(\mathbb{P}_{X})}^{2}=\mathbb{P}_{X}(A)=\frac{\epsilon }{2}.
    \end{equation*}

    Suppose our objective now is to distinguish between the cases when $f ^{*}=f _{0}$ and $f ^{*}=f _{1}$. Since $f _{0}(x)=f _{1}(x)=0$ for $x \notin A$, any efforts must rely on the samples that are supported on $A$. However, consider the counting statistic 
    \begin{equation*}
        N _{A}=\sum\limits _{i=1}^{n}\mathbf{1}_{\left\{\tilde{X} _{i}\in A\right\}}.
    \end{equation*}
    By construction, $N _{A}\sim \operatorname{Bin}\left(n,\epsilon /2\right)$ and consequently by the Chernoff bound, we have 
    \begin{equation*}
        \mathbb{P}\left(N _{A}>\epsilon n\right)=\mathbb{P}\left(N _{A}>2 \mathbb{E}N _{A}\right)\le \exp\left\{-\frac{1}{6}\epsilon n\right\}.
    \end{equation*}
    Consider a simple strategy of $\mathcal{C}$ that $\mathcal{C}$ replaces each sample on $A$ with $0$ when $N _{A}\le \lfloor \epsilon N \rfloor $ and otherwise leaves the samples unchanged. Denote $Q _{0}$ and $Q _{1}$ as the joint distributions of $(X,Y)$ given $f ^{*}=f _{0}$ and $f ^{*}=f _{1}$, respectively. Then by the argument above, we know that on the event $\left\{N _{A}\le \epsilon n\right\}$, the behavior of samples is exactly the same under both hypotheses. Therefore, we have 
    \begin{equation*}
        \operatorname{TV}(Q _{0},Q _{1})\le \exp\left\{-\frac{1}{6}\epsilon n\right\}.
    \end{equation*}
    Now from Le Cam's two-point lower bound (\cite{10.1214/aos/1176348114}), for any estimator $\hat{f}$, we have 
    \begin{equation*}
        \mathcal{R}_{n,\epsilon}(\mathcal{F};\mathcal{P}_{\xi})\ge \frac{\left\lVert f _{1}-f _{0}\right\rVert_{L _{2}(\mathbb{P}_{X})}^{2}}{8}\left(1-\operatorname{TV}(Q _{0},Q _{1})\right)\ge \frac{\epsilon }{16}\left(1-\exp\left\{-\frac{1}{6}\epsilon n\right\}\right)\gtrsim \epsilon ,
    \end{equation*}
    provided $n$ is sufficiently large.
\end{example}

In the next section, we establish a matching upper bound using tree-based Algorithms \ref{algorithm: construction of the pruned tree} and \ref{algorithm: theoretical algorithm}. This therefore determines the minimax rate of \eqref{eq: minimax risk} up to constant factors under the Bernstein's condition \eqref{eq: bernstein condition} and finite conditional second moment \eqref{eq: finite conditional second moment of xi} of $\xi $.

\subsection{Upper Bound}\label{subsection: upper bound under the heavy-tailed setting}

To obtain an upper bound corresponding to the critical radius $\eta ^{*}$ appearing in Section \ref{subsection: lower bound under the heavy-tailed setting}, we construct a hierarchy of localized approximations to $\mathcal{F}$. Each successive level uses finer packing resolutions within smaller neighborhoods, allowing the estimation procedure to explore the function class at progressively decreasing scales. This multiscale strategy builds on the construction developed for constrained mean estimation in the Gaussian sequence model in \cite{neykov2023}. We first introduce the graph notation and then describe the construction.

Let $G=\left(V,E\right)$ be a directed graph whose vertices are labeled by functions in $\mathcal{F}$. For each vertex $u \in V$, define its parent set and offspring set by
\begin{equation*}
    \mathcal{P}(u):=\left\{v \in V:\left(v,u\right)\in E\right\}, \qquad \mathcal{O}(u):=\left\{v \in V:\left(u,v\right)\in E\right\}.
\end{equation*}
The vertices are organized into levels, with $\mathcal{L}(k)$ denoting the collection of vertices at level $k$. The root occupies level $1$, and every edge connects a vertex at level $k$ to a vertex at level $k+1$. Since the pruning step may assign several parents to the same vertex, the resulting structure is a layered directed acyclic graph, which we refer to as the pruned tree. When computing distances, we identify each vertex with its function label; distinct vertices may nevertheless carry the same label.

Choose an arbitrary root $u _{1}^{1}\in\mathcal{F}$ and set $\mathcal{L}(1)=\left\{u ^{1} _{1}\right\}$. Because $d$ is the diameter of $\mathcal{F}$, we have $B _{2}(u ^{1} _{1},d)\cap\mathcal{F}=\mathcal{F}$. Construct a maximal $(d/c)$-packing of this set, take its elements as the vertices of $\mathcal{L}(2)$, and connect the root to each of them. The second level is retained without pruning.

For $k\ge 3$, suppose that the preceding level $\mathcal{L}(k-1)$ has already been finalized. For every $u\in\mathcal{L}(k-1)$, construct a maximal $d/(2^{k-1}c)$-packing of $B _{2}(u,d/2^{k-2})\cap\mathcal{F}$. Each packing point generates a preliminary child of $u$. If the same function occurs in packings associated with different parents, we initially represent these occurrences by separate vertices, so that every preliminary vertex has a unique parent.

Arrange the preliminary vertices according to a fixed ordering, write them as $u _{1}^{k},\dots,u _{N _{k}}^{k}$, and initialize the ordered collection of unprocessed vertices as $\mathcal{U}_{k}=\left[u _{1}^{k},\dots,u _{N _{k}}^{k}\right]$. The pruning procedure repeatedly retains the first unprocessed vertex and merges nearby candidates into it. Specifically, if $u _{i}^{k}$ is the first element of $\mathcal{U}_{k}$, define
\begin{equation*}
    \mathcal{T}_{k}(u _{i}^{k}):=\left\{u _{j}^{k}\in\mathcal{U}_{k}\;\middle|\;j\neq i,\ \left\lVert u _{j}^{k}-u _{i}^{k}\right\rVert_{L _{2}(\mathbb{P}_{X})}\le\frac{d}{2^{k-1}c}\right\}.
\end{equation*}
We retain $u _{i}^{k}$ as the representative of this group. For every $u _{j}^{k}\in\mathcal{T}_{k}(u _{i}^{k})$, let $v _{j}$ be its unique parent. Replace the edge $v _{j}\to u _{j}^{k}$ with $v _{j}\to u _{i}^{k}$, adding the latter only if it is not already present, and discard $u _{j}^{k}$ from the graph. Thus, the retained representative inherits the parent connections of all vertices merged into it.

After these edge reassignments, remove $\left\{u _{i}^{k}\right\}\cup\mathcal{T}_{k}(u _{i}^{k})$ from $\mathcal{U}_{k}$, while keeping $u _{i}^{k}$ in the graph. Continue until no unprocessed vertices remain. The surviving vertices form the finalized level
\begin{equation*}
    \mathcal{L}(k)=\bigcup\limits_{u\in\mathcal{L}(k-1)}\mathcal{O}(u).
\end{equation*}
By construction, any two distinct vertices retained at this level are separated by more than $d/(2^{k-1}c)$, and every discarded candidate lies within this distance of its retained representative.

Applying the same construction recursively defines the complete hierarchy. Algorithm \ref{algorithm: construction of the pruned tree} summarizes this procedure, whose pruning mechanism follows \cite{10.1214/25-AOS2576}. The subsequent Algorithm \ref{algorithm: theoretical algorithm} uses this hierarchy with a different selection criterion tailored to the contaminated regression model we study.

\begin{algorithm}[ht]
    \caption{Construction of the pruned tree}\label{algorithm: construction of the pruned tree}

    \KwInput{Function class $\mathcal{F}$ with diameter $d>0$, root $f _{1}^{1}\in\mathcal{F}$, and constant $c$}
    \KwOutput{The pruned graph $G$}

    Initialize $G$ with the root $f _{1}^{1}$ and set $\mathcal{L}(1)\gets\left\{f _{1}^{1}\right\}$\;

    Construct a maximal $(d/c)$-packing of $B _{2}(f _{1}^{1},d)\cap\mathcal{F}$ and create its vertices as $\mathcal{L}(2)$\;

    Add these vertices to $G$ and connect $f _{1}^{1}$ to each of them\;

    $k\gets 3$\;

    \While{\textnormal{True}}{
        $\mathcal{U}_{k}\gets\varnothing$\;

        \For{$u\in\mathcal{L}(k-1)$}{
            Construct a maximal $d/(2^{k-1}c)$-packing of $B _{2}(u,d/2^{k-2})\cap\mathcal{F}$\;

            Create a separate vertex for each packing point and add these vertices to $G$ and $\mathcal{U}_{k}$\;

            Add an edge from $u$ to each newly created vertex\;
        }

        Arrange $\mathcal{U}_{k}$ according to the fixed ordering\;

        \While{$\mathcal{U}_{k}\neq\varnothing$}{
            Let $\tilde{u}$ be the first vertex in $\mathcal{U}_{k}$\;

            $\mathcal{T}_{k}(\tilde{u})\gets\left\{u\in\mathcal{U}_{k}\;\middle|\;u\neq\tilde{u},\ \left\lVert u-\tilde{u}\right\rVert_{L _{2}(\mathbb{P}_{X})}\le\frac{d}{2^{k-1}c}\right\}$\;

            \For{$u^{\prime}\in\mathcal{T}_{k}(\tilde{u})$}{
                Let $v$ be the unique parent of $u^{\prime}$\;

                Remove the edge $v\to u^{\prime}$\;

                Add the edge $v\to\tilde{u}$ if it is not already present\;

                Delete $u^{\prime}$ from $G$\;
            }

            Remove $\left\{\tilde{u}\right\}\cup\mathcal{T}_{k}(\tilde{u})$ from $\mathcal{U}_{k}$\;
        }

        $\mathcal{L}(k)\gets\bigcup\limits_{u\in\mathcal{L}(k-1)}\mathcal{O}(u)$\;

        $k\gets k+1$\;
    }
\end{algorithm}

The iterative construction and pruning procedures described above yield the following crucial lemmas from \cite{10.1214/25-AOS2576}, which establish the covering and packing properties at each level of the directed tree $G$. Specifically, the nodes at each level $k$ constitute a valid covering and packing of the constraint set $\mathcal{F}$. Moreover, for any node $u ^{k-1}$ at level $k-1$, its offspring set $\mathcal{O}(u ^{k-1})$ forms an effective cover for the localized neighborhood $\mathcal{B}(u ^{k-1}, d/2^{k-2}) \cap \mathcal{F}$. These geometric structural properties play a foundational role in our subsequent analysis of the estimator's performance.

\begin{lemma}[Lemma 3.1 of \cite{10.1214/25-AOS2576}]\label{lemma: covering and packing properties of the constructed tree}
For any $J \ge 3$ and $u ^{J-1}\in \mathcal{L}(J-1)$, we have:\\
\hspace*{0.5em}(1), $\mathcal{L}(J)$ forms a $d/(2^{J-2}c)$-covering of $\mathcal{F}$;\\ 
\hspace*{0.5em}(2), $\mathcal{L}(J)$ forms a $d/(2^{J-1}c)$-packing of $\mathcal{F}$;\\ 
\hspace*{0.5em}(3), the offspring $\mathcal{O}(u ^{J-1})$ form a $d/(2^{J-2}c)$-covering of the set $\mathcal{B}(u ^{J-1}, d/2^{J-2}) \cap \mathcal{F}$;\\ 
\hspace*{0.5em}(4), the cardinality of $\mathcal{O}(u ^{J-1})$ is upper bounded by $\mathcal{M}_{\mathcal{F}}^{\mathrm{loc}}(d/2^{J-2}, 2c)$ for $J \ge 2$.
\end{lemma}

\begin{lemma}[Lemma 3.2 of \cite{10.1214/25-AOS2576}]\label{lemma: bound on the cardinality}
For any $f \in \mathcal{F}$ and $J \ge 2$, $\mathcal{L}(J - 1) \cap B _{2}(f, d/2^{J-2})$ has a cardinality upper bounded by $\mathcal{M}_{\mathcal{F}} ^{\text{loc}}(d/2 ^{J-2},c)\le \mathcal{M}_{\mathcal{F}} ^{\text{loc}}(d/2 ^{J-2},2c)$.
\end{lemma}

\begin{lemma}[Lemma 3.3 of \cite{10.1214/25-AOS2576}]\label{lemma: bound on the distance between the parent and child}
Let $\left[f _{1},f _{2},\dots\right]$ be the nodes of an infinite path in the graph $G$, i.e., $f _{J+1} \in \mathcal{O}(f _J)$ for any $J \in \mathbb{N}$. Then for any integers $1 \le J ^{\prime }\le J$, we have $\left\lVert f _{J ^{\prime }}-f _{J}\right\rVert_{L _{2}(\mathbb{P}_{X})}^{} \le \frac{d(2 + 4c)}{c 2^{J'}}$.
\end{lemma}

We refer the reader to the original source for the proofs. Algorithm \ref{algorithm: construction of the pruned tree} provides increasingly fine approximations of $\mathcal F$ and provides the hierarchy used by Algorithm \ref{algorithm: theoretical algorithm} to approximate \(f^*\) at progressively finer scales.

Intuitively, for any candidate $g \in \mathcal{F}$ that is close enough to the true function $f ^{*}$, we expect the mean squared error (MSE) of the estimation of $Y$ using $g$ to be small, at least compared with the function that is significantly farther. However, the adversary $\mathcal{C}$ breaks this intuition as it is known that the sample mean is not robust to the corruption. This motivates a more robust approach to estimate the MSE in the presence of the adversary---trimmed mean. The mathematical study of trimmed means dates back at least to \cite{daniell1920}, which investigated the removal of extreme observations within a framework of weighted order statistics. More recently, finite-sample guarantees under heavy tails and adversarial contamination have been established for modified trimmed-mean procedures. \cite{10.1214/20-AOS1961} prove the effectiveness of such procedures in both the univariate and multivariate settings. These developments provide methodological foundations for using trimming in the present regression problem.

\begin{lemma}[Trimmed mean, \cite{10.1214/20-AOS1961}]\label{lemma: trimmed mean}
Let $X$ be a real-valued random variable with distribution $\mathbb{P}_{X}$, mean $\mu $, and finite variance $\sigma _{X}^{2}$. Suppose that the observations $X _{1},\dots,X _{n}$ are obtained from $n$ independent copies of $X$ after a strong adversary $\mathcal{C}$ replaces at most $\left\lfloor \epsilon n\right\rfloor$ observations. Assume that $n$ is even, and let $\rho \in \left(0,1\right)$ satisfy $\rho \ge \exp\left\{-n/2\right\}/4$. Define
\begin{equation*}
    \bar{X}:=X-\mu ,\qquad \zeta :=8\epsilon +24\frac{\log\left(\frac{4}{\rho }\right)}{n},
\end{equation*}
and assume that $\zeta <1/4$. For $t\in\left(0,1\right)$, define the tail functional
\begin{equation*}
    \mathcal{E}(t,X):=\max\left\{\mathbb{E}\left[\left|\bar{X}\right|\mathbf{1}_{\left\{\bar{X}\le Q _{t/2}(\bar{X})\right\}}\right],\mathbb{E}\left[\left|\bar{X}\right|\mathbf{1}_{\left\{\bar{X}\ge Q _{1-t/2}(\bar{X})\right\}}\right]\right\},
\end{equation*}
where $Q _{u}(\bar{X})$ denotes the $u$-quantile of $\bar{X}$.

Partition the observed sample into two halves, $\left(X _{1,1},\dots,X _{1,n/2}\right)$ and $\left(X _{2,1},\dots,X _{2,n/2}\right)$. Let $X _{1,1}^{*}\le \dots\le X _{1,n/2}^{*}$ denote the order statistics of the first half, and set
\begin{equation*}
    \alpha :=X _{1,\left\lceil \zeta n/2\right\rceil}^{*},\qquad \beta :=X _{1,\left\lfloor \left(1-\zeta \right)n/2\right\rfloor}^{*}.
\end{equation*}
Define the truncation function
\begin{equation*}
    \phi _{\alpha ,\beta }(x):=
    \begin{cases}
        \alpha , & x\le \alpha ,\\
        x, & \alpha <x\le \beta ,\\
        \beta , & x>\beta ,
    \end{cases}
\end{equation*}
and the modified trimmed mean estimator
\begin{equation*}
    \operatorname{Tm}(X):=\frac{2}{n}\sum\limits _{i=1}^{n/2}\phi _{\alpha ,\beta }\left(X _{2,i}\right).
\end{equation*}
Then the following bounds hold:
\begin{enumerate}
    \item With probability at least $1-\rho $,
    \begin{equation}\label{eq: concentration inequality 1 of trimmed mean}
        \left|\operatorname{Tm}(X)-\mu \right|\le 4\mathcal{E}(4\zeta ,X)+2\sigma _{X}\sqrt{\frac{2\log\left(\frac{4}{\rho }\right)}{n}}.
    \end{equation}
    \item With probability at least $1-4\exp\left\{-\zeta n/24\right\}$,
    \begin{equation}\label{eq: concentration inequality 2 of trimmed mean}
        \left|\operatorname{Tm}(X)-\mu \right|\le 10\sqrt{\zeta }\sigma _{X}.
    \end{equation}
\end{enumerate}
\end{lemma}

\begin{remark*}[]
    We will later use $\mathbb{P}\left(\bar{X}\ge Q _{1-t/2}(\bar{X})\right)\le t/2$ to bound $\mathcal{E} (t,X)$. However, this is not true in general for random variables with atoms. As mentioned in \cite{10.1214/20-AOS1961}, this issue can be amended by adding a sufficiently small absolutely continuous perturbation to $X$ (for example, $rG$ for sufficiently small $r$ and $G$ is a standard Gaussian random variable).
\end{remark*}

We now introduce the pairwise comparison statistic underlying Algorithm \ref{algorithm: theoretical algorithm}. For any $f,g\in\mathcal{F}$ satisfying $\left\lVert f-g\right\rVert_{L _{2}(\mathbb{P}_{X})}>0$, define the normalized squared-loss difference
\begin{equation}\label{eq: definition of the pairwise comparison statistic}
    Z _{f,g}:=\frac{\left[Y-f(X)\right]^{2}-\left[Y-g(X)\right]^{2}}{\left\lVert f-g\right\rVert_{L _{2}(\mathbb{P}_{X})}}.
\end{equation}

Let $\operatorname{Tm}_{Z _{f,g}}(X,Y)$ denote the estimator from Lemma \ref{lemma: trimmed mean} applied to the normalized loss differences computed from the observed sample $\left\{\left(X _{i},Y _{i}\right)\right\}_{i=1}^{n}$. Define the pairwise comparison test
\begin{equation}\label{eq: pairwise comparison test}
    \psi _{f,g}(X,Y):=\mathbf{1}_{\left\{\operatorname{Tm}_{Z _{f,g}}(X,Y)>0\right\}}.
\end{equation}
We say that $f$ is dominated by $g$, denoted by $f\prec g$ or $g\succ f$, whenever $\psi _{f,g}(X,Y)=1$. This relation need not be transitive: $f\prec g$ and $g\prec h$ do not necessarily imply $f\prec h$. Moreover, for any two candidates $f,g \in \mathcal{F}$, \eqref{eq: pairwise comparison test} prevents the occurrence that both $f \prec g$ and $g \prec f$ happen simultaneously. To see this, let $z _{(1)}\le z _{(2)}\le \dots \le z _{(n/2)}$ be the order statistics of the first half of $Z _{f,g}$ and $\omega _{1},\omega _{2},\dots,\omega _{n/2}$ be the second half of $Z _{f,g}$. Let $k=\lceil \zeta n/2 \rceil $. By definition of the trimmed mean estimator, we have $\operatorname{Tm}(Z _{f,g})=\frac{2}{n}\sum\limits _{i=1}^{n/2}\phi _{z _{(k)},z _{(n/2-k)}}(\omega _{i})$. Since $Z _{f,g}=-Z _{g,f}$, we then have by the property of $\phi $: 
\begin{equation*}
    \operatorname{Tm}(Z _{g,f})=\frac{2}{n}\sum\limits _{i=1}^{n/2}\phi _{-z _{(n/2+1-k)},-z _{(k+1)}}(-\omega _{i})=-\frac{2}{n}\sum\limits _{i=1}^{n/2}\phi _{z _{(k+1)},z _{(n/2+1-k)}}(\omega _{i})\le -\operatorname{Tm}(Z _{f,g}).
\end{equation*}
Hence, $\operatorname{Tm}(Z _{f,g})>0$ implies that $\operatorname{Tm}(Z _{g,f})<0$.

Applying Lemma \ref{lemma: trimmed mean} to the normalized loss differences yields the following concentration bound, which is a key ingredient in establishing the upper bound for the estimator constructed by Algorithm \ref{algorithm: theoretical algorithm}.

\begin{lemma}\label{lemma: trimmed mean on the current problem}
    Suppose that $\mathcal{F}$ satisfies the Bernstein-type moment condition \eqref{eq: bernstein condition} and that the noise $\tilde{\xi }$ satisfies Condition \eqref{eq: finite conditional second moment of xi}. Let $f,g\in\mathcal{F}$ be fixed functions, and suppose that, for some $\delta >0$,
    \begin{equation*}
        \left\lVert g-f ^{*}\right\rVert_{L _{2}(\mathbb{P}_{X})}\le \delta ,\qquad \left\lVert f-g\right\rVert_{L _{2}(\mathbb{P}_{X})}\ge C\delta ,
    \end{equation*}
    where $C :=c/2-1\ge 3$ is a universal constant. Let $\rho $ and the corresponding trimming parameter $\zeta $ satisfy the conditions of Lemma \ref{lemma: trimmed mean}. Then, with probability at least $1-\rho $,
    \begin{equation}\label{eq: concentration inequality in the current problem}
        \left|\operatorname{Tm}_{Z _{f,g}}(X,Y)-\mathbb{E}Z _{f,g}\right|\lesssim \sqrt{L+\sigma ^{2}}\left[\sqrt{\epsilon }+\sqrt{\log\left(\frac{4}{\rho }\right)/n}\right],
    \end{equation}
    where the expectation is taken under the uncontaminated distribution and the implicit constant is universal.
\end{lemma}

By bounding $\mathbb{E}Z _{f,g}$ and properly selecting the constants in Lemma \ref{lemma: trimmed mean on the current problem}, we obtain the following corollary, which is analogous to Lemma 3.8 of \cite{10.1214/25-EJS2419}.

\begin{corollary}[]\label{corollary: concrete application of the trimmed mean to bound the error probability}
    Under the same settings of Lemma \ref{lemma: trimmed mean on the current problem}, additionally assume that $\epsilon \le \frac{1}{64}$ and $a \le \frac{1}{192}$. If $\delta ^{2}\gtrsim _{L,\sigma }\epsilon $, we have 
    \begin{equation*}
        \mathbb{P}\left(\psi _{f,g}(X,Y)=0\right)\le 4 \exp\left\{-\frac{an\delta ^{2}}{L+\sigma ^{2}}\right\}.
    \end{equation*}
    Equivalently, the pairwise comparison test \eqref{eq: pairwise comparison test} correctly selects $g$ over $f$ with probability greater than $1-4 \exp\left\{-an\delta ^{2}/(L+\sigma ^{2})\right\}$
\end{corollary}

The proofs of Lemma \ref{lemma: trimmed mean on the current problem} and Corollary \ref{corollary: concrete application of the trimmed mean to bound the error probability} are provided in Appendix \ref{subsection: proof of trimmed mean on the current problem} and Appendix \ref{subsection: proof of concrete application of the trimmed mean to bound the error probability}, respectively.

To describe the selection procedure on the pruned tree constructed in Algorithm \ref{algorithm: construction of the pruned tree}, we first introduce some notation. For a nonempty finite set $S\subset\mathcal{F}$, a scale $\delta >0$, and a function $f\in S$, define the set of admissible competitors of $f$ by
\begin{equation*}
    S _{f}:=\left\{g\in S:g\succ f,\ \left\lVert f-g\right\rVert_{L _{2}(\mathbb{P}_{X})}\ge C\delta \right\},
\end{equation*}
where $C$ is the constant appearing in Lemma \ref{lemma: trimmed mean on the current problem}. The dependence of $S _{f}$ on $S$ and $\delta $ is suppressed for convenience. We define the \textit{discrepancy} of $f$ relative to $S$ at scale $\delta $ as
\begin{equation}\label{eq: discrepancy}
    H(\delta ,f,S):=
    \begin{cases}
        \displaystyle\max\limits _{g\in S _{f}}\left\lVert f-g\right\rVert_{L _{2}(\mathbb{P}_{X})}, & \text{if }S _{f}\neq\varnothing ,\\
        0, & \text{if }S _{f}=\varnothing .
    \end{cases}
\end{equation}

Combining Corollary \ref{corollary: concrete application of the trimmed mean to bound the error probability} with the discrepancy criterion \eqref{eq: discrepancy} yields the following selection guarantee. It ensures that the probability of incorrect selection is small when at least one function from $S$ is close to the true function $f ^{*}$.

\begin{lemma}\label{lemma: crucial lemma for correct selection}
Suppose that the conditions of Corollary \ref{corollary: concrete application of the trimmed mean to bound the error probability} hold at scale $\delta >0$. Let $S=\left\{f _{1},\dots,f _{M}\right\}\subset\mathcal{F}$ be a fixed finite set satisfying
\begin{equation*}
    \min\limits _{i\in[M]}\left\lVert f _{i}-f ^{*}\right\rVert_{L _{2}(\mathbb{P}_{X})}\le\delta .
\end{equation*}
Select
\begin{equation*}
    f _{S}\in\mathop{\arg\min}\limits _{f\in S}H(\delta ,f,S),
\end{equation*}
with ties resolved according to a fixed ordering of $S$. Then
\begin{equation*}
    \mathbb{P}\left(\left\lVert f _{S}-f ^{*}\right\rVert_{L _{2}(\mathbb{P}_{X})}\ge\left(C+1\right)\delta \right)\le 4M\exp\left\{-\frac{an\delta ^{2}}{L+\sigma ^{2}}\right\},
\end{equation*}
where $a>0$ is the constant appearing in Corollary \ref{corollary: concrete application of the trimmed mean to bound the error probability}.
\end{lemma}

The proof follows from Corollary \ref{corollary: concrete application of the trimmed mean to bound the error probability}, the triangle inequality, and a union bound and is given in Appendix \ref{subsection: proof of the crucial lemma for correct selection}.

We now describe the selection procedure on the pruned tree constructed in Algorithm \ref{algorithm: construction of the pruned tree}. Although the tree may have infinitely many levels, the procedure traverses only finitely many of them. Starting from the root $f _{1}^{1}$, we set $f _{1}^{*}:=f _{1}^{1}$. At each iteration $k\ge 1$, we select $f _{k+1}^{*}$ from the offspring set $\mathcal{O}(f _{k}^{*})$ by minimizing the discrepancy $H\left(\delta =\frac{d}{2^{k}(C+1)},f,\mathcal{O}(f _{k}^{*})\right)$ over $f\in\mathcal{O}(f _{k}^{*})$. Ties are resolved by selecting the first minimizer according to a fixed ordering of the nodes. The procedure is summarized in Algorithm \ref{algorithm: theoretical algorithm}, where the stopping index $J^{*}$ will be specified later.

\begin{algorithm}[ht]
\caption{Sequential selection on the pruned tree}\label{algorithm: theoretical algorithm}
\SetKwComment{Comment}{/* }{ */}
\KwInput{Function class $\mathcal{F}$, observations $\left(X _{1},Y _{1}\right),\dots,\left(X _{n},Y _{n}\right)$, the pruned tree $G$ from Algorithm \ref{algorithm: construction of the pruned tree}, and a stopping index $J^{*}$}
\KwOutput{The sequence $\mathcal{F}^{*}=\left[f _{1}^{*},\dots,f _{J^{*}}^{*}\right]$}

$f _{1}^{*}\gets f _{1}^{1}$\;

$k\gets 1$\;

$\mathcal{F}^{*}\gets\left[f _{1}^{*}\right]$\;

\While{$k< J^{*}$}{
    Select $f _{k+1}^{*}\in\mathop{\arg\min}\limits _{f\in\mathcal{O}(f _{k}^{*})}H\left(\frac{d}{2^{k}(C+1)},f,\mathcal{O}(f _{k}^{*})\right)$ \Comment*[r]{Choose the first minimizer in the fixed node ordering}

    Append $f _{k+1}^{*}$ to $\mathcal{F}^{*}$\;

    $k\gets k+1$\;
}

\Return{$\mathcal{F}^{*}=\left[f _{1}^{*},\dots,f _{J ^{*}}^{*}\right]$}\;
\end{algorithm}

\begin{corollary}[One-step error bound along the tree traversal]
\label{corollary: probability of error during the trajectory of the tree traversal}
Under the assumptions of the preceding tournament selection lemma, denote by $f_j^*$ the node selected at level $j$ by Algorithm \ref{algorithm: theoretical algorithm}, with $f_1^*$ being the root. Thus, $f_j^*$ is selected after $j-1$ iterations. For each $j\ge 1$, define
\begin{equation*}
    \delta_j:=\frac{d}{2^{j-2}c}=\frac{d}{2^{j-1}(C+1)},\qquad\eta_j:=\sqrt{\frac{a}{L+\sigma^2}}\,\delta_j,
\end{equation*}
and, for each level $k\ge1$, define the failure event
\begin{equation}\label{eq: failure event}
    A_k:=\left\{\left\lVert f _{k}^{*}-f ^{*}\right\rVert_{L _{2}(\mathbb{P}_{X})}^{}>\frac{d}{2^{k-1}}=(C+1)\delta _{k}\right\}.
\end{equation}
Suppose that $\epsilon \le b\delta _{j}^{2}/(L+\sigma ^{2})$, where $b>0$ is a sufficiently small constant determined by the preceding testing bound and independent of $d$. Then
\begin{equation*}
    \mathbb{P}\left(A _{j-1}^{\complement }\cap A _{j}\right)\le4\left[\mathcal{M}_{\mathcal{F}}^{\mathrm{loc}}\left(\frac{d}{2^{j-2}},2c\right)\right]^2\exp\left\{-n\eta_j^2\right\}.
\end{equation*}
\end{corollary}

We note that Corollary \ref{corollary: probability of error during the trajectory of the tree traversal} is a direct deduction from Lemma \ref{lemma: bound on the cardinality}, Algorithm \ref{algorithm: theoretical algorithm}, and Corollary \ref{corollary: concrete application of the trimmed mean to bound the error probability}. We provide its proof in Appendix \ref{subsection: proof of probability of error during the trajectory of the tree traversal}.

With the help of Corollary \ref{corollary: probability of error during the trajectory of the tree traversal}, we can ``telescope'' the chain $\mathcal{F}^{*}$ to analyse the increment of the traversal. As a consequence, we are finally prepared to present the main theorem induced by Algorithm \ref{algorithm: theoretical algorithm}, which also serves as an upper bound for the minimax rate of \eqref{eq: minimax risk} uniformly over the admissible star-shaped function class $\mathcal{F}$ and adversary $\mathcal{C}$.

\begin{theorem}[]\label{theorem: main theorem of the algorithm}
    Let $J ^{*}$ be the maximal integer such that
    \begin{equation}\label{eq: definition of J^*}
        n\eta _{j}^{2}>4 \log \left[\mathcal{M}_{\mathcal{F}}^{\text{loc}}\left(\frac{d}{2 ^{j-2}},2c\right)\right]\vee \log 2.
    \end{equation}
    If this condition never holds, we let $J ^{*}=1$.
    
    For the nonparametric regression problem formulated in Section \ref{subsection: problem formulation}, suppose that $\epsilon \le \frac{1}{64}$ as required in Corollary \ref{corollary: concrete application of the trimmed mean to bound the error probability}. Let $f _{J ^{*}}^{*}$ be the $J ^{*}$th element of the sequence $\mathcal{F}^{*}$ after $J ^{*}-1$ iterations of Algorithm \ref{algorithm: theoretical algorithm}. Then we have 
    \begin{equation}\label{eq: upper bound of the theoretical algorithm under the bernstein condition}
        \sup\limits _{f ^{*}\in \mathcal{F}}\sup\limits _{\xi :\mathbb{P}_{\xi }\in \mathcal{P}_{\xi }^{\text{b}}}\sup\limits _{\mathcal{C}\in \mathfrak{C}_{n,\epsilon }}\mathbb{E}_{\tilde{X},\tilde{Y}}\left\lVert f _{J ^{*}}^{*}(\mathcal{C}(\tilde{X},\tilde{Y}))-f ^{*}\right\rVert_{L _{2}(\mathbb{P}_{X})}^{2}\lesssim \left(\eta _{J ^{*}}^{2}\vee \epsilon \right)\wedge d ^{2}.
    \end{equation}
    Here a constant that only depends on $L,\sigma ,c$ is omitted.
\end{theorem}

\begin{remark*}[Finiteness of the tree construction]\label{remark: finiteness of the tree construction}
Boundedness of $\mathcal{F}$ in $L_{2}\left(\mathbb{P}_{X}\right)$ does not, by itself, guarantee that its packing numbers are finite. For example, the unit ball of an infinite-dimensional Hilbert space is bounded but contains an infinite collection of points with pairwise distances bounded away from zero. Thus, boundedness alone would not justify constructing a finite tree to an arbitrary depth.

For our estimator, however, it suffices to construct the tree through the stopping level $J^{*}$. When $J^{*}\ge 2$, the local entropy condition defining $J^{*}$ guarantees that the local packing number controlling this level is finite. By the monotonicity of local packing numbers under star-shapedness established in Lemma \ref{lemma: nonincreasingness of the local entropy}, the corresponding packing numbers at all preceding levels are finite as well. Consequently, each node has finitely many offspring, and induction over the levels shows that the entire truncated tree is finite. When $J^{*}=1$, the estimator simply returns the root, and no packing construction is needed. Therefore, the construction does not require an additional assumption that $\mathcal{F}$ is totally bounded. This observation ensures that the truncated construction is finite, although its computational cost can still be prohibitive.
\end{remark*}

Theorem \ref{theorem: main theorem of the algorithm} implies that $\mathcal{R}_{n,\epsilon }(\mathcal{F};\mathcal{P}_{\xi})\lesssim _{L,\sigma }\left(\eta _{J ^{*}}^{2}\vee \epsilon \right)\wedge d ^{2}$. Recalling the \hyperref[def_critical_radius]{definition} of $\eta ^{*}$, we obtain the following corollary that the upper bound \eqref{eq: upper bound of the theoretical algorithm under the bernstein condition} matches the lower bound \eqref{eq: lower bound under the bernstein condition}, characterizing the exact minimax rate of $\mathcal{R}_{n,\epsilon}(\mathcal{F};\mathcal{P}_{\xi})$ under \eqref{eq: bernstein condition} and \eqref{eq: finite conditional second moment of xi}.

\begin{corollary}[]\label{corollary: minimax rate under the bernstein condition}
    For the nonparametric regression problem formulated in Section \ref{subsection: problem formulation}, under Conditions \eqref{eq: bernstein condition} and \eqref{eq: finite conditional second moment of xi}, the minimax rate of $\mathcal{R}_{n,\epsilon}(\mathcal{F};\mathcal{P}_{\xi})$ is 
    \begin{equation}\label{eq: minimax rate under the bernstein condition}
        \mathcal{R}_{n,\epsilon}(\mathcal{F};\mathcal{P}_{\xi})\asymp _{L,\sigma }\left(\left(\eta ^{*}\right)^{2}\vee \epsilon\right)\wedge d ^{2}.
    \end{equation}
\end{corollary}

Corollary \ref{corollary: minimax rate under the bernstein condition} originates from Theorem \ref{theorem: main theorem of the algorithm} and the \hyperref[def_critical_radius]{definition} of $\eta ^{*}$. We defer the proofs of Theorem \ref{theorem: main theorem of the algorithm} and Corollary \ref{corollary: minimax rate under the bernstein condition} to Appendices \ref{subsection: proof of main theorem of the algorithm} and \ref{subsection: proof of minimax rate under the bernstein condition}, respectively.

\begin{remark*}[]
    Recall that we previously argued that the conditional second-moment assumption \eqref{eq: finite conditional second moment of xi} is essential to the minimax rate. Intuitively, such a condition prevents the noise variance from concentrating on regions where competing regression functions differ. The following example shows that an unconditional second-moment bound alone is insufficient for both the variance estimate used in the proof and the resulting minimax guarantee.
\end{remark*}

\begin{example}[Insufficiency of unconditional variance control]\label{example: insufficiency of unconditional variance control}
Let $X\sim\mathrm{U}[0,1]$ and consider the fixed function class previously discussed in Example \ref{example: optimality of epsilon under finite second moment}:
\begin{equation*}
    \mathcal{F}:=\left\{f_{\theta,t}(x)=\theta\mathbf{1}_{[0,t]}(x):\theta,t\in[0,1]\right\}.
\end{equation*}
We know that this class is star-shaped about zero and has diameter $d=1$. Moreover, the Bernstein condition holds with $L=1$. Also, in this example, we consider the uncontaminated case, where $\epsilon =0$ and $\left(\tilde{X}_{i},\tilde{Y}_{i},\tilde{\xi }_{i}\right)=\left(X _{i},Y _{i},\xi _{i}\right), i=1,2,\dots,n$.

For $q\in(0,1/2)$, write $h_q:=\mathbf{1}_{[0,q]}$ and let $G\sim \mathcal{N}(0,1)$ be independent of $X$. Define
\begin{equation*}
    \xi_q:=\frac{\sigma}{\sqrt{2q}}h_q(X)G.
\end{equation*}
Then
\begin{equation*}
    \mathbb{E}[\xi_q\mid X]=0,\qquad \mathbb{E}\xi_q^2=\frac{\sigma^2}{2}<\sigma^2,\qquad \mathbb{E}[\xi_q^2\mid X]=\frac{\sigma^2}{2q}h_q(X).
\end{equation*}
Thus, the unconditional second moment is uniformly bounded, while the conditional variance on $[0,q]$ diverges as $q\downarrow 0$.

Take $f^*=g=0$ and $f=h_q$. Since $\|h_q\|_{L_2(\mathbb{P}_X)}=\sqrt q$, the pairwise comparison statistic satisfies
\begin{equation*}
    Z_{h_q,0}:=\frac{(Y-h_q(X))^2-Y^2}{\|h_q\|_{L_2(\mathbb{P}_X)}}=\frac{h_q(X)}{\sqrt q}(1-2\xi_q).
\end{equation*}
Consequently,
\begin{equation*}
    \mathbb{E}Z_{h_q,0}=\sqrt q,\qquad \operatorname{Var}(Z_{h_q,0})=1-q+\frac{2\sigma^2}{q}.
\end{equation*}
Hence the variance cannot be bounded uniformly by a constant depending only on $L$ and $\sigma$, although the function class and its diameter remain fixed.

The same construction also demonstrates failure of the critical-radius minimax guarantee. Let $\tilde{\mathcal{P}}_\xi$ denote the noise class obtained by requiring only $\mathbb{E}[\xi\mid X]=0$ and $\mathbb{E}\xi^2\le\sigma^2$. Let $\mathbb{P}_0$ and $\mathbb{P} _1$ be the laws of $(X,Y)=(\tilde{X},\tilde{Y})$ under $f_0=0$ and $f_1=h_q$, respectively, using the same noise $\xi_q$ under both hypotheses. The two conditional response distributions coincide outside $[0,q]$ and are Gaussian with common variance $\sigma^2/(2q)$ inside $[0,q]$. Therefore,
\begin{equation*}
    \left\lVert f _{1}-f _{0}\right\rVert_{L _{2}(\mathbb{P}_{X})}^{2}=q,\qquad \mathrm{KL}\left(\mathbb{P}_0^{\otimes n}\,\middle\|\,\mathbb{P}_1^{\otimes n}\right)=nq\frac{1}{2(\sigma^2/(2q))}=\frac{nq^2}{\sigma^2}.
\end{equation*}
Choose $q=\frac{1}{4}\wedge\frac{\sigma}{4\sqrt n}$. The preceding divergence is at most $1/16$, so Pinsker's inequality and the two-point testing bound give
\begin{equation*}
    \mathcal{R}_{n,0}(\mathcal{F};\tilde{\mathcal{P}}_\xi)\ge\frac{q}{8}\left(1-\mathrm{TV}\left(P_0^{\otimes n},P_1^{\otimes n}\right)\right)\ge\frac{q}{16}\gtrsim1\wedge\frac{\sigma}{\sqrt n}.
\end{equation*}
The dependence of $q$ on $n$ is permitted because the minimax risk takes a supremum over the noise class for each sample size; $\mathcal{F}$ and $\mathbb{P}_X$ remain fixed.

On the other hand, this function class has small metric entropy. Indeed, suppose $t ^{\prime }\le t$. Then we have
\begin{equation*}
    \left\lVert f_{\theta,t}-f_{\theta',t'}\right\rVert_{L_2(\mathbb{P}_X)}^{2}=t ^{\prime }(\theta -\theta ^{\prime })^{2}+\theta ^{2}(t-t ^{\prime })\le \left|\theta -\theta ^{\prime }\right|^{2}+\left|t-t ^{\prime }\right|.
\end{equation*}
Discretizing $\theta$ and $t$ with mesh sizes $\rho/4$ and $\rho^2/16$, respectively, gives a covering with radius strictly smaller than $\rho/2$ and at most a constant multiple of $\rho^{-3}$ elements (think about the $t$-$\theta $ grid). Consequently,
\begin{equation*}
    \mathcal{M}(\rho,\mathcal{F})\lesssim\rho^{-3},\qquad 0<\rho\le1,
\end{equation*}
We argue the following inequality:
\begin{equation*}
    (\eta^*)^2\lesssim_c\frac{\log(en)}{n}.
\end{equation*}
To see this, take any $r>0$ satisfying $nr ^{2}\le \log \mathcal{M}_{\mathcal{F}}^{\text{loc}}(r,c)$. We split the region of $r$ as follows:\\ 
\hspace*{0.5em}(i), $\bm{r< n ^{-1/2}}$: in this case, we trivially have $r ^{2}\le \frac{1}{n}\le \frac{\log (en)}{n}$;\\ 
\hspace*{0.5em}(ii), $\bm{r \ge n ^{-1/2}}$: in this case, every $r/c$-packing of $B_2(f,r)\cap\mathcal{F}$ is also a $1/c \sqrt{n}$-packing of the whole class: its points belong to $\mathcal{F}$, and their pairwise distances are at least $r/c\ge1/(c\sqrt n)$. Therefore,
\begin{equation*}
    \mathcal M_{\mathcal F}^{\mathrm{loc}}(r,c)\le \mathcal M\left(\frac{1}{c\sqrt n},\mathcal F\right)\lesssim c^3 n^{3/2}.
\end{equation*}
Taking logarithms gives
\begin{equation*}
    nr^2 \lesssim \log(c^3)+\frac{3}{2}\log n\lesssim \log(en).
\end{equation*}
Since $r$ is any value that satisfies $nr ^{2}\le \log \mathcal{M}_{\mathcal{F}}^{\text{loc}}(r,c)$, taking the limit yields the desired inequality.

For every fixed $\sigma>0$, the minimax lower bound above is asymptotically larger than the minimax rate determined by the critical radius. Therefore, replacing conditional variance control by an unconditional second-moment bound can invalidate the claimed minimax rate even when $\epsilon=0$.
\end{example}

\section{Minimax Rate under Sub-Gaussian Assumptions}\label{section: minimax rate under sub-gaussian assumptions}

In this section, we study the minimax rate under the sub-Gaussian conditions on $\xi $ and the increment of functions in $\mathcal{F}$. Namely, we assume that there exist (known) constants $L _{\text{sg}}$ and $L _{\xi }$ such that for each integer $q \ge 1$:
\begin{equation}\label{eq: sub-gaussian condition}
    \begin{aligned}
        &\left\lVert f-g\right\rVert_{L _{q}(\mathbb{P}_{X})}^{}\le L _{\text{sg}}\sqrt{q}\left\lVert f-g\right\rVert_{L _{2}(\mathbb{P}_{X})}^{},\\ 
        &\left\lVert \xi \right\rVert_{q}^{}\le L _{\xi } \sqrt{q}\left\lVert \xi \right\rVert_{2}^{},
    \end{aligned}
\end{equation}
Recall that we previously assumed an upper bound on the conditional variance of $\xi $ and demonstrated that it is necessary for the derivations of minimax rate \eqref{eq: minimax rate under the bernstein condition} under the Bernstein's condition \eqref{eq: bernstein condition}. However, we find that, under the stronger sub-Gaussian type conditions \eqref{eq: sub-gaussian condition}, only the more general marginal control on the second moment of $\xi $ is required. Namely, in this section, $\sigma ^{2}$ is the variance proxy of $\xi $ as we assume
\begin{equation}\label{eq: marginal control on the second moment of xi}
    \mathbb{E}\xi ^{2}\le \sigma ^{2}.
\end{equation}
We denote the distribution family of the noise satisfying \eqref{eq: sub-gaussian condition} and \eqref{eq: marginal control on the second moment of xi} as $\mathcal{P}_{\xi }^{\text{sg}}$.

Though we say that \eqref{eq: sub-gaussian condition} is stronger compared with \eqref{eq: bernstein condition}, the implication from \eqref{eq: sub-gaussian condition} to \eqref{eq: bernstein condition} requires some care regarding the dependence of the constants on the diameter. Assume without loss of generality that $L_{\mathrm{sg}}\ge 1$. The first inequality in \eqref{eq: sub-gaussian condition} gives, for every integer $q\ge 1$,
\begin{equation*}
    \|f-g\|_{L_{2q}(\mathbb{P}_X)}^{2q}\le (2q)^q L_{\mathrm{sg}}^{2q}\left\lVert f-g\right\rVert_{L _{2}(\mathbb{P}_{X})}^{2q}\le q!\,(2eL_{\mathrm{sg}}^2)^q d^{2q-2}\left\lVert f-g\right\rVert_{L _{2}(\mathbb{P}_{X})}^{2},
\end{equation*}
where we used $\left\lVert f-g\right\rVert_{L _{2}(\mathbb{P}_{X})}^{}\le d$ and $q^q\le e^q q!$. Setting $B=2eL_{\mathrm{sg}}^2$, we have $B\ge 1$ and $2B^q\le (2B^2)^{q-1}$ for $q\ge 2$. Consequently,
\begin{equation*}
    \|f-g\|_{L_{2q}(\mathbb{P}_X)}^{2q}\le \frac{q!}{2}L^{q-1}\|f-g\|_{L_2(\mathbb{P}_X)}^2,\qquad L=8e^2L_{\mathrm{sg}}^4d^2.
\end{equation*}
Thus, for a class of finite diameter, \eqref{eq: sub-gaussian condition} implies \eqref{eq: bernstein condition}, but the resulting moment constant is allowed to depend on $d$. The noise assumption is not needed for this implication.

This dependence reflects an essential difference between the two conditions. Recall that \eqref{eq: bernstein condition} implies $d^2\le L$. Therefore, fixing $L$ automatically bounds the diameter. In contrast, the relative sub-Gaussian condition is invariant under multiplication of every function by the same positive scalar, so fixing $L_{\mathrm{sg}}$ does not control $d$, as shown by the following example.

\begin{example}[]\label{example: unboundedness of sub-gaussian function class}
    Consider an uncontaminated scenario. Let $X\sim\operatorname{Bernoulli}(1/2)$ and, for any $d>0$, consider the convex, complete, star-shaped class
    \begin{equation*}
        \mathcal{F}_d=\{f_\theta:\theta\in[0,d]\},\qquad f_\theta(x)=\sqrt{2}\,\theta x.
    \end{equation*}
    For any $\theta,\theta'\in[0,d]$,
    \begin{equation*}
        \|f_\theta-f_{\theta'}\|_{L_2(\mathbb{P}_X)}=|\theta-\theta'|,
    \end{equation*}
    so the diameter of $\mathcal{F}_d$ is exactly $d$. Moreover, for every $p\ge 2$,
    \begin{equation*}
        \|f_\theta-f_{\theta'}\|_{L_p(\mathbb{P}_X)}=2^{1/2-1/p}|\theta-\theta'|\le \sqrt{p}\,\|f_\theta-f_{\theta'}\|_{L_2(\mathbb{P}_X)}.
    \end{equation*}
    Hence the relative sub-Gaussian condition holds with $L_{\mathrm{sg}}=1$, uniformly over all $d>0$. However, applying the general moment condition with $p=2$ to the endpoints $f_0$ and $f_d$ requires
    \begin{equation*}
        2d^4=\|f_d-f_0\|_{L_4(\mathbb{P}_X)}^4\le L\|f_d-f_0\|_{L_2(\mathbb{P}_X)}^2=Ld^2.
    \end{equation*}
    Thus any admissible moment constant must satisfy $L\ge 2d^2$, even though $L_{\mathrm{sg}}$ remains equal to one.
\end{example}

This example also illustrates why an additional diameter bound is relevant to a uniform statistical guarantee. Consider the clean regression model $Y=f_\theta(X)+\xi$, with $\xi\sim \mathcal{N}(0,1)$ independent of $X$. On the event
\begin{equation*}
    A_n=\{X_1=\cdots=X_n=0\},\qquad \mathbb{P}(A_n)=2^{-n},
\end{equation*}
the conditional distribution of the observations is identical for all $\theta$. In particular, the endpoints $f_0$ and $f_d$ cannot be distinguished on this event. Using
\begin{equation*}
    \frac{\|u-f_0\|_2^2+\|u-f_d\|_2^2}{2}=\left\|u-\frac{f_0+f_d}{2}\right\|_2^2+\frac{d^2}{4},
\end{equation*}
we obtain, even without contamination,
\begin{equation*}
    \inf_{\hat{f}}\sup_{\theta\in[0,d]}\mathbb{E}_\theta\|\hat{f}-f_\theta\|_2^2\ge \frac{d^2}{4}\,2^{-n}.
\end{equation*}
Meanwhile, $\mathcal{F}_d$ is isometric to the interval $[0,d]$, whose local packing entropy is bounded by a constant for a fixed packing ratio. Consequently, $\left(\eta ^{*}\right)^2\asymp n^{-1}$ whenever $d\ge n^{-1/2}$. Allowing $d$ to grow arbitrarily therefore rules out an upper bound of order $\left(\eta ^{*}\right)^2$ with a multiplicative constant depending only on the fixed sub-Gaussian and noise parameters.

Accordingly, in the sub-Gaussian setting we impose a common, known bound
\begin{equation*}
    d\le M,
\end{equation*}
where $M$ is fixed independently of $n, p$ and of the particular class under consideration. Constants may depend on $M$, but remain uniform over the actual diameters $d\le M$. This is stronger than merely requiring each individual class to have finite diameter. Under this assumption, the conversion above is valid with the uniform choice $L=8e^2L_{\mathrm{sg}}^4M^2$.

\subsection{Lower Bound}\label{subsection: lower bound under the sub-gaussian setting}

We now derive the lower bound of the minimax rate of \eqref{eq: minimax risk} relevant to the corruption rate $\epsilon $ in the following lemma, complementing Lemma \ref{lemma: general minimax lower bound} under the sub-Gaussian setting, where the admissible distributions include $\mathcal{N}(0,\sigma ^{2})$

\begin{lemma}[Lower bound in the sub-Gaussian case] \label{lemma: sub-gaussian lower bound on epsilon}
For the nonparametric regression problem formulated in Section \ref{section: introduction}, suppose that the admissible noise distributions include $\mathcal{N}(0, \sigma^2)$. Then
\begin{equation*}
    \mathcal R_{n,\epsilon}(\mathcal F;\mathcal P_\xi^{\mathrm{sg}})\ge \frac{1}{96}\bigl(\sigma^2\epsilon^2\wedge d^2\bigr).
\end{equation*}
This lower bound already holds when the noise distribution is fixed to $\mathcal{N}(0,\sigma^2)$.
\end{lemma}

Similar to Lemma \ref{lemma: lower bound regarding epsilon under the bernstein's condition}, the proof is based on a two-point construction. This time, we select two functions in $\mathcal{F}$ whose $L_2\left(\mathbb{P}_X\right)$-distance is of order $\sigma\epsilon\wedge d$. By Pinsker's inequality and a coupling argument, we argue that the probability that the number of different samples from the two hypotheses does not exceed the corruption budget $\lfloor \epsilon n \rfloor $ is at least $3/4$, in which case the adversary can tamper with the original dataset to make them identical, yielding the results.

Again, combining Lemma \ref{lemma: sub-gaussian lower bound on epsilon} with the general lower bound \eqref{eq: general lower bound} in the uncontaminated case, we conclude that the minimax rate is lower bounded by 
\begin{equation}\label{eq: lower bound under the sub-gaussian condition}
    \mathcal R_{n,\epsilon}(\mathcal F;\mathcal P_\xi^{\mathrm{sg}})\gtrsim \left[(\eta ^{*})^2\vee\epsilon^2\right]\wedge d^2,
\end{equation}
where the implicit constant is independent of $n$, $\epsilon$, and the actual diameter $d$.

Lemma \ref{lemma: sub-gaussian lower bound on epsilon} suggests a minimax rate in the order of $\epsilon ^{2}$. We note here that Example \ref{example: optimality of epsilon under finite second moment}, with a minimax lower bound $\epsilon $, does not meet the sub-Gaussian requirement \eqref{eq: sub-gaussian condition}. Indeed, consider $h _{t}=\mathbf{1}_{\left\{x \in [0,t]\right\}}$ for some $t \in [0,1]$. For $p>2$, we have 
\begin{equation*}
    \left\lVert h _{t}\right\rVert_{L _{p}(\mathbb{P}_{X})}^{}=t ^{\frac{1}{p}}, \left\lVert h _{t}\right\rVert_{L _{2}(\mathbb{P}_{X})}^{}=\sqrt{t}.
\end{equation*}
Therefore the ratio
\begin{equation*}
    \frac{\left\lVert h _{t}\right\rVert_{L _{p}(\mathbb{P}_{X})}^{}}{\sqrt{p}\left\lVert h _{t}\right\rVert_{L _{2}(\mathbb{P}_{X})}^{}}=\frac{t ^{\frac{1}{p}-\frac{1}{2}}}{\sqrt{p}}
\end{equation*}
diverges as $t \rightarrow 0$, which implies that there is no universal $L _{\text{sg}}$ for which \eqref{eq: sub-gaussian condition} holds uniformly. 

\subsection{Upper Bound}\label{subsection: upper bound under the sub-gaussian setting}

The constructive design required for the upper bound is also the tree-based packings given in Algorithm \ref{algorithm: construction of the pruned tree} and the selection criterion given in Algorithm \ref{algorithm: theoretical algorithm}, described in Section \ref{subsection: upper bound under the heavy-tailed setting}. However, as shown in the following lemmas, the sub-Gaussian conditions \eqref{eq: sub-gaussian condition} allow us to deduce stronger properties for the pairwise comparison statistic $Z _{f,g}$, which consequently leads to the wider region of feasibility of Algorithm \ref{algorithm: theoretical algorithm}.

The following lemma indicates that under \eqref{eq: sub-gaussian condition}, $Z _{f,g}$ behaves as a sub-exponential random variable for each fixed pair of $f,g$. We defer the proof to Appendix \ref{subsection: proof of sub-exponential of Z_f,g}. The reader can verify that the marginal second moment control $\mathbb{E}\xi ^{2}\le \sigma ^{2}$ is sufficient in the proof as we leverage H\"{o}lder's inequality.

\begin{lemma}[]\label{lemma: sub-exponential of Z_f,g}
    For the nonparametric regression problem formulated in Section \ref{subsection: problem formulation} and the pairwise comparison statistic $Z _{f,g}$ defined in \eqref{eq: definition of the pairwise comparison statistic}, under Conditions \eqref{eq: sub-gaussian condition} and \eqref{eq: marginal control on the second moment of xi}, and
    \begin{equation*}
        \left\lVert g-f^*\right\rVert_{L_2(\mathbb P_X)}\le\delta,\qquad \left\lVert f-g\right\rVert_{L_2(\mathbb P_X)}\ge C\delta 
    \end{equation*}
    for some $\delta >0, C>3$, we have 
    \begin{equation}
        \left\lVert Z _{f,g}-\mathbb{E}Z _{f,g}\right\rVert_{L _{p}(\mathbb{P}_{\xi ,X})}^{}\le K(\left\lVert f-g\right\rVert_{L _{2}(\mathbb{P}_{X})}^{}+\sigma )p
    \end{equation}
    for any $p \ge 2$, where $K$ is a constant depending on $L _{\text{sg}}, L _{\xi }$, and $C$.
\end{lemma}

Recall that the effectiveness of trimmed mean is fundamental for $Z _{f,g}$ to correctly select between two candidates with one close to $f ^{*}$ while the other is sufficiently far. As remarked in \cite{10.1214/20-AOS1961}, their result (Lemma \ref{lemma: trimmed mean}) can be reinforced when the underlying distribution satisfies stronger moment control. The following lemma establishes one such reinforcement under the sub-Gaussian conditions \eqref{eq: sub-gaussian condition}.

\begin{lemma}[]\label{lemma: validity of the robust trimmed mean test under sub-gaussian condition}
    For two functions $f,g \in \mathcal{F}$, assume that $g$ is close to $f ^{*}$ while $f$ is far from $f ^{*}$ with respect to $\delta >0$:
    \begin{equation*}
        \left\lVert g-f ^{*}\right\rVert_{L _{2}(\mathbb{P}_{X})}\le \delta ,\qquad \left\lVert f-g\right\rVert_{L _{2}(\mathbb{P}_{X})}\ge C\delta .
    \end{equation*}
    Suppose additionally we have \eqref{eq: sub-gaussian condition}, \eqref{eq: marginal control on the second moment of xi}, and $\epsilon \le \min\limits \left\{1/128,\epsilon _{0}\right\}$ for some $\epsilon _{0}$ depending only on $L _{\text{sg}},L _{\xi },C$ (will be determined later). Then for $t _{\delta }:=a \min\limits \left\{\delta ^{2}/\sigma ^{2},1\right\}$, $a \le \min\limits \left\{\frac{1}{384},\frac{1}{1024K ^{2}e ^{2}}\right\}, K=8L _{\text{sg}}(L _{\text{sg}}+L _{\xi })$, as long as $\delta \ge \frac{3456Ke}{C}\sigma \epsilon \log \left(\frac{e}{\epsilon }\right)$, the test $\psi _{f,g}(X,Y)$ from \eqref{eq: pairwise comparison test} can correctly select $g$ over $f$ with probability higher than $1-4\exp\left\{-nt _{\delta }\right\}$.
\end{lemma}

The proof of Lemma \ref{lemma: validity of the robust trimmed mean test under sub-gaussian condition} relies on sharply bounding $\mathcal{E}(4\zeta ,Z _{f,g})$ upon \eqref{eq: sub-gaussian condition}, using the $L _{p}(\mathbb{P}_{X,\xi })$-norm control derived in Lemma \ref{lemma: sub-exponential of Z_f,g}. We defer the proof to Appendix \ref{subsection: proof of validity of the robust trimmed mean test under sub-gaussian condition} for the reader's reference.

Analogous to Corollary \ref{corollary: concrete application of the trimmed mean to bound the error probability}, we conclude the error probability of one-step incorrect selection along the trajectory of the tree traversal. The following argument is straightforward since any nonvacuous comparison must have $\delta \le d \le M$ and consequently we have $t _{\delta }=a \min\limits \left\{\delta ^{2}/\sigma ^{2},1\right\}\ge a\delta ^{2}/(M ^{2}+\sigma ^{2})$.

\begin{corollary}[]\label{corollary: all-scale delta version of the validity of the robust trimmed mean test}
    Under the same conditions and assumptions as Lemma \ref{lemma: validity of the robust trimmed mean test under sub-gaussian condition}, we have 
    \begin{equation}\label{eq: error probability of selection under the sub-gaussian case}
        \sup\limits _{f ^{*}\in \mathcal{F}}\sup\limits _{\left\lVert g-f ^{*}\right\rVert_{L _{2}(\mathbb{P}_{X})}^{}\le \delta ,\left\lVert f-g\right\rVert_{L _{2}(\mathbb{P}_{X})}^{}\ge C\delta }\mathbb{P}\left(\psi _{f,g}(X,Y)=0\right)\le 4 \exp\left\{-\frac{an\delta ^{2}}{M ^{2}+\sigma ^{2}}\right\}.
    \end{equation}
    Equivalently, the probability of incorrect selection scales with $n\delta ^{2}$.  When $4 \exp\left\{-an\delta ^{2}/(M ^{2}+\sigma ^{2})\right\}\ge 1$, this inequality trivially holds.
\end{corollary}

We comment here that once we obtain Lemma \ref{lemma: validity of the robust trimmed mean test under sub-gaussian condition} and Corollary \ref{corollary: all-scale delta version of the validity of the robust trimmed mean test}, the rest of the argument parallels that derived in Section \ref{section: minimax rate under heavy-tailed assumptions}. Specifically, an alternative version of Lemma \ref{lemma: crucial lemma for correct selection}, Corollary \ref{corollary: probability of error during the trajectory of the tree traversal}, Theorem \ref{theorem: main theorem of the algorithm}, and Corollary \ref{corollary: minimax rate under the bernstein condition} and their proofs continue to hold under the sub-Gaussian case. The only, but also the crucial difference is that the feasible region of Lemma \ref{lemma: validity of the robust trimmed mean test under sub-gaussian condition} is now $\delta \ge \frac{3456K}{C}\sigma \epsilon \log \left(\frac{e}{\epsilon }\right)$ compared with $\delta ^{2}\gtrsim _{L,\sigma }\epsilon $ required in Corollary \ref{corollary: concrete application of the trimmed mean to bound the error probability}. Such difference directly changes the performance of Algorithm \ref{algorithm: theoretical algorithm}, where the term regarding $\epsilon $ changes from $\epsilon $ to $\epsilon ^{2}\log ^{2}\left(\frac{e}{\epsilon }\right)$. Formally, we have the following Theorem for the sub-Gaussian case.

\begin{theorem}[]
    Define $J ^{*}$ in the same way as in Theorem \ref{theorem: main theorem of the algorithm}, where $\delta _{j}$, $\eta _{j}$ are similarly defined as: 
    \begin{equation*}
        \delta _{j}=\frac{d}{2 ^{j-1}(C+1)}, \qquad \eta _{j}:=\sqrt{\frac{a}{M ^{2}+\sigma ^{2}}}\delta _{j}, \qquad a \le \min\limits \left\{\frac{1}{384}, \frac{1}{1024K ^{2}e ^{2}}\right\}.
    \end{equation*}
    For the nonparametric regression problem defined in Section \ref{subsection: problem formulation}, suppose that the function class $\mathcal{F}$ and the noise $\xi $ satisfy the sub-Gaussian conditions \eqref{eq: sub-gaussian condition} and the marginal control \eqref{eq: marginal control on the second moment of xi}. For the output of Algorithm \ref{algorithm: theoretical algorithm} $f _{J ^{*}}^{*}$ after $J ^{*}-1$ iterations, we have 
    \begin{equation}\label{eq: upper bound of the theoretical algorithm under the sub-gaussian condition}
    \sup\limits _{f ^{*}\in \mathcal{F}}\sup\limits _{\xi :\mathbb{P}_{\xi }\in \mathcal{P}_{\xi }^{\text{sg}}}\sup\limits _{\mathcal{C}\in \mathfrak{C}_{n,\epsilon }}\mathbb{E}_{\tilde{X},\tilde{Y}}\left\lVert f _{J ^{*}}^{*}(\mathcal{C}(\tilde{X},\tilde{Y}))-f ^{*}\right\rVert_{L _{2}(\mathbb{P}_{X})}^{2}\lesssim \left(\eta _{J ^{*}}^{2}\vee \epsilon ^{2}\log ^{2}\left(\frac{e}{\epsilon }\right)\right)\wedge d ^{2}.
\end{equation}

Combining with \eqref{eq: lower bound under the sub-gaussian condition}, we conclude that the minimax rate of \eqref{eq: minimax risk} is characterized by \eqref{eq: lower bound under the sub-gaussian condition} and \eqref{eq: upper bound of the theoretical algorithm under the sub-gaussian condition} up to a logarithmic factor $\log ^{2}\left(\frac{e}{\epsilon }\right)$.
\end{theorem}

\section{Examples}\label{section: examples}

\subsection{Sparse linear regression}\label{subsection: examples sparse linear regression}

Sparse linear regression provides a classical example of a star-shaped function class that need not be convex. Its uncontaminated minimax rates have been studied extensively; see, for example, \cite{raskutti2011}. Robust estimation under Huber's contamination model is considered in \cite{gao2020}. We show below how our general results specialize to sparse linear regression with Gaussian design under strong adversarial contamination.

Let $\tilde{X}\sim \mathcal{N}(0,\mathbf{I}_{p})$ be a Gaussian random vector (when $\tilde{X}$ follows a sub-Gaussian distribution with zero mean and isotropic covariance matrix, the following arguments hold with minor modification), and let $s$, with $1\le s\le p$, denote the sparsity level of the model, and fix $0<R\le M/2$. Consider the parameter set and associated function class
\begin{equation*}
    \Theta_{s,R}:=\left\{\beta\in\mathbb{R}^{p}:\left\lVert\beta\right\rVert_{0}\le s,\ \left\lVert\beta\right\rVert_{2}\le R\right\},\qquad \mathcal{F}_{s,R}:=\left\{f_{\beta}:f_{\beta}(x)=x ^{\top }\beta,\ \beta\in\Theta_{s,R}\right\}.
\end{equation*}
The class $\mathcal{F}_{s,R}$ is star-shaped about the zero function and has diameter $d=2R\le M$. Moreover, for any $\beta,\gamma\in\Theta_{s,R}$, Gaussianity and isotropy of $X$ give
\begin{equation*}
    f_{\beta}(X)-f_{\gamma}(X)\sim\mathcal{N}\left(0,\left\lVert\beta-\gamma\right\rVert_{2}^{2}\right).
\end{equation*}
Therefore, we have a general expression for $L _{q}(\mathbb{P}_{X})$-norm for any integer $q \ge 2$:
\begin{equation*}
    \left\lVert f _{\beta }-f _{\gamma }\right\rVert_{L _{2}(\mathbb{P}_{X})}^{}=\left\lVert \beta -\gamma \right\rVert_{2}^{},\qquad \left\lVert f_\beta-f_\gamma\right\rVert_{L_q(\mathbb P_X)}=\left\lVert\beta-\gamma\right\rVert_2\left\lVert G\right\rVert_q,
\end{equation*}
where $\left\lVert G\right\rVert_{q}^{}=\left(\mathbb{E}\left|G\right| ^{q}\right)^{1/q}, G \sim \mathcal{N}(0,1)$. From integration by parts and H\"older's inequality:
\begin{equation*}
    \mathbb E\left|G\right|^q=(q-1)\mathbb E\left|G\right|^{q-2}\le(q-1)\left(\mathbb E\left|G\right|^q\right)^{(q-2)/q},
\end{equation*}
where we conclude $\left\lVert G\right\rVert_{q}^{2}\le q-1<q (q \ge 2)$. When $q=1$, we directly have $\left\lVert G\right\rVert_{1}^{}\le \left\lVert G\right\rVert_{2}^{}$ Hence, we have $\left\lVert f _{\beta }-f _{\gamma }\right\rVert_{L _{q}(\mathbb{P}_{X})}^{}\le \sqrt{q}\left\lVert f _{\beta }-f _{\gamma }\right\rVert_{L _{2}(\mathbb{P}_{X})}^{}$, which means the first condition in \eqref{eq: sub-gaussian condition} holds with $L _{\text{sg}}=1$. By the moment comparison in Section \ref{section: minimax rate under sub-gaussian assumptions} and the bound $d\le M$, the Bernstein condition \eqref{eq: bernstein condition} also holds with $L=8e^{2}M^{2}$.

It remains to calculate the local entropy. Write $k:=s\log\left(ep/s\right)$. Fix a center $\beta_0\in\Theta_{s,R}$ and a radius $t>0$. For each coordinate set $S\subseteq[p]$ of cardinality $s$, the vectors supported on $S$ and lying within distance $t$ of $\beta_0$ are contained in an $s$-dimensional Euclidean ball of radius $t$, centered at the orthogonal projection of $\beta_0$ onto that coordinate subspace. The volumetric packing bound therefore gives
\begin{equation*}
    \mathcal{M}_{\mathcal{F}_{s,R}}^{\mathrm{loc}}\left(t,c\right)
    \le
    \binom{p}{s}\left(1+2c\right)^{s}.
\end{equation*}
Taking logarithms yields, uniformly over $t>0$,
\begin{equation*}
    \log\mathcal{M}_{\mathcal{F}_{s,R}}^{\mathrm{loc}}\left(t,c\right)
    \le
    s\log\left(\frac{ep}{s}\right)+s\log\left(1+2c\right)
    \le C_c k,
\end{equation*}
where $C_c>0$ depends only on the fixed packing constant $c$.

For the reverse inequality, the standard sparse packing construction, combined with a Euclidean packing on a fixed support when $s>p/2$, provides universal constants $a,b>0$ and a finite set
\begin{equation*}
    \mathcal{V}\subseteq
    \left\{v\in\mathbb{R}^{p}:\left\lVert v\right\rVert_{0}\le s,\ \left\lVert v\right\rVert_{2}\le1\right\}
\end{equation*}
such that
\begin{equation*}
    \log\left|\mathcal{V}\right|\ge ak,
    \qquad
    \left\lVert v-w\right\rVert_{2}\ge b
    \quad\text{for all distinct }v,w\in\mathcal{V};
\end{equation*}
see the sparse packing arguments in \cite{raskutti2011}. For $0<t\le R$, the functions $\left\{f_{tv}:v\in\mathcal{V}\right\}$ belong to $B_2(0,t)\cap\mathcal{F}_{s,R}$ and have pairwise distances at least $bt$. Taking the universal constant $c$ sufficiently large that $1/c\le b$, we obtain
\begin{equation*}
    \log\mathcal{M}_{\mathcal{F}_{s,R}}^{\mathrm{loc}}\left(t,c\right)
    \ge ak,
    \qquad 0<t\le R.
\end{equation*}

These bounds identify the critical radius, including the regime where the diameter determines the risk. Indeed, the upper entropy bound implies $\left(\eta^*\right)^2\le C_c k/n$. Conversely, setting $t_0:=R\wedge\sqrt{ak/n}$ gives $nt_0^{2}\le ak\le \log\mathcal{M}_{\mathcal{F}_{s,R}}^{\mathrm{loc}}\left(t_0,c\right),$ so the definition of $\eta^*$ implies $\eta^*\ge t_0$. Since $d=2R$, it follows that
\begin{equation*}
    \left(\eta^*\right)^2\wedge d^2\asymp\frac{s\log\left(ep/s\right)}{n}\wedge d^2.
\end{equation*}
Therefore, for fixed $\sigma>0$ and the strong contamination model from Section \ref{subsection: problem formulation}, under the conditional moment assumptions \eqref{eq: finite conditional second moment of xi}, we have
\begin{equation*}
    \mathcal{R}_{n,\epsilon}\left(\mathcal{F}_{s,R};\mathcal{P}_{\xi}^{\mathrm{b}}\right)\asymp\left[\frac{s\log\left(ep/s\right)}{n}\vee\epsilon\right]\wedge d^{2}.
\end{equation*}
Under the sub-Gaussian noise condition in \eqref{eq: sub-gaussian condition}, together with $\mathbb{E}\left[\xi\mid X\right]=0$ and $\mathbb{E}\xi^{2}\le\sigma^{2}$, we instead obtain
\begin{equation*}
    \begin{aligned}
        \left[\frac{s\log\left(ep/s\right)}{n}\vee\epsilon^{2}\right]\wedge d^{2}&\lesssim \mathcal{R}_{n,\epsilon}\left(\mathcal{F}_{s,R};\mathcal{P}_{\xi}^{\mathrm{sg}}\right)\\
        &\lesssim \left[\frac{s\log\left(ep/s\right)}{n}\vee\epsilon^{2}\log^{2}\left(\frac{e}{\epsilon}\right)\right]\wedge d^{2}.
    \end{aligned}
\end{equation*}
The sub-Gaussian bounds hold for $0\le\epsilon\le\epsilon_0$, where $\epsilon_0>0$ is a sufficiently small constant as required in Section \ref{section: minimax rate under sub-gaussian assumptions}; the logarithmic term is interpreted as zero when $\epsilon=0$. The implicit constants depend only on $M,\sigma,c$ and, in the sub-Gaussian case, $L_{\xi}$. In particular, they are independent of $n,p,s,\epsilon$ and the actual diameter $d=2R$.

The heavy-tailed conclusion concerns the supremum over the full conditional-moment noise family $\mathcal{P}_{\xi}^{\mathrm{b}}$, which allows dependence between $\xi$ and $X$. In the sub-Gaussian setting, the bounds retain a logarithmic gap in their dependence on $\epsilon$. For comparison, \cite{gao2020} obtains an $\epsilon^{2}$ contamination term for sparse regression with Gaussian errors under Huber's contamination model while we consider a more general class of noise distributions and the strong contamination model.

\subsection{Low-rank trace regression}\label{subsection: examples low-rank trace regression}

Low-rank trace regression is a matrix analogue of linear regression in which the unknown coefficient matrix has a low-rank structure. Classical results on estimation and minimax rates in this model can be found in \cite{rohde2011,candes2011}. We consider Gaussian design matrices, for which the prediction distance coincides with the Frobenius distance between the coefficient matrices.

Let $m_1,m_2\ge 1$, $1\le r\le \min\left\{m_1,m_2\right\}$, and $0<R\le M/2$, where $M$ is a fixed constant. Suppose that the entries of $X\in\mathbb{R}^{m_1\times m_2}$ are independent standard Gaussian random variables, and consider the regression model
\begin{equation*}
    Y=\operatorname{tr}\left(X^{\top}A^{*}\right)+\xi,
    \qquad
    A^{*}\in\Theta_{r,R}:=\left\{A\in\mathbb{R}^{m_1\times m_2}:\operatorname{rank}\left(A\right)\le r,\ \left\lVert A\right\rVert_{\mathrm F}^{}\le R\right\}.
\end{equation*}
Define the corresponding function class and its effective dimension by
\begin{equation*}
    \mathcal{F}_{r,R}:=\left\{f_A:f_A\left(X\right)=\operatorname{tr}\left(X^{\top}A\right),\ A\in\Theta_{r,R}\right\},
    \qquad
    k_r:=r\left(m_1+m_2-r\right).
\end{equation*}
Here $\left\lVert A\right\rVert_{\mathrm F}^{2}=\sum_{j=1}^{m_1}\sum_{k=1}^{m_2}A_{jk}^{2}$. The class $\mathcal{F}_{r,R}$ is star-shaped around the zero function and has diameter $d=2R\le M$, as witnessed by the two rank-one matrices $Re_1e_1^{\top}$ and $-Re_1e_1^{\top}$, with the two coordinate vectors taken in their respective ambient spaces. Moreover, for any $A,B\in\Theta_{r,R}$, independence and standardization of the Gaussian entries give
\begin{equation*}
    \left\lVert f_A-f_B\right\rVert_{L_2\left(\mathbb{P}_X\right)}^{2}=\mathbb{E}\left[\operatorname{tr}\left(X^{\top}\left(A-B\right)\right)^{2}\right]=\left\lVert A-B\right\rVert_{\mathrm F}^{2}.
\end{equation*}
Using the same argument as in the preceding example, $f_A\left(X\right)-f_B\left(X\right)$ is a centered Gaussian random variable with variance $\left\lVert A-B\right\rVert_{\mathrm F}^{2}$. Hence, for every $q\ge 3$,
\begin{equation*}
    \left\lVert f_A-f_B\right\rVert_{L_q\left(\mathbb{P}_X\right)}^{}\le\sqrt{q}\left\lVert A-B\right\rVert_{\mathrm F}^{}=\sqrt{q}\left\lVert f_A-f_B\right\rVert_{L_2\left(\mathbb{P}_X\right)}^{}.
\end{equation*}

We now calculate the local entropy. Fix the constant $c\ge2$ appearing in the definition of $\eta^{*}$, and write
\begin{equation*}
    H_c\left(t\right):=\log\mathcal{M}_{\mathcal{F}_{r,R}}^{\mathrm{loc}}\left(t,c\right).
\end{equation*}
We will show that there are positive constants $a_c,b_c$, depending only on $c$, such that
\begin{equation*}
    H_c\left(t\right)\le b_c k_r
    \quad\text{for every }t>0,
    \qquad
    H_c\left(t\right)\ge a_c k_r
    \quad\text{for }0<t\le R.
\end{equation*}

For the upper bound, let $\mathcal{N}\left(u,\mathcal{A}\right)$ denote the covering number of a matrix set $\mathcal{A}$ by Frobenius balls of radius $u$. For $1\le k\le\min\left\{m_1,m_2\right\}$, define
\begin{equation*}
    \mathcal{S}_k:=\left\{A:\operatorname{rank}\left(A\right)\le k,\ \left\lVert A\right\rVert_{\mathrm F}^{}=1\right\},
    \qquad
    \mathcal{K}_k:=\left\{A:\operatorname{rank}\left(A\right)\le k,\ \left\lVert A\right\rVert_{\mathrm F}^{}\le1\right\}.
\end{equation*}
The low-rank covering bound of \cite{candes2011} states that, for $0<u<1$,
\begin{equation*}
    \mathcal{N}\left(u,\mathcal{S}_k\right)
    \le
    \left(\frac{9}{u}\right)^{\left(m_1+m_2+1\right)k}.
\end{equation*}
Their result concerns the unit sphere. To obtain the corresponding bound for the unit ball, take a $u/2$-net of $\mathcal{S}_k$ and a grid of $[0,1]$ with mesh at most $u/2$ and cardinality at most $4/u$. Every nonzero matrix in $\mathcal{K}_k$ can be written as $aU$, where $a\in[0,1]$ and $U\in\mathcal{S}_k$. Choosing a grid point $b$ and a net point $V$ with $\left|a-b\right|\le u/2$ and $\left\lVert U-V\right\rVert_{\mathrm F}^{}\le u/2$ yields
\begin{equation*}
    \left\lVert aU-bV\right\rVert_{\mathrm F}^{}
    \le
    a\left\lVert U-V\right\rVert_{\mathrm F}^{}
    +
    \left|a-b\right|\left\lVert V\right\rVert_{\mathrm F}^{}
    \le u.
\end{equation*}
Including zero in the radial grid therefore gives
\begin{equation*}
    \mathcal{N}\left(u,\mathcal{K}_k\right)
    \le
    \frac{4}{u}\left(\frac{18}{u}\right)^{\left(m_1+m_2+1\right)k}
    \le
    \left(\frac{18}{u}\right)^{\left(m_1+m_2+1\right)k+1}.
\end{equation*}

Fix any center $A_0\in\Theta_{r,R}$ and radius $t>0$. Differences between two matrices of rank at most $r$ have rank at most $2r$. Thus, setting $r_0:=\min\left\{2r,m_1,m_2\right\}$, we have
\begin{equation*}
    \left\{A-A_0:A\in\Theta_{r,R},\ \left\lVert A-A_0\right\rVert_{\mathrm F}^{}\le t\right\}
    \subseteq t\mathcal{K}_{r_0}.
\end{equation*}
Applying the preceding covering bound with $u=1/\left(3c\right)$ shows that $t\mathcal{K}_{r_0}$ can be covered by at most $\left(54c\right)^{\left(m_1+m_2+1\right)r_0+1}$ Frobenius balls of radius $t/\left(3c\right)$. Each such ball contains at most one point of a $t/c$-separated set, because its diameter is strictly smaller than $t/c$. Taking the supremum over $A_0$ and using the isometry established above yields
\begin{equation*}
    H_c\left(t\right)
    \le
    \left[\left(m_1+m_2+1\right)r_0+1\right]\log\left(54c\right)
    \le
    7k_r\log\left(54c\right).
\end{equation*}
For the last inequality, we used $r_0\le2r$, $m_1+m_2+1\le3\left(m_1+m_2\right)/2$, and $k_r=r\left(m_1+m_2-r\right)\ge\frac{r\left(m_1+m_2\right)}{2}\ge1$.

This proves the upper bound for every $t>0$. In particular, the covering resolution after rescaling is the fixed number $1/\left(3c\right)$, which explains why no factor involving $\log\left(R/t\right)$ appears in the local entropy.

For the lower bound, suppose without loss of generality that $m_1\ge m_2$; otherwise, transpose all matrices. Consider the linear subspace
\begin{equation*}
    \mathcal{V}:=\left\{A=\left[B\ \ 0\right]:B\in\mathbb{R}^{m_1\times r}\right\}.
\end{equation*}
Every matrix in $\mathcal{V}$ has rank at most $r$, and $\mathcal{V}$ has dimension $D=m_1r$. Therefore, for $0<t\le R$, the set
\begin{equation*}
    \left\{A\in\mathcal{V}:\left\lVert A\right\rVert_{\mathrm F}^{}\le t\right\}
\end{equation*}
is contained in the local parameter set centered at zero and is isometric to a Euclidean ball of radius $t$ in $\mathbb{R}^{D}$.

Take a maximal $t/c$-separated subset of this Euclidean ball, and denote its cardinality by $N$. Maximality implies that balls of radius $t/c$ around these points cover the radius-$t$ ball. Comparing $D$-dimensional volumes gives
\begin{equation*}
    N\ge
    \frac{\operatorname{vol}_D\left(t\mathbb{B}_2^{D}\right)}
    {\operatorname{vol}_D\left(\left(t/c\right)\mathbb{B}_2^{D}\right)}
    =c^{D},
\end{equation*}
where $\mathbb{B}_2^{D}$ is the Euclidean unit ball. Since the zero function belongs to $\mathcal{F}_{r,R}$, this packing supplies a lower bound for its local packing number. Consequently,
\begin{equation*}
    H_c\left(t\right)
    \ge
    D\log c
    \ge
    \frac{k_r}{2}\log c,
    \qquad 0<t\le R,
\end{equation*}
where we used $k_r\le r\left(m_1+m_2\right)\le2m_1r=2D$. Combining the two bounds proves
\begin{equation*}
    \log\mathcal{M}_{\mathcal{F}_{r,R}}^{\mathrm{loc}}\left(t,c\right)
    \asymp_c
    r\left(m_1+m_2-r\right),
    \qquad 0<t\le R.
\end{equation*}

It remains to calculate the critical radius. Set $a_c=\left(\log c\right)/2$ and $b_c=7\log\left(54c\right)$. By the definition of $\eta^{*}$ and the entropy upper bound valid at every radius, we have $\left(\eta^{*}\right)^{2}\le\frac{b_c k_r}{n}$.

For the reverse direction, take $t_0:=R\wedge\sqrt{\frac{a_c k_r}{n}}$. Since $0<t_0\le R$, the entropy lower bound gives $nt_0^{2}\le a_c k_r\le H_c\left(t_0\right)$. Thus $t_0$ belongs to the set defining $\eta^{*}$, and $\eta^{*}\ge t_0$. Recalling that $d=2R$, we conclude that
\begin{equation*}
    \left(\eta^{*}\right)^{2}\wedge d^{2}
    \asymp_c
    \frac{r\left(m_1+m_2-r\right)}{n}\wedge d^{2}.
\end{equation*}

Substituting this critical radius into Corollary \ref{corollary: minimax rate under the bernstein condition} proves the rate of \eqref{eq: minimax risk} for the conditional second-moment noise family. Specifically, under the strong corruption model, for sufficiently small $\epsilon$, if the noise satisfies Condition \eqref{eq: finite conditional second moment of xi}, then
\begin{equation*}
    \mathcal{R}_{n,\epsilon}\left(\mathcal{F}_{r,R};\mathcal{P}_{\xi}^{\mathrm b}\right)
    \asymp
    \left[\frac{k_r}{n}\vee\epsilon\right]\wedge d^{2}.
\end{equation*}
Applying the upper and lower bounds from Section \ref{section: minimax rate under sub-gaussian assumptions} controls the rate under the sub-Gaussian condition \eqref{eq: sub-gaussian condition}:
\begin{equation*}
    \left[\frac{k_r}{n}\vee\epsilon^{2}\right]\wedge d^{2}
    \lesssim
    \mathcal{R}_{n,\epsilon}\left(\mathcal{F}_{r,R};\mathcal{P}_{\xi}^{\mathrm{sg}}\right)
    \lesssim
    \left[\frac{k_r}{n}\vee\epsilon^{2}\log^{2}\left(\frac{e}{\epsilon}\right)\right]\wedge d^{2}.
\end{equation*}
The corruption term is interpreted as zero when $\epsilon=0$. The implicit constants may depend on the fixed parameters $M$, $\sigma>0$, and the sub-Gaussian noise constant, but do not depend on $n$, $m_1$, $m_2$, $r$, $R$, $d$, or $\epsilon$.

\subsection{H\"older smooth regression}\label{subsection: examples holder smooth regression}

We next consider regression over a H\"older ball on a fixed-dimensional domain. Throughout this example, the smoothness parameter $\alpha>0$ and the dimension $m\ge1$ are fixed. Constants may depend on $\alpha$ and $m$, but remain independent of the sample size, corruption proportion, and actual diameter of the function class.

Let $k:=\lceil\alpha\rceil-1$ and $\beta:=\alpha-k\in\left(0,1\right]$. For a function with continuous partial derivatives through order $k$ on $\left[0,1\right]^m$, define
\begin{equation*}
    \left\lVert f\right\rVert_{C^\alpha}^{}
    :=
    \max_{\left|\nu\right|\le k}\left\lVert D^\nu f\right\rVert_{\infty}^{}
    +
    \max_{\left|\nu\right|=k}
    \sup_{\substack{x,y\in\left[0,1\right]^m\\x\ne y}}
    \frac{\left|D^\nu f\left(x\right)-D^\nu f\left(y\right)\right|}
    {\left\lVert x-y\right\rVert_2^{\beta}},
\end{equation*}
where $\nu$ is a multi-index and $\left|\nu\right|$ is its total order. This convention includes integer smoothness parameters by requiring the derivatives of order $\alpha-1$ to be Lipschitz continuous.

Suppose that $X$ is uniformly distributed on $\left[0,1\right]^m$ and
\begin{equation*}
    Y=f^{*}\left(X\right)+\xi,
    \qquad
    f^{*}\in\mathcal{F}_{\alpha,R}
    :=
    \left\{f:\left\lVert f\right\rVert_{C^\alpha}^{}\le R\right\},
    \qquad
    0<R\le M/2,
\end{equation*}
where $M$ is a fixed constant. Assume that the noise belongs to the family $\mathcal{P}_{\xi}^{\mathrm b}$ satisfying
\begin{equation*}
    \mathbb{E}\left[\xi\mid X\right]=0,
    \qquad
    \mathbb{E}\left[\xi^{2}\mid X\right]\le\sigma^{2}.
\end{equation*}

The class $\mathcal{F}_{\alpha,R}$ is a norm ball and contains the zero function, so it is convex and star-shaped. Since $\left\lVert f\right\rVert_{\infty}^{}\le R$ for every $f\in\mathcal{F}_{\alpha,R}$, its $L_2\left(\mathbb{P}_X\right)$-diameter is at most $2R$. Equality follows by considering the two constant functions $f\equiv R$ and $g\equiv-R$.

Moreover, for every integer $q\ge 2$ and $f,g\in\mathcal{F}_{\alpha,R}$,
\begin{equation*}
    \begin{aligned}
        \left\lVert f-g\right\rVert_{L_{2q}\left(\mathbb{P}_X\right)}^{2q}&\le\left\lVert f-g\right\rVert_{\infty}^{2q-2}\left\lVert f-g\right\rVert_{L_2\left(\mathbb{P}_X\right)}^{2}\\
        &\le M^{2q-2}\left\lVert f-g\right\rVert_{L_2\left(\mathbb{P}_X\right)}^{2}\\
        &\le\frac{q!}{2}\left(M^{2}\right)^{q-1}\left\lVert f-g\right\rVert_{L_2\left(\mathbb{P}_X\right)}^{2}.
    \end{aligned}
\end{equation*}
Thus Condition \eqref{eq: bernstein condition} holds with $L=M^{2}$, independently of the actual diameter $d$.

We next calculate the local entropy. Fix the sufficiently large constant $c$ used in the definition of $\eta^{*}$, taking $c\ge4$, and write
\begin{equation*}
    H_c\left(t\right)
    :=
    \log\mathcal{M}_{\mathcal{F}_{\alpha,R}}^{\mathrm{loc}}\left(t,c\right).
\end{equation*}
We will establish
\begin{equation*}
    H_c\left(t\right)
    \asymp_{\alpha,m,c}
    \left(\frac{R}{t}\right)^{m/\alpha},
    \qquad 0<t\le R.
\end{equation*}

For the upper bound, let $\mathcal{N}\left(u,\mathcal{F},\left\lVert\cdot\right\rVert_{\infty}^{}\right)$ denote the covering number in the uniform norm. Theorem~2.7.1 of \cite{vandervaart1996}, applied to the unit H\"older ball on $\left[0,1\right]^m$, gives
\begin{equation*}
    \log\mathcal{N}\left(u,\mathcal{F}_{\alpha,1},\left\lVert\cdot\right\rVert_{\infty}^{}\right)
    \le
    A_{\alpha,m}u^{-m/\alpha},
    \qquad 0<u\le1,
\end{equation*}
where $A_{\alpha,m}$ depends only on $\alpha$ and $m$. Scaling the function class by $R$ yields
\begin{equation*}
    \log\mathcal{N}\left(u,\mathcal{F}_{\alpha,R},\left\lVert\cdot\right\rVert_{\infty}^{}\right)
    \le
    A_{\alpha,m}\left(\frac{R}{u}\right)^{m/\alpha},
    \qquad 0<u\le R.
\end{equation*}
For $u\ge R$, the class is covered by the single uniform-norm ball centered at zero with radius $u$, so the same upper bound remains valid.

Because $\mathbb{P}_X$ is a probability measure, $\left\lVert f-g\right\rVert_{L_2\left(\mathbb{P}_X\right)}^{}\le\left\lVert f-g\right\rVert_{\infty}^{}$. Therefore, a uniform-norm ball of radius $t/\left(3c\right)$ contains at most one point of a $t/c$-separated set in $L_2\left(\mathbb{P}_X\right)$. Every local packing is also a packing of the entire class, and hence
\begin{equation*}
    \begin{aligned}
        H_c\left(t\right)
        &\le
        \log\mathcal{N}\left(\frac{t}{3c},\mathcal{F}_{\alpha,R},\left\lVert\cdot\right\rVert_{\infty}^{}\right)\\
        &\le
        A_{\alpha,m}\left(\frac{3cR}{t}\right)^{m/\alpha}
        =
        b_{\alpha,m,c}\left(\frac{R}{t}\right)^{m/\alpha},
        \qquad t>0.
    \end{aligned}
\end{equation*}

For the lower bound, we construct a packing inside the single local ball centered at zero. Choose a smooth function $\psi$ compactly supported in $\left(1/4,3/4\right)^m$, normalized so that
\begin{equation*}
    \int_{\mathbb{R}^{m}}\psi\left(x\right)^{2}\,dx=1.
\end{equation*}
For an integer $K\ge1$, define the index set $\mathcal{J}_K:=\left\{0,\ldots,K-1\right\}^{m}$ and let $N_K:=K^{m}$. For each binary vector $\theta=\left(\theta_j\right)_{j\in\mathcal{J}_K}\in\left\{0,1\right\}^{N_K}$, set
\begin{equation*}
    g_{\theta}\left(x\right)
    :=
    t\sum_{j\in\mathcal{J}_K}\theta_j\psi\left(Kx-j\right).
\end{equation*}
The bumps have disjoint supports, each lying strictly inside one of the grid cubes of side length $K^{-1}$. Differentiation contributes a factor $K^{\left|\nu\right|}$, and the $\beta$-H\"older seminorm of a derivative of order $k$ contributes an additional factor $K^\beta$. The separation between supports controls differences involving distinct bumps. Consequently, there exists a constant $C_0\ge1$, depending only on $\alpha,m$ and the fixed bump $\psi$, such that $\left\lVert g_{\theta}\right\rVert_{C^\alpha}^{} \le C_0tK^\alpha$ uniformly over $K$ and $\theta$. Since $\psi$ is fixed once $\alpha,m$ are fixed, we regard $C_0$ as depending only on $\alpha,m$.

Disjointness of the supports and the normalization of $\psi$ also give the exact identities
\begin{equation*}
    \left\lVert g_{\theta}\right\rVert_{L_2\left(\mathbb{P}_X\right)}^{2}=t^{2}K^{-m}\sum_{j\in\mathcal{J}_K}\theta_j^{2}\le t^{2},
\end{equation*}
and
\begin{equation*}
    \left\lVert g_{\theta}-g_{\theta'}\right\rVert_{L_2\left(\mathbb{P}_X\right)}^{2}=t^{2}K^{-m}\rho_{\mathrm H}\left(\theta,\theta'\right),
\end{equation*}
where $\rho_{\mathrm H}$ denotes Hamming distance.

First suppose that $0<t\le R/C_0$, and choose $K:=\left\lfloor\left(\frac{R}{C_0t}\right)^{1/\alpha}\right\rfloor$. Then $K\ge1$ and $C_0tK^\alpha\le R$, so every $g_{\theta}$ belongs to $\mathcal{F}_{\alpha,R}\cap B_2\left(0,t\right)$. The Varshamov--Gilbert bound, as stated in Lemma~2.9 of \cite{tsybakov2009}, supplies a subset $\Theta_K\subseteq\left\{0,1\right\}^{N_K}$ satisfying
\begin{equation*}
    \log\left|\Theta_K\right|
    \ge \frac{\log2}{8}N_K,
    \qquad
    \rho_{\mathrm H}\left(\theta,\theta'\right)\ge\frac{N_K}{8}
    \quad\text{for distinct }\theta,\theta'\in\Theta_K.
\end{equation*}
The cited lemma applies when $N_K\ge8$; when $N_K<8$, the all-zero and all-one vectors give these same inequalities. It follows that $\left\lVert g_{\theta}-g_{\theta'}\right\rVert_{L_2\left(\mathbb{P}_X\right)}^{}\ge\frac{t}{\sqrt{8}}>\frac{t}{c}$. Thus these functions form an admissible local packing. Using $\lfloor x\rfloor\ge x/2$ for $x\ge1$, we obtain
\begin{equation*}
    \begin{aligned}
        H_c\left(t\right)
        &\ge
        \frac{\log2}{8}K^{m}\\
        &\ge
        \frac{\log2}{8\cdot2^{m}C_0^{m/\alpha}}
        \left(\frac{R}{t}\right)^{m/\alpha},
        \qquad 0<t\le R/C_0.
    \end{aligned}
\end{equation*}

For the remaining radii $R/C_0<t\le R$, the constant functions $f\equiv t$ and $g\equiv-t$ both belong to $\mathcal{F}_{\alpha,R}\cap B_2\left(0,t\right)$ and are $2t$ apart. Hence $H_c\left(t\right)\ge\log2$. Since $\left(R/t\right)^{m/\alpha}\le C_0^{m/\alpha}$ in this range, this also gives a lower bound of the required order. Combining the two ranges, we have
\begin{equation*}
    a_{\alpha,m}\left(\frac{R}{t}\right)^{m/\alpha}
    \le H_c\left(t\right)
    \le b_{\alpha,m,c}\left(\frac{R}{t}\right)^{m/\alpha},
    \qquad 0<t\le R,
\end{equation*}
for positive constants independent of $R$ and $t$. The upper bound continues to hold for every $t>0$.

We now calculate the critical radius. Define $v_n:=R^{\frac{2m}{2\alpha+m}}n^{-\frac{2\alpha}{2\alpha+m}}$. Whenever $nt^{2}\le H_c\left(t\right)$, the entropy upper bound implies $nt^{2+m/\alpha}\le b_{\alpha,m,c}R^{m/\alpha}$,
and therefore $\left(\eta^{*}\right)^{2}\le b_{\alpha,m,c}^{\frac{2\alpha}{2\alpha+m}}v_n$. Conversely, set $t_0:=R\wedge\left(a_{\alpha,m}R^{m/\alpha}/n\right)^{\frac{\alpha}{2\alpha+m}}$.
Then $0<t_0\le R$ and
\begin{equation*}
    nt_0^{2}
    \le
    a_{\alpha,m}\left(\frac{R}{t_0}\right)^{m/\alpha}
    \le H_c\left(t_0\right).
\end{equation*}
Thus $\eta^{*}\ge t_0$. Recalling $d=2R$, we conclude that
\begin{equation*}
    \left(\eta^{*}\right)^{2}\wedge d^{2}
    \asymp_{\alpha,m,c}
    R^{\frac{2m}{2\alpha+m}}n^{-\frac{2\alpha}{2\alpha+m}}
    \wedge d^{2}.
\end{equation*}
Finally, by Corollary \ref{corollary: minimax rate under the bernstein condition}, under the strong corruption model, for sufficiently small $\epsilon$, we have:
\begin{equation*}
    \mathcal{R}_{n,\epsilon}\left(\mathcal{F}_{\alpha,R};\mathcal{P}_{\xi}^{\mathrm b}\right)
    \asymp_{\alpha,m,M,\sigma}
    \left[
        R^{\frac{2m}{2\alpha+m}}
        n^{-\frac{2\alpha}{2\alpha+m}}
        \vee\epsilon
    \right]\wedge d^{2}.
\end{equation*}
In particular, for fixed $R>0$, the uncontaminated component has the classical nonparametric rate $n^{-2\alpha/\left(2\alpha+m\right)}$.
The asserted minimax rate now follows from Corollary \ref{corollary: minimax rate under the bernstein condition}.

\begin{remark*}[]
    Although the functions in $\mathcal{F}_{\alpha,R}$ are uniformly bounded, the class does not satisfy the relative sub-Gaussian increment condition \eqref{eq: sub-gaussian condition}. Indeed, using the bump $\psi$ in the argument above, let
    \begin{equation*}
        h_K\left(x\right):=\frac{R}{C_0}K^{-\alpha}\psi\left(Kx\right),
        \qquad K\ge1.
    \end{equation*}
    The same scaling argument shows that $h_K\in\mathcal{F}_{\alpha,R}$, whereas
    \begin{equation*}
        \frac{\left\lVert h_K\right\rVert_{L_4\left(\mathbb{P}_X\right)}^{}}
        {\left\lVert h_K\right\rVert_{L_2\left(\mathbb{P}_X\right)}^{}}
        =
        K^{m/4}
        \frac{\left\lVert\psi\right\rVert_{L_4\left(\mathbb{R}^{m}\right)}^{}}
        {\left\lVert\psi\right\rVert_{L_2\left(\mathbb{R}^{m}\right)}^{}}
        \longrightarrow\infty.
    \end{equation*}
    Since the zero function also belongs to the class, no finite $L_{\mathrm{sg}}$ can make the required increment inequality hold uniformly, even at $q=4$. Consequently, the sub-Gaussian minimax theorem cannot be applied to this full H\"older ball merely by imposing sub-Gaussian noise.
\end{remark*}

\begin{remark*}[]
    Section 3.1 of \cite{gao2020} studies a related regression problem with uniform design on $\left[0,1\right]$, independent standard Gaussian noise, and a Sobolev ball defined through the Fourier coefficients of the regression function. Under Huber's $\epsilon$-contamination model, Theorem 3.1 therein reports a high-probability bound of order $n^{-2\alpha/(2\alpha+1)}\vee\epsilon^2$ for the integrated squared error of a regression depth estimator. For fixed $R>0$, the uncontaminated term for our H\"older ball is of order $n^{-2\alpha/(2\alpha+m)}$, so the two exponents agree when $m=1$. The contamination terms differ: in our setting, the minimax risk has a term of order $\epsilon$, subject to the diameter cap $d^2$. This linear dependence is a worst-case statement over conditionally centered noise laws satisfying $\mathbb{E}\left[\xi^2\mid X\right]\le\sigma^2$ and adversaries from the strong contamination model. It does not assert that the same linear dependence is necessary when the noise law is fixed to be Gaussian. We also acknowledge one limitation of our approach that we are not able to obtain the $\epsilon ^{2}$ term on the contamination even when the noise is Gaussian, as what is achieved by \cite{gao2020}, since by the last remark the function class $\mathcal{F}_{\alpha ,R}$ in general does not satisfy the sub-Gaussian increment conditions \eqref{eq: sub-gaussian condition}. 
\end{remark*}

\section{Discussion}\label{section: discussion}

\subsection{Summary}\label{subsection: summary}

We study nonparametric regression over star-shaped function classes under strong adversarial contamination, allowing the adversary to replace both covariates and responses. Our results relate the minimax squared $L_{2}(\mathbb{P}_{X})$ risk to the local metric entropy of the function class and the contamination proportion. Under the Bernstein-type increment condition and a bounded conditional second moment of the noise, we establish matching upper and lower bounds of order $\left[(\eta^{*})^{2}\vee\epsilon\right]\wedge d^{2}$. Under relative sub-Gaussian increment and noise conditions, together with a common diameter bound, the contamination contribution lies between $\epsilon^{2}$ and $\epsilon^{2}\log^{2}(e/\epsilon)$. The upper bounds follow from a hierarchy of local packings combined with robust pairwise comparisons based on modified trimmed means. The examples of sparse linear regression, low-rank trace regression, and H\"older smooth regression illustrate how the local-entropy characterization yields explicit rates. Our analysis also distinguishes the roles of conditional and marginal noise assumptions: conditional second-moment control is essential for the general guarantee under the Bernstein-type condition, whereas marginal second-moment control suffices under the stronger sub-Gaussian assumptions.

\subsection{Future Research Directions}\label{subsection: future research directions}

An immediate question is the exact dependence on $\epsilon$ in the sub-Gaussian setting. The present upper bound uses sub-exponential moment control of normalized squared-loss differences, while the contamination lower bound is obtained from a Gaussian noise submodel. These arguments leave open whether either of the logarithmic factors can be removed. Resolving this question may require sharper comparisons that exploit conditional centering more directly, or stronger lower bounds using the full family of admissible sub-Gaussian noise distributions. The total-variation modulus approach developed by \citet{chen2016} provides a useful perspective on how contamination sensitivity depends on the statistical model. Applying this perspective here would require accounting for the fixed design distribution, the restrictions on the conditional noise laws, and the finite replacement budget of the strong adversary.

A second direction is computational efficiency. The packing construction establishes statistical achievability, but constructing the candidate sets and performing all pairwise comparisons can be computationally prohibitive. It is therefore natural to seek efficient procedures for structured subclasses, including the sparse and low-rank models considered in Section \ref{section: examples}. Robust iterative hard thresholding in \citet{liu2020} and covariate filtering in \citet{pensia2025} demonstrate how additional structure can support tractable robust regression methods. \cite{li2026minimaxityefficiencyexponentialfamily} study estimation of the natural parameter in uncontaminated exponential-family models under star-shaped constraints and develop a computationally efficient estimator under additional geometric assumptions. In the related robust signal-detection problems, \cite{li2026robustsignaldetectionquadratically,li2026efficientrobustconstrainedsignal} develop efficient robust testing procedures with matching minimax upper bound up to a logarithmic factor. A corresponding goal in our setting is to replace explicit packing searches by suitable optimization procedures while retaining the localization properties and comparison guarantees used in the analysis. Identifying the geometric or computational assumptions needed for such a replacement would help distinguish statistical achievability from efficient achievability.

Another direction concerns unknown design geometry and tuning parameters. Our construction assumes access to the $L_{2}(\mathbb{P}_{X})$ geometry, or to an equivalent norm with controlled comparison constants. When such an oracle is unavailable, the geometry must be estimated from covariates that may themselves be corrupted. Developing a joint analysis of robust distance estimation and function selection would extend the scope of the present results. It is also of interest to construct procedures that adapt to an unknown contamination proportion or unknown moment bounds. Such adaptation requires calibrating the pairwise comparisons without knowing the scale at which contamination prevents further reliable refinement. An important question is whether this can be achieved without changing the local-entropy rate or introducing additional logarithmic losses.

Finally, the two distributional regimes considered here suggest studying intermediate moment assumptions and more restricted noise families. Bounds on conditional moments of order greater than two may yield contamination penalties between the finite-variance and sub-Gaussian regimes. Requiring independence between the noise and the covariates, or restricting contamination to the responses, may also change the least favorable distributions and the attainable rates. Characterizing these changes would clarify which features of the contamination penalty are determined by the geometry of the function class and which arise from the permitted dependence and tail behavior of the noise.

\section*{AI Usage Disclosure}

Throughout this study and the write-up of this manuscript, ChatGPT 6 and 6.1 were used mainly to find and correct the grammar mistakes, improve the word usage, and polish the language. They were also leveraged to search and provide some insights on the studies that share similar background with the current work where the authors are unaware of. Some constants in proofs of some lemmas, theorems, and corollaries are also validated by AI for correctness. The authors are responsible for the novelty of the idea of this work, the main framework of, and the derivations of the lemmas, theorems, corollaries, and algorithms in the manuscript.

\bibliographystyle{abbrvnat}
\bibliography{biblio}
\newpage

\appendix

\section{Proof of Section \ref{section: general lower bound}}\label{section: proof of lower bound of minimax rate}

\subsection{Proof of Lemma \ref{lemma: control on mutual information}}\label{subsection: proof of control on mutual information}

\begin{proof}
    Denote $\mathbb{P}_{\tilde{X};j}, \mathbb{P}_{\tilde{Y}\left\lvert\right. \tilde{X};j}, \mathbb{P}_{\tilde{X},\tilde{Y};j}$ as the law of $\tilde{X}$, $\tilde{Y}$ conditioned on $\tilde{X}$, and $(\tilde{X},\tilde{Y})$ when $J=j$ (i.e., $f ^{*}=f _{j}$). Let $\bar{\mathbb{P}}_{\tilde{X},\tilde{Y}}:=\frac{1}{m}\sum\limits _{j=1}^{m}\mathbb{P}_{\tilde{X},\tilde{Y};j}$ be the mixture distribution of $(\tilde{X},\tilde{Y})$ over the index $J$. From the definition above, we know $\mathbb{P}_{\tilde{X},\tilde{Y};j}=\mathbb{P}_{\tilde{X};j}\times \mathbb{P}_{\tilde{Y}\left\lvert\right. \tilde{X};j}$. By the definition of mutual information, we have 
    \begin{equation*}
            I(\tilde{X},\tilde{Y};J)=\frac{1}{m}\sum\limits _{j=1}^{m}\operatorname{KL}(\mathbb{P}_{\tilde{X},\tilde{Y};j}\left\lVert\right. \bar{\mathbb{P}}_{\tilde{X},\tilde{Y}})\overset{\text{(i)}}{\le }\frac{1}{m}\sum\limits _{j=1}^{m}\operatorname{KL}(\mathbb{P}_{\tilde{X},\tilde{Y};j}\left\lVert\right. \mathbb{P}^{h}_{\tilde{X},\tilde{Y}}).
    \end{equation*}
    Here $\mathbb{P}^{h}_{\tilde{X},\tilde{Y}}$ represents the joint probability distribution of $(\tilde{X},\tilde{Y})$ when $f ^{*}=h \in \mathcal{F}$, and (i) leverages the property that $\bar{\mathbb{P}}_{\tilde{X},\tilde{Y}}$ minimizes the objective $\frac{1}{m}\sum\limits _{j=1}^{m}\operatorname{KL}(\mathbb{P}_{\tilde{X},\tilde{Y};j}\left\lVert\right. \mathbb{Q})$ over arbitrary probability distribution $\mathbb{Q}$. Now by the definition of KL divergence, we further have 
    \begin{equation*}
        \begin{aligned}
            I(\tilde{X},\tilde{Y};J)&\le \frac{1}{m}\sum\limits _{j=1}^{m}\mathbb{E}_{\tilde{X},\tilde{Y};j}\left(\log \frac{d \mathbb{P}_{\tilde{X},\tilde{Y};j}}{d\mathbb{P}^{h}_{\tilde{X},\tilde{Y}}}\right)\\ 
            &=\frac{1}{m}\sum\limits _{j=1}^{m}\mathbb{E}_{\tilde{X}}\mathbb{E}_{\tilde{Y}\left\lvert\right. \tilde{X};j}\left(\log \frac{d \mathbb{P}_{\tilde{Y}\left\lvert\right. \tilde{X};j}}{d\mathbb{P}^{h}_{\tilde{Y}\left\lvert\right. \tilde{X}}}\right)\\ 
            &=\frac{1}{m}\sum\limits _{j=1}^{m}\sum\limits _{i=1}^{n}\mathbb{E}_{\tilde{X}_{i}}\mathbb{E}_{\tilde{Y}_{i}\left\lvert\right. \tilde{X}_{i};j}\left(\log \frac{d \mathbb{P}_{\tilde{Y}_{i}\left\lvert\right. \tilde{X}_{i};j}}{d\mathbb{P}^{h}_{\tilde{Y}_{i}\left\lvert\right. \tilde{X}_{i}}}\right)\\ 
            &\overset{\text{(i)}}{\le }\frac{1}{2m \sigma ^{2}}\sum\limits _{j=1}^{m}\sum\limits _{i=1}^{n}\left\lVert f _{j}-h\right\rVert_{L _{2}(\mathbb{P}_{X})}^{2}\\ 
            &\le \frac{n}{2 \sigma ^{2}}\max\limits _{j}\left\lVert f _{j}-h\right\rVert_{L _{2}(\mathbb{P}_{X})}^{2},
        \end{aligned}
    \end{equation*}
    where (i) uses the property that given $\tilde{X}_{i}$, $\tilde{Y}_{i}$ is distributed as $\mathcal{N}(f ^{*}(\tilde{X}_{i}),\sigma ^{2})$ and the KL divergence between the distributions $\mathcal{N}\left(f _{j}(\tilde{X}_{i}),\sigma ^{2}\right)$ and $\mathcal{N}\left(h(\tilde{X}_{i}),\sigma ^{2}\right)$ is $\frac{\left[f ^{*}(\tilde{X}_{i})-h(\tilde{X}_{i})\right]^{2}}{2 \sigma ^{2}}$, followed by the definition of $L _{2}(\mathbb{P}_{X})$-distance.
\end{proof}

\subsection{Proof of Lemma \ref{lemma: general minimax lower bound}}\label{subsection: proof of general minimax lower bound}

\begin{proof}
    For any $\eta $ satisfying $\log \mathcal{M}_{\mathcal{F}}^{\text{loc}}(\eta ,c)\ge 4 \left(\frac{n \eta ^{2}}{2 \sigma ^{2}}\vee \log 2\right)$ (in the case that $\log \mathcal{M}_{\mathcal{F}}^{\text{loc}}(\eta ,c)$ is infinite, choose any packing set that is finite and satisfies the inequality), choose some $f \in \mathcal{F}$ admitting a $\eta /c$-packing set $\left\{f _{1},\dots,f _{\mathcal{M}(\eta /c,B _{2}(f,\eta )\cap \mathcal{F})}\right\}$ of $B _{2}(f, \eta )\cap \mathcal{F}$. (If the local packing number is infinite, choose any finite subset of the packing set whose cardinality satisfies the condition above.) Therefore, from Theorem \ref{theorem: fano inequality} and Lemma \ref{lemma: control on mutual information}, we know that 
    \begin{equation*}
        \mathcal{R}_{n,\epsilon}(\mathcal{F};\mathcal{P}_{\xi})\ge \frac{\eta ^{2}}{4c ^{2}}\left(1-\frac{I(\tilde{X},\tilde{Y};J)+\log 2}{\log \mathcal{M}_{\mathcal{F}}^{\text{loc}}(\eta ,c)}\right)\ge \frac{\eta ^{2}}{4c ^{2}}\left(1-\frac{\frac{n \eta ^{2}}{2 \sigma ^{2}}+\log 2}{\log \mathcal{M}_{\mathcal{F}}^{\text{loc}}(\eta ,c)}\right)\ge \frac{\eta ^{2}}{8c ^{2}}.
    \end{equation*}
    The proof is completed.
\end{proof}

\section{Proof of Section \ref{section: minimax rate under heavy-tailed assumptions}}\label{section: proof of minimax rate under heavy-tailed assumptions}

\subsection{Proof of Lemma \ref{lemma: lower bound regarding epsilon under the bernstein's condition}}\label{subsection: proof of lower bound regarding epsilon under the bernstein's condition}

\begin{proof}
The conclusion is trivial if $\epsilon=0$, $\sigma=0$, or $d=0$. We therefore assume that these quantities are positive. Let $f_{\mathrm{c}}\in\mathcal{F}$ be a star center of $\mathcal{F}$. By the definition of the diameter and the triangle inequality, there exists $\bar{f}\in\mathcal{F}$ such that $\left\lVert \bar{f}-f _{c}\right\rVert_{L _{2}(\mathbb{P}_{X})}^{}>d/3$ (also see Lemma 1.3 of \cite{10.1214/25-AOS2576}). Set
\begin{equation*}
    r^2:=\epsilon\sigma^2\wedge\frac{d^2}{9},\qquad f_0:=f_{\mathrm{c}},\qquad f_1:=f_{\mathrm{c}}+\frac{r}{\left\lVert \bar{f}-f _{c}\right\rVert_{L _{2}(\mathbb{P}_{X})}^{}}\left(\bar{f}-f_{\mathrm{c}}\right).
\end{equation*}
Star-shapedness ensures that $f_0,f_1\in\mathcal{F}$, and $\left\lVert f _{1}-f _{0}\right\rVert_{L _{2}(\mathbb{P}_{X})}^{}=r$. Define
\begin{equation*}
    q(x):=\frac{\left[f _{1}(x)-f _{0}(x)\right]^{2}}{2 \sigma ^{2}+\left[f _{1}(x)-f _{0}(x)\right]^{2}}.
\end{equation*}
We construct two admissible regression models with regression functions $f_0$ and $f_1$. Let $X\sim\mathbb{P}_X$, draw $B\mid X\sim\operatorname{Bernoulli}(q(X))$, and let $R$ be a Rademacher random variable independent of $(X,B)$. For $j\in\{0,1\}$, define
\begin{equation*}
    Y^{(j)}:=
    \left\{
    \begin{tabular}{ll}
        $\frac{f _{0}(X)+f _{1}(X)}{2}+(2j-1)\dfrac{f _{1}(X)-f _{0}(X)}{2q(X)}B$, & $f _{0}(X)\neq f _{1}(X)$,\\
        $\frac{f _{0}(X)+f _{1}(X)}{2}+\frac{\sigma }{\sqrt{2}} R$, & $f _{0}(X)=f _{1}(X)$.
    \end{tabular}
    \right.
\end{equation*}
Denote the law of $(X,Y^{(j)})$ by $\mathbb{P} _{j}$ and set $\xi_j:=Y^{(j)}-f_j(X)$.

When $f _{0}(X)\neq f _{1}(X)$, we have
\begin{equation*}
    \mathbb{E}\left[Y^{(j)}\mid X\right]=\frac{f _{0}(X)+f _{1}(X)}{2}+(2j-1)\frac{f _{1}(X)-f _{0}(X)}{2}=f _{j}(X),
\end{equation*}
and
\begin{equation*}
    \mathbb{E}\left[\xi_j^2\mid X\right]=\frac{\left[f _{1}(X)-f _{0}(X)\right]^{2}}{4q(X)^2}q(X)\left(1-q(X)\right)=\frac{\sigma ^{2}}{2}.
\end{equation*}
Similar conditional mean and variance identities of $\xi _{j}$ hold on $\left\{f _{1}(X)=f _{0}(X)\right\}$ by the definition of $R$. Consequently, both models satisfy
\begin{equation*}
    \mathbb{E}\left[\xi_j\mid X\right]=0,
    \qquad
    \mathbb{E}\left[\xi_j^2\mid X\right]=\frac{\sigma ^{2}}{2}
    \quad\text{a.s.}
\end{equation*}
Thus both conditional noise laws are admissible. Moreover, the unconditional probability of the event $\{B=1\}$ satisfies
\begin{equation*}
    \bar q:=\mathbb{E}\left[q(X)\right]\le\frac{\left\lVert f _{1}-f _{0}\right\rVert_{L _{2}(\mathbb{P}_{X})}^{2}}{2 \sigma ^{2}}=\frac{r ^{2}}{2 \sigma ^{2}}\le \frac{\epsilon}{2}.
\end{equation*}
To construct the samples, take independent copies $(\tilde{X}_i,B_i,R_i)_{i=1}^n$ of $(X,B,R)$ and use the same copies to generate both
\begin{equation*}
    \tilde{\mathcal{D}}_j:=\left\{\left(\tilde{X}_i,\tilde{Y}_i^{(j)}\right)\right\}_{i=1}^n,\qquad j\in\{0,1\}.
\end{equation*}
Each marginal original sample has law $\mathbb{P}_j^{\otimes n}$.

Let $m:=\lfloor\epsilon n\rfloor$. Consider the following deterministic adversary $\mathcal{C}_0$. Given a sample $\mathcal{D}=\{(x_i,y_i)\}_{i=1}^n$, define
\begin{equation*}
    I(\mathcal{D}):=\left\{i:f _{0}(x _{i})\neq f _{1}(x _{i}),y_i\neq \frac{f _{0}(x _{i})+f _{1}(x _{i})}{2}\right\}.
\end{equation*}
If $\left|I(\mathcal{D})\right|\le m$, replace $y_i$ by $\frac{f _{0}(x _{i})+f _{1}(x _{i})}{2}$ for every $i\in I(\mathcal{D})$. Otherwise, leave the sample unchanged. This adversary changes at most $m$ observations for every input and is therefore admissible. Under either constructed model,
\begin{equation*}
    I\left(\tilde{\mathcal{D}}_j\right)=\left\{i:B_i=1\right\}\quad\text{a.s.}
\end{equation*}
Set $N:=\sum\limits _{i=1}^{n}B _{i}$ and $\mathcal{E}:=\{N\le m\}$. Since $\mathbb{E}N=n\bar q$, Markov's inequality gives
\begin{equation*}
    \mathbb{P}\left(\mathcal{E}^{\mathrm{c}}\right)=\mathbb{P}\left(N\ge m+1\right)\le\frac{n\bar q}{m+1}\le\frac{n\epsilon}{2(m+1)}<\frac12.
\end{equation*}
In particular, $\mathbb{P}(\mathcal{E})\ge1/2$, including when $\lfloor\epsilon n\rfloor=0$.

Define the corrupted samples by $\mathcal{D}_j:=\mathcal{C}_0(\tilde{\mathcal{D}}_j)$. On $\mathcal{E}$, all observations corresponding to $B_i=1$ are replaced by the common response $\frac{f _{0}(X _{i})+f _{1}(X _{i})}{2}$. All remaining observations already coincide under the coupling. Hence
\begin{equation*}
    \mathcal{D}_0=\mathcal{D}_1\qquad\text{on }\mathcal{E}.
\end{equation*}

Fix any estimator $\hat{f}$. If the estimator is randomized, couple its auxiliary randomness identically under the two models. On $\mathcal{E}$, its outputs therefore coincide. For every possible common output $\hat{f}$, we have 
\begin{equation*}
    \left\lVert \hat{f}-f _{0}\right\rVert_{L _{2}(\mathbb{P}_{X})}^{2}+\left\lVert \hat{f}-f _{1}\right\rVert_{L _{2}(\mathbb{P}_{X})}^{2}=2 \left\lVert \hat{f}-\frac{f _{0}+f _{1}}{2}\right\rVert_{L _{2}(\mathbb{P}_{X})}^{2}+\frac{1}{2}\left\lVert f _{1}-f _{0}\right\rVert_{L _{2}(\mathbb{P}_{X})}^{2}\ge \frac{r ^{2}}{2}.
\end{equation*}
Writing $\mathbb{E}_{\mathrm{cpl}}$ for expectation under the coupling, we obtain
\begin{equation*}
    \max_{j\in\{0,1\}}\mathbb{E}_{P_j^{\otimes n}}\left[\left\|\hat{f}\left(\mathcal{C}_0(\tilde{\mathcal{D}}_j)\right)-f_j\right\|_2^2\right]\ge\frac12\mathbb{E}_{\mathrm{cpl}}\left[\left(\|\hat{f}(\mathcal{D}_0)-f_0\|_2^2+\|\hat{f}(\mathcal{D}_1)-f_1\|_2^2\right)\mathbf{1}_{\mathcal{E}}\right]\ge \frac{r^2}{8}.
\end{equation*}
Both regression models and the adversary $\mathcal{C}_0$ are included in the supremum defining the minimax risk. Taking the infimum over $\hat{f}$ therefore gives
\begin{equation*}
    \mathcal{R}_{n,\epsilon}\left(\mathcal{F};\mathcal{P}_{\xi}\right)
    \ge \frac18\left(\epsilon\sigma^2\wedge\frac{d^2}{9}\right)\ge \frac{1}{72}\left(\epsilon\sigma^2\wedge d^2\right).
\end{equation*}
This proves the claim with an absolute constant independent of $d$, $n$, $\epsilon$, and $\sigma$.
\end{proof}

\subsection{Proof of Lemma \ref{lemma: trimmed mean on the current problem}}\label{subsection: proof of trimmed mean on the current problem}

\begin{proof}
    We compute the variance of $Z _{f,g}$ (under the uncontaminated samples, hence $\left(\left(\tilde{X}_{i},\tilde{Y}_{i},\tilde{\xi }_{i}\right)\right)_{i=1}^{n}=\left(\left(X _{i},Y _{i},\xi _{i}\right)\right)_{i=1}^{n}$). By rewriting $Z _{f,g}$, we have 
    \begin{equation*}
        \begin{aligned}
            Z _{f,g}&=\frac{\left[Y-f(X)\right]^{2}-\left[Y-g(X)\right]^{2}}{\left\lVert f-g\right\rVert_{L _{2}(\mathbb{P}_{X})}^{}}\\ 
            &=\frac{\left[g(X)-f(X)\right]\left[2Y-f(X)-g(X)\right]}{\left\lVert f-g\right\rVert_{L _{2}(\mathbb{P}_{X})}^{}}\\ 
            &=\frac{\left[g(X)-f(X)\right]}{\left\lVert f-g\right\rVert_{L _{2}(\mathbb{P}_{X})}^{}}\cdot \left[2f ^{*}(X)-f(X)-g(X)\right]+\frac{2\left[g(X)-f(X)\right]}{\left\lVert f-g\right\rVert_{L _{2}(\mathbb{P}_{X})}^{}}\cdot \xi .
        \end{aligned}
    \end{equation*}
    For the variance of $Z _{f,g}$, we have 
    \begin{equation*}
        \begin{aligned}
            \operatorname{Var}\left(Z _{f,g}\right)&\le \mathbb{E}Z _{f,g}^{2}\\ 
            &\overset{\text{(i)}}{=} \mathbb{E}\left[\frac{\left[g(X)-f(X)\right]^{2}}{\left\lVert f-g\right\rVert_{L _{2}(\mathbb{P}_{X})}^{2}}\cdot \left[2f ^{*}(X)-f(X)-g(X)\right]^{2}+\frac{4\left[g(X)-f(X)\right]^{2}}{\left\lVert f-g\right\rVert_{L _{2}(\mathbb{P}_{X})}^{2}}\cdot \xi ^{2}\right]\\ 
            &\overset{\text{(ii)}}{\le }\mathbb{E}\left[2 \cdot \frac{\left[g(X)-f(X)\right]^{2}}{\left\lVert f-g\right\rVert_{L _{2}(\mathbb{P}_{X})}^{2}}\cdot \left(\left[f ^{*}(X)-f(X)\right]^{2}+\left[f ^{*}(X)-g(X)\right]^{2}\right)+\frac{4\left[g(X)-f(X)\right]^{2}}{\left\lVert f-g\right\rVert_{L _{2}(\mathbb{P}_{X})}^{2}}\cdot \xi ^{2}\right]\\ 
            &\overset{\text{(iii)}}{\le } 2 \cdot \frac{\left\lVert f-g\right\rVert_{L _{4}(\mathbb{P}_{X})}^{2}}{\left\lVert f-g\right\rVert_{L _{2}(\mathbb{P}_{X})}^{2}}\left(\left\lVert f ^{*}-f\right\rVert_{L _{4}(\mathbb{P}_{X})}^{2}+\left\lVert f ^{*}-g\right\rVert_{L _{4}(\mathbb{P}_{X})}^{2}\right)+4 \sigma ^{2},
        \end{aligned}
    \end{equation*}
    where (i) follows from the fact that $\mathbb{E}\left[\xi \left\lvert\right. X\right]=0$; (ii) is from the inequality $(a+b)^{2}\le 2a ^{2}+2b ^{2},\forall a,b \in \mathbb{R}$, and (iii) is from $\mathbb{E}\left[\xi ^{2}\left\lvert\right. X\right]\le \sigma ^{2}$ and the Cauchy--Schwarz inequality.

    From the moment condition \eqref{eq: bernstein condition}, we have 
    \begin{equation*}
        \begin{aligned}
            &\left\lVert f-g\right\rVert_{L _{4}(\mathbb{P}_{X})}^{4}\le L \left\lVert f-g\right\rVert_{L _{2}(\mathbb{P}_{X})}^{2}\\ 
            &\left\lVert f ^{*}-g\right\rVert_{L _{4}(\mathbb{P}_{X})}^{4}\le L \left\lVert f ^{*}-g\right\rVert_{L _{2}(\mathbb{P}_{X})}^{2}\\ 
            &\left\lVert f ^{*}-f\right\rVert_{L _{4}(\mathbb{P}_{X})}^{4}\le L \left\lVert f ^{*}-f\right\rVert_{L _{2}(\mathbb{P}_{X})}^{2}.
        \end{aligned}
    \end{equation*}
    Therefore, we have
    \begin{equation*}
        \begin{aligned}
            \operatorname{Var}\left(Z _{f,g}\right)&\le 2L \cdot \frac{\left\lVert f ^{*}-f\right\rVert_{L _{2}(\mathbb{P}_{X})}^{}+\left\lVert f ^{*}-g\right\rVert_{L _{2}(\mathbb{P}_{X})}^{}}{\left\lVert f-g\right\rVert_{L _{2}(\mathbb{P}_{X})}^{}}+4 \sigma ^{2}\\ 
            &\overset{\text{(i)}}{\le } 2L \cdot \frac{\left\lVert f-g\right\rVert_{L _{2}(\mathbb{P}_{X})}^{}+2\left\lVert f ^{*}-g\right\rVert_{L _{2}(\mathbb{P}_{X})}^{}}{\left\lVert f-g\right\rVert_{L _{2(\mathbb{P}_{X})}}^{}}+4 \sigma ^{2},
        \end{aligned}
    \end{equation*}
    where (i) is from the triangle inequality. Note that we assume $\left\lVert f ^{*}-g\right\rVert_{L _{2}(\mathbb{P}_{X})}^{}\le \delta $ while $\left\lVert f-g\right\rVert_{L _{2}(\mathbb{P}_{X})}^{}\ge C\delta $. Therefore, we have 
    \begin{equation*}
        \operatorname{Var}\left(Z _{f,g}\right)\le 2L \left(1+\frac{2}{C}\right)+4 \sigma ^{2}.
    \end{equation*}
    By \eqref{eq: concentration inequality 2 of trimmed mean} in Lemma \ref{lemma: trimmed mean}, we have 
    \begin{equation*}
        \left|\operatorname{Tm}_{Z _{f,g}}(X,Y)-\mathbb{E}Z _{f,g}\right|\le 10 \sqrt{2L(1+\frac{2}{C})+4 \sigma ^{2}}\cdot \sqrt{8\epsilon +24 \frac{\log \left(\frac{4}{\rho }\right)}{n}}.
    \end{equation*}
    The proof is completed.
\end{proof}

\subsection{Proof of Corollary \ref{corollary: concrete application of the trimmed mean to bound the error probability}}\label{subsection: proof of concrete application of the trimmed mean to bound the error probability}

\begin{proof}
    We first verify the requirements of Lemma \ref{lemma: trimmed mean on the current problem}. The inequality trivially holds if $\rho =4 \exp\left\{-an\delta ^{2}/(L+\sigma ^{2})\right\}\ge 1$. Otherwise, we have 
    \begin{equation*}
        4 \exp\left\{-\frac{an\delta ^{2}}{L+\sigma ^{2}}\right\}\overset{\text{(i)}}{\ge }4 \exp\left\{-\frac{1}{2}n\right\}>\frac{1}{4}\exp\left\{-\frac{1}{2}n\right\},
    \end{equation*}
    where (i) is from the boundedness that $d ^{2}\le L$, $a \le \frac{1}{192}$, and that every nonvacuous comparison must have $\delta ^{2}\le d ^{2}\le L$. Also, the quantile $\zeta $ satisfies that 
    \begin{equation*}
        \zeta =8\epsilon +24 \frac{\log \left(\frac{4}{\rho }\right)}{n}=8\epsilon +24 \frac{a\delta ^{2}}{L+\sigma ^{2}}\le 8\epsilon +24a \le \frac{1}{4},
    \end{equation*}
    given that $\epsilon \le \frac{1}{64}$ and $a \le \frac{1}{192}$. Second, we calculate $\mathbb{E}Z _{f,g}$. For $f,g \in \mathcal{F}$ such that
    \begin{equation*}
        \left\lVert f-g\right\rVert_{L _{2}(\mathbb{P}_{X})}^{}\ge C\delta, \qquad \left\lVert f ^{*}-g\right\rVert_{L _{2(\mathbb{P}_{X})}}^{}\le \delta ,
    \end{equation*}
    the trimmed mean test \eqref{eq: pairwise comparison test} correctly selects $g$ if $\mathbb{E}Z _{f,g}$ dominates the variation term (RHS) in \eqref{eq: concentration inequality in the current problem}. By the \hyperref[subsection: proof of general minimax lower bound]{proof} of Lemma \ref{lemma: trimmed mean on the current problem}, it suffices that
    \begin{equation*}
        \mathbb{E}Z _{f,g}\ge 10 \sqrt{2L(1+\frac{2}{C})+4 \sigma ^{2}}\cdot \sqrt{8\epsilon +24 \frac{\log \left(\frac{4}{\rho }\right)}{n}}\asymp \sqrt{L+\sigma ^{2}}\left[\sqrt{\epsilon }+\sqrt{\log \left(\frac{4}{\rho }\right)/n}\right].
    \end{equation*}
    For $\mathbb{E}Z _{f,g}$, we have
    \begin{equation*}
        \begin{aligned}
            \mathbb{E}Z _{f,g}&=\mathbb{E}\left[\frac{\left[\xi +f ^{*}(X)-f(X)\right]^{2}-\left[\xi +f ^{*}(X)-g(X)\right]^{2}}{\left\lVert f-g\right\rVert_{L _{2}(\mathbb{P}_{X})}^{}}\right]\\ 
            &=\frac{\left(\left\lVert f ^{*}-f\right\rVert_{L _{2}(\mathbb{P}_{X})}^{}-\left\lVert f ^{*}-g\right\rVert_{L _{2}(\mathbb{P}_{X})}^{}\right)\left(\left\lVert f ^{*}-f\right\rVert_{L _{2}(\mathbb{P}_{X})}^{}+\left\lVert f ^{*}-g\right\rVert_{L _{2}(\mathbb{P}_{X})}^{}\right)}{\left\lVert f-g\right\rVert_{L _{2}(\mathbb{P}_{X})}^{}}\\ 
            &\overset{\text{(i)}}{\ge }\left\lVert f-g\right\rVert_{L _{2}(\mathbb{P}_{X})}^{}-2 \left\lVert f ^{*}-g\right\rVert_{L _{2}(\mathbb{P}_{X})}^{}\\ 
            &\ge (C-2)\delta ,
        \end{aligned}
    \end{equation*}
    where (i) is again from the triangle inequality. Therefore, to let the variation term be dominated by the true signal, we require
    \begin{equation*}
        \left\{
        \begin{tabular}{l}
            $(C-2)\delta \gtrsim \sqrt{L+\sigma ^{2}}\sqrt{\epsilon } \Leftrightarrow \delta ^{2}\gtrsim _{L,\sigma }\epsilon $ (the desired condition in Corollary \ref{corollary: concrete application of the trimmed mean to bound the error probability}),\\ 
            $(C-2)\delta \gtrsim \sqrt{L+\sigma ^{2}}\sqrt{\log (4/\rho )/n} \Leftrightarrow (C-2)\delta \gtrsim \sqrt{a}\delta $ (already achieved by our choice of $a$).
        \end{tabular}
        \right.
    \end{equation*}
\end{proof}

\subsection{Proof of Lemma \ref{lemma: crucial lemma for correct selection}}\label{subsection: proof of the crucial lemma for correct selection}

\begin{proof}
    The proof follows from an application of the union bound. Denote $f _{i ^{*}}=\mathop{\arg \min}\limits _{i \in [M]}\left\lVert f _{i}-f ^{*}\right\rVert_{L _{2}(\mathbb{P}_{X})}^{}$, then by the assumption of the lemma, we have $\left\lVert f _{i ^{*}}-f ^{*}\right\rVert_{L _{2}(\mathbb{P}_{X})}^{}\le \delta $. Consider any $f _{i}\in S$ such that $\left\lVert f _{i}-f ^{*}\right\rVert_{L _{2}(\mathbb{P}_{X})}^{}\ge (C+1)\delta $, then by the triangle inequality, we know $\left\lVert f _{i ^{*}}-f _{i}\right\rVert_{L _{2}(\mathbb{P}_{X})}^{}\ge C\delta $. By Corollary \ref{corollary: concrete application of the trimmed mean to bound the error probability}, we have 
    \begin{equation*}
        \mathbb{P}(f _{i ^{*}} \text{ fails to dominate } f _{i})\le 4 \exp\left\{-\frac{an\delta ^{2}}{L+\sigma ^{2}}\right\}.
    \end{equation*}
    Consequently, we have by the union bound 
    \begin{equation*}
        \mathbb{P}(f _{i ^{*}}\succ f _{i}, \forall f _{i}, \left\lVert f _{i ^{*}}-f _{i}\right\rVert_{L _{2}(\mathbb{P}_{X})}^{}\ge C\delta )\ge 1-4M \exp\left\{-\frac{an\delta ^{2}}{L+\sigma ^{2}}\right\}.
    \end{equation*}
    However, on the event $\left\{f _{i ^{*}}\succ f _{i}, \forall f _{i}, \left\lVert f _{i ^{*}}-f _{i}\right\rVert_{L _{2}(\mathbb{P}_{X})}^{}\ge C\delta \right\}$, we have $H (\delta ,f _{i},S)>0$ for all $f _{i}$ such that $\left\lVert f _{i}-f ^{*}\right\rVert_{L _{2}(\mathbb{P}_{X})}^{}\ge (C+1)\delta $ and $H (\delta ,f _{i ^{*}},S)=0$ (by definition of $H$ and the triangle inequality), which further implies that $\left\lVert f _{S}-f ^{*}\right\rVert_{L _{2}(\mathbb{P}_{X})}^{}<(C+1)\delta $. Therefore, we have 
    \begin{equation*}
        \mathbb{P}\left(\left\lVert f _{S}-f ^{*}\right\rVert_{L _{2}(\mathbb{P}_{X})}^{}\ge (C+1)\delta \right)=1-\mathbb{P}\left(\left\lVert f _{S}-f ^{*}\right\rVert_{L _{2}(\mathbb{P}_{X})}^{}< (C+1)\delta \right)\le 4M \exp\left\{-\frac{an\delta ^{2}}{L+\sigma ^{2}}\right\}.
    \end{equation*}
    The proof is completed.
\end{proof}

\subsection{Proof of Corollary \ref{corollary: probability of error during the trajectory of the tree traversal}}\label{subsection: proof of probability of error during the trajectory of the tree traversal}

\begin{proof}
    Let's directly evaluate $\mathbb{P}\left(A _{j-1}^{\complement }\cap A _{j}\right)$ for $j \ge 3$:
    \begin{equation*}
        \begin{aligned}
            \mathbb{P}\left(A _{j-1}^{\complement }\cap A _{j}\right)&=\mathbb{P}\left(\left\lVert f _{j-1}^{*}-f ^{*}\right\rVert_{L _{2}(\mathbb{P}_{X})}^{}\le \frac{d}{2 ^{j-2}},\left\lVert f _{j}^{*}-f ^{*}\right\rVert_{L _{2}(\mathbb{P}_{X})}^{}> \frac{d}{2 ^{j-1}}\right)\\ 
            &=\sum\limits _{f \in \mathcal{L}(j-1)\cap B _{2}(f ^{*},\frac{d}{2 ^{j-2}})}^{}\mathbb{P}\left(f _{j-1}^{*}    =f, \left\lVert f _{j}^{*}-f ^{*}\right\rVert_{L _{2}(\mathbb{P}_{X})}^{}> \frac{d}{2 ^{j-1}}\right)\\ 
            &\le \sum\limits _{f \in \mathcal{L}(j-1)\cap B _{2}(f ^{*},\frac{d}{2 ^{j-2}})}^{}\mathbb{P}\left(\left\lVert f _{j}^{*}-f ^{*}\right\rVert_{L _{2}(\mathbb{P}_{X})}^{}> \frac{d}{2 ^{j-1}}\right)\\ 
            &\overset{\text{(i)}}{\le }\sum\limits _{f \in \mathcal{L}(j-1)\cap B _{2}(f ^{*},\frac{d}{2 ^{j-2}})}^{}4 \mathcal{M}_{\mathcal{F}}^{\text{loc}}\left(\frac{d}{2 ^{j-2}},2c\right)\exp\left\{-n \eta _{j}^{2}\right\}\\ 
            &\overset{\text{(ii)}}{\le } 4 \left[\mathcal{M}_{\mathcal{F}}^{\text{loc}}\left(\frac{d}{2 ^{j-2}},2c\right)\right]^{2}\exp\left\{-n \eta  _{j}^{2}\right\}.
        \end{aligned}
    \end{equation*}
    The term $\mathcal{M}_{\mathcal{F}}^{\text{loc}}$ in (i) is from (4) of Lemma \ref{lemma: covering and packing properties of the constructed tree} and exponential term is from Lemma \ref{lemma: crucial lemma for correct selection} by observing $d/2 ^{j-1}=(C+1)\delta _{j}$. (ii) originates from Lemma \ref{lemma: bound on the cardinality}. When $i=2$, since $\mathcal{L}(1)$ only consists of a single point and $d>0$, we have 
    \begin{equation*}
        \mathbb{P}\left(A _{1}^{\complement }\cap A _{2}\right)\le 4 \mathcal{M}_{\mathcal{F}}^{\text{loc}}\left(d,c\right)\exp\left\{-n \eta _{2}^{2}\right\}\le  4 \left[\mathcal{M}_{\mathcal{F}}^{\text{loc}}\left(d,2c\right)\right]^{2}\exp\left\{-n \eta  _{2}^{2}\right\}.
    \end{equation*}
    For $i=1$, we trivially have by the definition of diameter:
    \begin{equation*}
        \mathbb{P}(A _{1})=\mathbb{P}(\left\lVert f _{1}^{*}-f ^{*}\right\rVert_{L _{2}(\mathbb{P}_{X})}^{}>d)=0.
    \end{equation*}
    The proof is completed.
\end{proof}

\subsection{Proof of Theorem \ref{theorem: main theorem of the algorithm}}\label{subsection: proof of main theorem of the algorithm}

\begin{proof}
    If $d=0$ or the sample size $n$ is overly small, then \eqref{eq: definition of J^*} might never hold, in which case $J ^{*}=1$ and $\eta _{J ^{*}}=\eta _{1}\asymp d$. Therefore, the upper bound trivially holds. Assume in the following that $d>0$ and \eqref{eq: definition of J^*} holds for some $J \ge 2$.

    For any $J \ge 2$ and a sequence of measurable sets $A _{1},\dots,A _{J}$, we have 
    \begin{equation*}
        \begin{aligned}
            A _{J}&=\left(A _{J}\cap A _{1}\right)\cup \left(A _{J}\cap A _{1}^{\complement }\right)\\ 
            &=\left(A _{J}\cap A _{1}\right)\cup \left(A _{J}\cap A _{1}^{\complement }\cap A _{2}\right)\cup \left(A _{J}\cap A _{1}^{\complement }\cap A _{2}^{\complement }\right)\\ 
            &=\dots\\ 
            &=\left[\bigcup\limits _{i=1}^{J}\left(A _{J}\cap \left(\bigcap\limits _{j=1}^{i-1}A _{j}^{\complement }\right)\cap A _{i}\right)\right]\cup \left(A _{J}\cap \left(\bigcap\limits _{j=1}^{J}A _{i}^{\complement }\right)\right)\\ 
            &=\bigcup\limits _{i=1}^{J}\left(A _{J}\cap \left(\bigcap\limits _{j=1}^{i-1}A _{j}^{\complement }\right)\cap A _{i}\right)\\ 
            &\subseteq \bigcup\limits _{i=1}^{J}\left(A _{i-1}^{\complement }\cap A _{i}\right),
        \end{aligned}
    \end{equation*}
    where we define $A _{0}=\varnothing $. Therefore, we have 
    \begin{equation}\label{eq: inequality between probabilities of sets}
        \mathbb{P}(A _{J})\le \sum\limits _{i=1}^{J}\mathbb{P}\left(A _{i-1}^{\complement }\cap A _{i}\right).
    \end{equation}
    Applying Corollary \ref{corollary: probability of error during the trajectory of the tree traversal} yields that for $J \le J ^{*}$:
    \begin{equation}\label{eq: control on the probability of A_J}
        \begin{aligned}
            \mathbb{P}\left(A _{J}\right)&\le \sum\limits _{i=1}^{J}\mathbb{P}\left(A _{i-1}^{\complement }\cap A _{i}\right)\\ 
            &\le \sum\limits _{j=2}^{J}4 \left[\mathcal{M}_{\mathcal{F}}^{\text{loc}}\left(\frac{d}{2 ^{j-2}},2c\right)\right]^{2}\exp\left\{-n\eta _{j}^{2}\right\}\\ 
            &\overset{\text{(i)}}{\le } 4 \left[\mathcal{M}_{\mathcal{F}}^{\text{loc}}\left(\frac{d}{2 ^{J-2}},2c\right)\right]^{2}\sum\limits _{j=2}^{J}\exp\left\{-n\eta _{j}^{2}\right\}\\ 
            &\le 4\left[\mathcal{M}_{\mathcal{F}}^{\text{loc}}\left(\frac{d}{2 ^{J-2}},2c\right)\right]^{2}\frac{\exp\left\{-n\eta _{J}^{2}\right\}}{1-\exp\left\{-n\eta _{J}^{2}\right\}}\\ 
            &\overset{\text{(ii)}}{\le } 4 \exp\left\{\frac{n\eta _{J}^{2}}{2}\right\}\cdot \frac{\exp\left\{-n\eta _{J}^{2}\right\}}{1-\exp\left\{-n\eta _{J}^{2}\right\}}\\ 
            &\overset{\text{(iii)}}{\le } 8 \exp\left\{-\frac{n\eta _{J}^{2}}{2}\right\},
        \end{aligned}
    \end{equation}
    where (i) is from Lemma \ref{lemma: nonincreasingness of the local entropy}; (ii) is from the fact that $n\eta _{J}^{2}>4 \log \left[\mathcal{M}_{\mathcal{F}}^{\text{loc}}\left(\frac{d}{2 ^{J-2}},2c\right)\right]$ and (iii) is from $n\eta _{J}^{2}>\log 2$, both originate from Condition \eqref{eq: definition of J^*} for $J \le J ^{*}$.

    From Algorithm \ref{algorithm: theoretical algorithm} and the argument above, we know that the closeness between $f _{J}^{*}$ and $f ^{*}$ (on the rate of $\eta _{J}$) depends on the correct selection for $f _{J-1}^{*}$, and hence on all previous trajectories $f _{1}^{*},\dots,f _{J-1}^{*}$. This, surely, is not guaranteed to happen. However, we can control the error probability for large-scale deviation given the control on $\mathbb{P}(A _{J})$ from \eqref{eq: control on the probability of A_J} and the property that the sequence $f _{1}^{*},\dots,f _{J ^{*}}^{*},\dots$ is a Cauchy sequence (Lemma \ref{lemma: bound on the distance between the parent and child}). Consider the event $\left\{\left\lVert f _{J ^{*}}^{*}-f ^{*}\right\rVert_{L _{2}(\mathbb{P}_{X})}^{}\ge \frac{d(1+3c)}{2 ^{J-1}c}\right\}$ for some $J \le J ^{*}$, we have 
    \begin{equation*}
        \begin{aligned}
            \mathbb{P}\left(\left\lVert f _{J ^{*}}^{*}-f ^{*}\right\rVert_{L _{2}(\mathbb{P}_{X})}^{}\ge \frac{d(1+3c)}{2 ^{J-1}c}\right)&\le \mathbb{P}\left(\left\lVert f _{J ^{*}}^{*}-f _{J}^{*}\right\rVert_{L _{2}(\mathbb{P}_{X})}^{}+\left\lVert f _{J}^{*}-f ^{*}\right\rVert_{L _{2}(\mathbb{P}_{X})}^{}\ge \frac{d(1+3c)}{2 ^{J-1}c}\right)\\ 
            &\overset{\text{(i)}}{\le }\mathbb{P}\left(\left\lVert f _{J}^{*}-f ^{*}\right\rVert_{L _{2}(\mathbb{P}_{X})}^{}+\frac{d(2+4c)}{2 ^{J}\cdot c}\ge \frac{d(1+3c)}{2 ^{J-1}c}\right)\\ 
            &\le 8 \exp\left\{-\frac{n\eta _{J}^{2}}{2}\right\},
        \end{aligned}
    \end{equation*}
    where (i) is from Lemma \ref{lemma: bound on the distance between the parent and child}. Note that we can extend this inequality to the region $J \le 0$, where the LHS is trivially zero followed by $\frac{d(1+3c)}{2 ^{J-1}c}\ge 6d, J \le 0$. Define $\kappa =\frac{1+3c}{2}\sqrt{\frac{L+\sigma ^{2}}{a}}$, then $\frac{d(1+3c)}{2 ^{J-1}c}=\kappa \eta _{J}$. For any $x \ge \eta _{J ^{*}}$, there exists an integer $J \le J ^{*}$ such that $\eta _{J}\le x \le \eta _{J-1}=2\eta _{J}$. Therefore, we have 
    \begin{equation}\label{eq: probability control on the large-scale deviation}
        \mathbb{P}\left(\left\lVert f _{J ^{*}}^{*}-f ^{*}\right\rVert_{L _{2}(\mathbb{P}_{X})}^{}\ge \kappa x\right)\le \mathbb{P}\left(\left\lVert f _{J ^{*}}^{*}-f ^{*}\right\rVert_{L _{2}(\mathbb{P}_{X})}^{}\ge \kappa \eta _{J}\right)\le 8 \exp\left\{-\frac{n\eta _{J}^{2}}{2}\right\}\le 8 \exp\left\{-\frac{nx ^{2}}{8}\right\}.
    \end{equation}
    Based on \eqref{eq: probability control on the large-scale deviation}, we can well bound the expectation of the squared loss.
    \begin{equation}\label{eq: bound on the squared loss}
        \begin{aligned}
            \mathbb{E}_{\tilde{X},\tilde{Y}}\left\lVert f _{J ^{*}}^{*}(\mathcal{C}(\tilde{X},\tilde{Y}))-f ^{*}\right\rVert_{L _{2}(\mathbb{P}_{X})}^{2}&=2\displaystyle\int _{0}^{\infty }x \cdot \mathbb{P}\left(\left\lVert f _{J ^{*}}^{*}-f ^{*}\right\rVert_{L _{2}(\mathbb{P}_{X})}^{}\ge x\right)dx\\ 
            &=2\kappa ^{2}\displaystyle\int _{0}^{\infty }t \cdot \mathbb{P}\left(\left\lVert f _{J ^{*}}^{*}-f ^{*}\right\rVert_{L _{2}(\mathbb{P}_{X})}^{}\ge \kappa t\right)dt\\ 
            &\le 2\kappa ^{2}\displaystyle\int _{0}^{\eta _{J ^{*}}}tdt+2\kappa ^{2}\displaystyle\int _{\eta _{J ^{*}}}^{\infty }t \cdot \mathbb{P}\left(\left\lVert f _{J ^{*}}^{*}-f ^{*}\right\rVert_{L _{2}(\mathbb{P}_{X})}^{}\ge \kappa t\right)dt\\ 
            &\le \kappa ^{2}\eta _{J ^{*}}^{2}+2\kappa ^{2}\displaystyle\int _{\eta _{J ^{*}}}^{\infty }8t \exp\left\{-\frac{nt ^{2}}{8}\right\}dt\\ 
            &=\kappa ^{2}\eta _{J ^{*}}^{2}+\frac{64\kappa ^{2}}{n}\exp\left\{-\frac{n\eta _{J ^{*}}^{2}}{8}\right\}\\ 
            &\overset{\text{(i)}}{\le }\kappa ^{2}\left(1+\frac{64}{\log 2}\right)\eta _{J ^{*}}^{2}\\ 
            &=\frac{L+\sigma ^{2}}{a}\left(\frac{1+3c}{2}\right)^{2}\left(1+\frac{64}{\log 2}\right)\eta _{J ^{*}}^{2},
        \end{aligned}
    \end{equation}
    where (i) is from the definition of $J ^{*}$ \eqref{eq: definition of J^*}. As an important note, we should be aware that the bound on the squared loss \eqref{eq: bound on the squared loss} depends on the condition that $\epsilon \lesssim \eta _{J ^{*}}^{2}$. When $\epsilon \gtrsim \eta _{J ^{*}}^{2}$, this crucial basis -- Lemma \ref{lemma: trimmed mean on the current problem}, Lemma \ref{lemma: crucial lemma for correct selection}, and Corollary \ref{corollary: probability of error during the trajectory of the tree traversal} -- of the probability control on $A _{J}$ can fail before $J$ increases to $J ^{*}$ and $\eta _{J}$ decreases to $\eta _{J ^{*}}$. In this case, we can alternatively define
    \begin{equation*}
        J_\epsilon:=\max\left(\{1\}\cup\left\{2\le j\le J^{*}:\epsilon\le\frac{b\delta_j^{2}}{L+\sigma^{2}}\right\}\right)
    \end{equation*}
    and parallel the argument above by replacing $J ^{*}$ with $J _{\epsilon }$ since in this case $\eta _{J _{\epsilon }}^{2}\ge \eta _{J ^{*}}^{2}$. Consequently, a counterpart inequality of \eqref{eq: bound on the squared loss} on the squared loss holds: 
    \begin{equation}\label{eq: bound on the squared loss 2}
        \mathbb{E}_{\tilde{X},\tilde{Y}}\left\lVert f _{J ^{*}}^{*}(\mathcal{C}(\tilde{X},\tilde{Y}))-f ^{*}\right\rVert_{L _{2}(\mathbb{P}_{X})}^{2}\le \frac{L+\sigma ^{2}}{a}\left(\frac{1+3c}{2}\right)^{2}\left(1+\frac{64}{\log 2}\right)\epsilon .
    \end{equation}
    Combining \eqref{eq: bound on the squared loss} and \eqref{eq: bound on the squared loss 2}, we know that $\eta _{J ^{*}}^{2}\vee \epsilon $ is an upper bound for the minimax rate \eqref{eq: minimax risk}. Since $d ^{2}$ is always a valid upper bound and the RHS of \eqref{eq: bound on the squared loss} and \eqref{eq: bound on the squared loss 2} do not depend on $\mathcal{C}$, $f ^{*}$, and $\xi $, we conclude that 
    \begin{equation*}
        \mathcal R_{n,\epsilon}(\mathcal F;\mathcal P_\xi^b)\lesssim \left(\eta _{J ^{*}}^{2}\vee \epsilon \right)\wedge d ^{2},
    \end{equation*}
    where a constant that depends only on $L,\sigma ,c$ is omitted. The proof is completed.
\end{proof}

\subsection{Proof of Corollary \ref{corollary: minimax rate under the bernstein condition}}\label{subsection: proof of minimax rate under the bernstein condition}

\begin{proof}
The lower bound follows from \eqref{eq: lower bound under the bernstein condition}, so it remains to prove the upper bound. The case $d=0$ is trivial. Suppose $d>0$, and write
\begin{equation*}
    \kappa:=\sqrt{\frac{a}{L+\sigma^{2}}},\qquad H_c(r):=\log\mathcal{M}_{\mathcal{F}}^{\mathrm{loc}}(r,c),\qquad \eta_j=\frac{\kappa d}{2^{j-2}c},\quad j\ge1.
\end{equation*}

First, star-shapedness and the definition of the diameter imply that $\mathcal{F}$ contains a line segment of length $d/3$. Consequently, $H_c(r)\ge\log 2$ for every $0<r\le d/3$. Applying the definition of $\eta^*$ at $r=(d/3)\wedge\sqrt{(\log 2)/n}$ gives
\begin{equation*}
    (\eta^*)^{2}\ge\frac{d^{2}}{9}\wedge\frac{\log 2}{n},\qquad \frac{1}{n}\wedge d^{2}\lesssim(\eta^*)^{2}\wedge d^{2}.
\end{equation*}

We next compare the packing constants $c$ and $2c$. Take any finite $r/(2c)$-packing of $B_2(f,r)\cap\mathcal{F}$ and partition it using a maximal $r/c$-packing of $B_2(f,r)\cap\mathcal{F}$. There are at most $\mathcal{M}_{\mathcal{F}}^{\mathrm{loc}}(r,c)$ groups by maximality, each contained in a ball of radius $r/c$, and each group has at most $\mathcal{M}_{\mathcal{F}}^{\mathrm{loc}}(r/c,2)$ points (recall these packing points are $r/2c$-separated, and $(r/c)/(r/2c)=2$). Taking logarithms and using $c\ge2$ and Lemma \ref{lemma: nonincreasingness of the local entropy}, we obtain
\begin{equation*}
    H_{2c}(r)\le H_c(r)+H_2(r/c)\le 2H_c(r/c).
\end{equation*}

Let $s:=\eta_{J^*+1}=\eta_{J^*}/2$. By the definition of $J^*$, Condition \eqref{eq: definition of J^*} fails at $J^*+1$, including when the convention $J^*=1$ applies. The preceding packing comparison therefore yields
\begin{equation*}
    ns^{2}\le 4H_{2c}\left(\frac{cs}{\kappa}\right)\vee\log 2
    \le 8H_c\left(\frac{s}{\kappa}\right)\vee\log 2.
\end{equation*}
Since $2\eta^*>\eta^*$, the definition of the critical radius gives $H_c(2\eta^*)<4n(\eta^*)^{2}$. If $s\ge2\kappa\eta^*$, monotonicity of $H_c$ consequently implies $ns^{2}\le32n(\eta^*)^{2}\vee\log 2$. Otherwise, $s^{2}<4\kappa^{2}(\eta^*)^{2}$. Thus, in either case,
\begin{equation*}
    \eta_{J^*}^{2}=4s^{2}\lesssim_{L,\sigma}(\eta^*)^{2}\vee\frac{1}{n}.
\end{equation*}
Combining this with the critical-radius lower bound established above gives
\begin{equation*}
    \eta_{J^*}^{2}\wedge d^{2}\lesssim_{L,\sigma}\left((\eta^*)^{2}\vee\frac{1}{n}\right)\wedge d^{2}\lesssim_{L,\sigma}(\eta^*)^{2}\wedge d^{2}.
\end{equation*}

Finally, Theorem \eqref{theorem: main theorem of the algorithm} and the identity $(u\vee v)\wedge w=(u\wedge w)\vee(v\wedge w)$ yield
\begin{equation*}
    \mathcal{R}_{n,\epsilon}(\mathcal{F};\mathcal{P}_{\xi})\lesssim_{L,\sigma}\left(\eta_{J^*}^{2}\vee\epsilon\right)\wedge d^{2}\lesssim_{L,\sigma}\left((\eta^*)^{2}\vee\epsilon\right)\wedge d^{2}.
\end{equation*}
Together with \eqref{eq: lower bound under the bernstein condition}, this proves the corollary. Here $a$ and $c$ are fixed universal constants, and none of the implicit constants depends on $d$.
\end{proof}

\section{Proof of Section \ref{section: minimax rate under sub-gaussian assumptions}}\label{section: proof of minimax rate under sub-gaussian assumptions}

\subsection{Proof of Lemma \ref{lemma: sub-gaussian lower bound on epsilon}}\label{subsection: proof of sub-gaussian lower bound on epsilon}

\begin{proof}
Throughout the proof, write $\left\lVert \cdot \right\rVert_{2}^{}=\left\lVert \cdot \right\rVert_{L_2\left(\mathbb{P}_X\right)}^{}$. The case $d=0$ or $\epsilon =0$ is immediate. Hence, suppose that $d>0$ and $\epsilon>0$. Since Gaussian noise is an admissible submodel, it suffices to establish the contamination-dependent lower bound under $\xi\sim \mathcal{N}\left(0,\sigma^2\right)$, independently of $X$.

Let $f_{\mathrm{c}}$ be a star center of $\mathcal{F}$. The triangle inequality gives
\begin{equation*}
    d\le 2\sup_{f\in\mathcal{F}}\left\lVert f-f_{\mathrm{c}} \right\rVert_{2}^{}.
\end{equation*}
Thus, there exists $u\in\mathcal{F}$ with $\left\lVert u-f_{\mathrm{c}} \right\rVert_{2}^{}\ge d/3$. By star-shapedness, the entire segment joining $f_{\mathrm{c}}$ and $u$ belongs to $\mathcal{F}$. Setting
\begin{equation*}
    D:=\min\left\{\frac{d}{3},\frac{\sigma\epsilon}{2}\right\},
\end{equation*}
we can therefore choose $f_0,f_1\in\mathcal{F}$ such that $\left\lVert f_1-f_0 \right\rVert_{2}^{}=D$.

For $j\in\left\{0,1\right\}$, let $P_j$ denote the joint distribution of one clean observation under
\begin{equation*}
    \tilde{X}\sim\mathbb{P}_X,\qquad \tilde{Y}=f_j\left(\tilde{X}\right)+\xi,\qquad \xi\sim \mathcal{N}\left(0,\sigma^2\right),\qquad \xi\perp\tilde{X}.
\end{equation*}
Because the design distribution is the same under both hypotheses, the conditional Gaussian likelihoods yield
\begin{equation*}
    \operatorname{KL}\left(P_0\,\middle\|\,P_1\right)=\frac{1}{2\sigma^2}\mathbb{E}_{\tilde{X}}\left[\left(f_0\left(\tilde{X}\right)-f_1\left(\tilde{X}\right)\right)^2\right]=\frac{D^2}{2\sigma^2}.
\end{equation*}
By Pinsker's inequality,
\begin{equation*}
    q:=\operatorname{TV}\left(P_0,P_1\right)\le\sqrt{\frac{1}{2}\operatorname{KL}\left(P_0\,\middle\|\,P_1\right)}=\frac{D}{2\sigma}\le\frac{\epsilon}{4}.
\end{equation*}

We next construct an admissible corruption procedure that makes the two experiments difficult to distinguish. By maximal coupling, there exists a pair $\left(U,V\right)$ with marginal distributions $P_0$ and $P_1$ such that
\begin{equation*}
    \mathbb{P}\left(U\ne V\right)=q.
\end{equation*}
Take $n$ independent copies $\left(U_i,V_i\right)$ of this coupling, and define
\begin{equation*}
    \tilde{\mathcal{D}}_0:=\left(U_1,\ldots,U_n\right),\qquad \tilde{\mathcal{D}}_1:=\left(V_1,\ldots,V_n\right).
\end{equation*}
Then $\tilde{\mathcal{D}}_j\sim P_j^{\otimes n}$, and the number of differing observations satisfies
\begin{equation*}
    W:=\sum_{i=1}^{n}\mathbf{1}\left\{U_i\ne V_i\right\}\sim\operatorname{Bin}\left(n,q\right).
\end{equation*}

Under $f_0$, let $\mathcal{C}_0$ be the identity corruption map. Under $f_1$, the adversary first samples a companion dataset $\tilde{\mathcal{D}}_0$ from its conditional distribution given the observed $\tilde{\mathcal{D}}_1$ under the coupling above. If $W\le\epsilon n$, it replaces precisely the differing observations and outputs $\tilde{\mathcal{D}}_0$. If $W>\epsilon n$, it leaves $\tilde{\mathcal{D}}_1$ unchanged. Denote this procedure by $\mathcal{C}_1$.

Since $W$ is integer-valued, $\mathcal{C}_1$ always replaces at most $\left\lfloor\epsilon n\right\rfloor$ observations. Hence, it is an admissible strong corruption procedure. Its use of auxiliary randomness is harmless for a minimax lower bound: the supremum over deterministic corruption maps is at least the risk averaged over these random seeds.

Let
\begin{equation*}
    Q_0:=P_0^{\otimes n},\qquad Q_1:=\mathcal{L}\left(\mathcal{C}_1\left(\tilde{\mathcal{D}}_1\right)\right).
\end{equation*}
Under the coupling just constructed, the output $\mathcal{C}_1\left(\tilde{\mathcal{D}}_1\right)$ agrees with $\tilde{\mathcal{D}}_0$ whenever $W\le\epsilon n$. Therefore, the coupling inequality and Markov's inequality give
\begin{equation*}
    \operatorname{TV}\left(Q_0,Q_1\right)\le\mathbb{P}\left(W>\epsilon n\right)\le\frac{\mathbb{E}\left[W\right]}{\epsilon n}=\frac{q}{\epsilon}\le\frac{1}{4}.
\end{equation*}
This comparison uses the unconditional law of $\tilde{\mathcal{D}}_0$; it does not require that conditioning on $W\le\epsilon n$ preserves that law.

For completeness, given any estimator $\widehat{f}$, construct a test by selecting whichever of $f_0$ and $f_1$ is closer to $\widehat{f}$ in $L_2\left(\mathbb{P}_X\right)$, with ties resolved arbitrarily. Whenever this test is incorrect under hypothesis $j$, the triangle inequality implies $\left\lVert \widehat{f}-f_j \right\rVert_{2}^{}\ge D/2$. The sum of the two testing error probabilities is at least $1-\operatorname{TV}\left(Q_0,Q_1\right)$. Consequently,
\begin{equation*}
    \max_{j\in\left\{0,1\right\}}\mathbb{E}_{Q_j}\left[\left\lVert \widehat{f}-f_j \right\rVert_{2}^{2}\right]\ge\frac{D^2}{8}\left[1-\operatorname{TV}\left(Q_0,Q_1\right)\right]\ge\frac{3D^2}{32}.
\end{equation*}
Taking the infimum over $\widehat{f}$, and recalling that these two functions and corruption procedures are admissible, yields
\begin{equation*}
    \mathcal R_{n,\epsilon}(\mathcal F;\mathcal P_\xi^{\mathrm{sg}})\ge\frac{3}{32}\min\left\{\frac{d^2}{9},\frac{\sigma^2\epsilon^2}{4}\right\}\ge\frac{1}{96}\left(\sigma^2\epsilon^2\wedge d^2\right).
\end{equation*}
The proof is completed.
\end{proof}

\subsection{Proof of Lemma \ref{lemma: sub-exponential of Z_f,g}}\label{subsection: proof of sub-exponential of Z_f,g}

\begin{proof}
    Recall that
    \begin{equation*}
        Z _{f,g}:=\frac{\left[Y-f(X)\right]^{2}-\left[Y-g(X)\right]^{2}}{\left\lVert f-g\right\rVert_{L _{2}(\mathbb{P}_{X})}^{}}=\frac{g(X)-f(X)}{\left\lVert f-g\right\rVert_{L _{2}(\mathbb{P}_{X})}^{}}\left[\left[2f ^{*}(X)-f(X)-g(X)\right]+2\xi \right].
    \end{equation*}
    First, by a direct application of \eqref{eq: sub-gaussian condition}, we know 
    \begin{equation*}
        \left\lVert \frac{g-f}{\left\lVert f-g\right\rVert_{L _{2}(\mathbb{P}_{X})}^{}}\right\rVert_{L _{2p}(\mathbb{P}_{X})}^{}\le L _{\text{sg}}\sqrt{2p}.
    \end{equation*}
    Similarly, we have 
    \begin{equation*}
        \left\lVert 2f ^{*}-f-g\right\rVert_{L _{2p}(\mathbb{P}_{X})}^{}\le L _{\text{sg}}\sqrt{2p}\left(\left\lVert f ^{*}-f\right\rVert_{L _{2}(\mathbb{P}_{X})}^{}+\left\lVert f ^{*}-g\right\rVert_{L _{2}(\mathbb{P}_{X})}^{}\right)\le L _{\text{sg}}\sqrt{2p}\left(\left\lVert f-g\right\rVert_{L _{2}(\mathbb{P}_{X})}^{}+2\delta \right),
    \end{equation*}
    where the latter step uses the assumption that $\left\lVert f ^{*}-g\right\rVert_{L _{2}(\mathbb{P}_{X})}^{}\le \delta $ and the triangle inequality.

    Combining the last two arguments with H\"{o}lder's inequality and the assumption that $\left\lVert f-g \right\rVert_{L _{2}(\mathbb{P}_{X})}^{}\ge C\delta $, we come to 
    \begin{equation*}
        \left\lVert \frac{g-f}{\left\lVert f-g\right\rVert_{L _{2}(\mathbb{P}_{X})}^{}}\left(2f ^{*}-f-g\right)\right\rVert_{L _{p}(\mathbb{P}_{X})}^{}\le 2L _{\text{sg}}^{2}(\left\lVert f-g\right\rVert_{L _{2}(\mathbb{P}_{X})}^{}+2\delta )p \le 4L _{\text{sg}}^{2}\left\lVert f-g\right\rVert_{L _{2}(\mathbb{P}_{X})}^{}p.
    \end{equation*}
    By a parallel argument, we have 
    \begin{equation*}
        \left\lVert 2 \frac{g-f}{\left\lVert f-g\right\rVert_{L _{2}(\mathbb{P}_{X})}^{}}\xi \right\rVert_{L _{p}(\mathbb{P}_{\xi ,X})}^{}\le 2 \cdot L _{\text{sg}}\sqrt{2p}\cdot L _{\xi }\sqrt{2p}\sigma =4L _{\text{sg}}L _{\xi }\sigma p.
    \end{equation*}
    Therefore, we have 
    \begin{equation*}
        \left\lVert Z _{f,g}-\mathbb{E}Z _{f,g}\right\rVert_{L _{p}(\mathbb{P}_{\xi ,X})}^{}\le 2\left\lVert Z _{f,g}\right\rVert_{L _{p}(\mathbb{P}_{\xi ,X})}^{}\le Kp(\left\lVert f-g\right\rVert_{L _{2}(\mathbb{P}_{X})}^{}+\sigma ),
    \end{equation*}
    where $K=8L _{\text{sg}}(L _{\text{sg}}+L _{\xi })$ depends only on $L _{\text{sg}}$ and $L _{\xi }$.
\end{proof}

\subsection{Proof of Lemma \ref{lemma: validity of the robust trimmed mean test under sub-gaussian condition}}\label{subsection: proof of validity of the robust trimmed mean test under sub-gaussian condition}

\begin{proof}
    We first evaluate $\mathcal{E}(\eta ,Z _{f,g})$ (recall its definition from Lemma \ref{lemma: trimmed mean}). By H\"{o}lder's inequality and Lemma \ref{lemma: sub-exponential of Z_f,g}, we have 
    \begin{equation*}
        \mathbb{E}\left[\left|\bar{Z} _{f,g}\right|\mathbf{1}_{\left\{\bar{Z}_{f,g}\ge Q _{1-\eta /2}(\bar{Z})\right\}}\right]\le \left(\mathbb{E}\left|\bar{Z}_{f,g}\right|^{p}\right)^{\frac{1}{p}}\left(\frac{\eta }{2}\right)^{1-1/p}\le K \left(\left\lVert f-g\right\rVert_{L _{2}(\mathbb{P}_{X})}^{}+\sigma \right)p \left(\frac{\eta }{2}\right)^{1-\frac{1}{p}}
    \end{equation*}
    Select $p=2 \vee \log \left(\frac{2e}{\eta }\right)$. We have: 
    \begin{equation*}
        \left(\frac{\eta }{2}\right)^{-\frac{1}{p}}=\left\{
        \begin{tabular}{ll}
            $\exp\left\{\frac{1}{1-1/\log \left(\eta /2\right)}\right\}$&, $(\eta \le \frac{2}{e})$,\\ 
            $\sqrt{\frac{2}{\eta }}$&, $(\frac{2}{e}<\eta <1)$
        \end{tabular}
        \right.
    \end{equation*}
    One can easily verify that $\sqrt{2}\le \left(\frac{\eta }{2}\right)^{-1/p}\le e$. Therefore, we have:
    \begin{equation*}
        \begin{aligned}
            \mathcal{E}(\eta ,Z _{f,g})&\le \frac{Kpe\eta }{2}\left(\left\lVert f-g\right\rVert_{L _{2}(\mathbb{P}_{X})}^{}+\sigma \right)\\ 
            &\le Ke\eta \log \left(\frac{2e}{\eta }\right)\left(\left\lVert f-g\right\rVert_{L _{2}(\mathbb{P}_{X})}^{}+\sigma \right).
        \end{aligned}
    \end{equation*}
    We have $\eta =8 \epsilon +24t _{\delta }\le 8 \epsilon +24a \le \frac{1}{8}$ by our selection of $\epsilon $ and $a$. Also, $4 \exp\left\{-nt _{\delta }\right\}\ge 4 \exp\left\{-an\right\}\ge \exp\left\{-n/2\right\}/4$. Applying \eqref{eq: concentration inequality 1 of trimmed mean} in Lemma \ref{lemma: trimmed mean} with $X$ replaced by $Z _{f,g}$ and Lemma \ref{lemma: sub-exponential of Z_f,g} with $p=2$, we know that with probability at least $1-4 \exp\left\{-nt _{\delta }\right\}$, we have 
    \begin{equation*}
        \begin{aligned}
            \left|\operatorname{Tm}(Z _{f,g})-\mathbb{E}Z_{f,g}\right|&\le 16Ke\eta \log \left(\frac{e}{2\eta }\right)\left(\left\lVert f-g\right\rVert_{L _{2}(\mathbb{P}_{X})}^{}+\sigma \right)+4 \sqrt{2}K \left(\left\lVert f-g\right\rVert_{L _{2}(\mathbb{P}_{X})}^{}+\sigma \right)\sqrt{t _{\delta }}\\ 
            &\le 16Ke\left(\left\lVert f-g\right\rVert_{L _{2}(\mathbb{P}_{X})}^{}+\sigma \right)\left(\eta \log \left(\frac{e}{\eta }\right)+\sqrt{t _{\delta }}\right)\\ 
            &\overset{\text{(i)}}{\le }16Ke\left(\left\lVert f-g\right\rVert_{L _{2}(\mathbb{P}_{X})}^{}+\sigma \right)\left(8\epsilon \log \left(\frac{e}{\epsilon }\right)+24t _{\delta }\log \left(\frac{e}{t _{\delta }}\right)+\sqrt{t _{\delta }}\right)\\ 
            &\overset{\text{(ii)}}{\le }128Ke \underbrace{\left(\left\lVert f-g\right\rVert_{L _{2}(\mathbb{P}_{X})}^{}+\sigma \right)\epsilon \log \left(\frac{e}{\epsilon }\right)}_{\text{contamination term}}+384Ke\underbrace{\left(\left\lVert f-g\right\rVert_{L _{2}(\mathbb{P}_{X})}^{}+\sigma \right)\sqrt{t _{\delta }}}_{\text{confidence term}},
        \end{aligned}
    \end{equation*}
    where in (i) we use $(a+b)\log \left(e/(a+b)\right)\le a \log (e/a)+b \log (e/b)$ and (ii) follows from the fact that $x \log (e/x)\le \sqrt{x}$ for small $x$ (equivalently, sufficiently small $a$).

    For the confidence term, observe that
    \begin{equation*}
        \begin{aligned}
            (\left\lVert f-g\right\rVert_{L _{2}(\mathbb{P}_{X})}^{}+\sigma)\sqrt{t_\delta}&=\sqrt a\,(\left\lVert f-g\right\rVert_{L _{2}(\mathbb{P}_{X})}^{}+\sigma)\min\left\{\frac\delta\sigma,1\right\}\\ 
            &\le\sqrt a\,(\left\lVert f-g\right\rVert_{L _{2}(\mathbb{P}_{X})}^{}+\delta)\\ 
            &\le\sqrt a\left(1+\frac1C\right)\left\lVert f-g\right\rVert_{L _{2}(\mathbb{P}_{X})}^{}.
        \end{aligned}
    \end{equation*}
    The middle inequality holds in both regimes: for \(\delta\le\sigma\), expand the product; for \(\delta>\sigma\), use \(\sigma<\delta\). The last inequality follows from $\left\lVert f-g\right\rVert_{L _{2}(\mathbb{P}_{X})}^{}\ge C\delta $.

    For the contamination term, suppose $\delta\ge \frac{3456Ke}{C}\sigma \epsilon \log \left(\frac{e}{\epsilon }\right)$. Then
    \begin{equation*}
        \begin{aligned}
            (\left\lVert f-g\right\rVert_{L _{2}(\mathbb{P}_{X})}^{}+\sigma)\epsilon \log \left(\frac{e}{\epsilon }\right)&\le \left\lVert f-g\right\rVert_{L _{2}(\mathbb{P}_{X})}^{}\epsilon \log \left(\frac{e}{\epsilon }\right)+\frac{C\delta }{3456Ke}\\ 
            &\le \left\lVert f-g\right\rVert_{L _{2}(\mathbb{P}_{X})}^{}\left(\epsilon \log \left(\frac{e}{\epsilon }\right)+\frac{1}{3456Ke}\right).
        \end{aligned}
    \end{equation*}
    Combining these estimates,
    \begin{equation*}
        |\operatorname{TM}(Z_{f,g})-\mathbb{E}Z_{f,g}|\le 384Ke\left\lVert f-g\right\rVert_{L _{2}(\mathbb{P}_{X})}^{}\left[\epsilon \log \left(\frac{e}{\epsilon }\right)+\frac{1}{3456Ke}+\sqrt a\left(1+\frac1C\right)\right].
    \end{equation*}
    Take \(\epsilon_0\) small enough that $384Ke\epsilon _{0}\log \left(\frac{e}{\epsilon _{0}}\right)\le 1/9$. The bracketed expression, including $K$, is then at most $1/3$. Recall that $\mathbb{E}Z _{f,g}\ge \left\lVert f-g\right\rVert_{L _{2}(\mathbb{P}_{X})}^{}-2\delta \ge \left(1-\frac{2}{C}\right)\left\lVert f-g\right\rVert_{L _{2}(\mathbb{P}_{X})}^{}\ge 1/3 \left\lVert f-g\right\rVert_{L _{2}(\mathbb{P}_{X})}^{}$, which implies that the clean signal dominates the variation. Consequently, we correctly select $g$ over $f$. The proof is completed.
\end{proof}

\end{document}